
\magnification=\magstep1
%
%
\newif\iftitlepage     \titlepagetrue
\newtoks\titlepagehead \titlepagehead={\hfil}
\newtoks\titlepagefoot \titlepagefoot={\hfil}

\newtoks\runningauthor \runningauthor={\hfil}
\newtoks\runningtitle  \runningtitle={\hfil}

\newtoks\evenpagehead  \newtoks\oddpagehead
\evenpagehead={\hfil\the\runningauthor\hfil}
\oddpagehead={\hfil\the\runningtitle\hfil}

\newtoks\evenpagefoot  \evenpagefoot={\hfil\tenrm\folio\hfil}
\newtoks\oddpagefoot   \oddpagefoot={\hfil\tenrm\folio\hfil}

\headline={\iftitlepage\the\titlepagehead
\else\ifodd\pageno\the\oddpagehead
\else\the\evenpagehead\fi\fi}

\footline={\iftitlepage\the\titlepagefoot
\global\titlepagefalse
\else\ifodd\pageno\the\oddpagefoot
\else\the\evenpagefoot\fi\fi}

\font\tenbb=msbm10  
\font\sevenbb=msbm7 
\font\fivebb=msbm5  
\newfam\bbfam      
\textfont\bbfam=\tenbb
\scriptfont\bbfam=\sevenbb
\scriptscriptfont\bbfam=\fivebb
\def\bb{\fam\bbfam}

\font\tenfrak=eufm10  
\font\sevenfrak=eufm7 
\font\fivefrak=eufm5  
\newfam\frakfam      
\textfont\frakfam=\tenfrak
\scriptfont\frakfam=\sevenfrak
\scriptscriptfont\frakfam=\fivefrak
\def\frak{\fam\frakfam}

\font\sc=cmcsc10 
%
%
\def\monthname{\ifcase\month\or January\or February\or
    March\or April\or May\or June\or July\or August\or
    September\or October\or November\or December\fi}
\def\today{\monthname\ \number\day, \number\year}

\def\widedotfill{\leaders\hbox to 15pt{\hfil.\hfil}\hfill}

\def\q{{\bf q}}       
\def\n{{\bf n}}       
\def\b{{\bf b}}       
\def\s{{\bf s}}       

\def\eps{\varepsilon}
\def\La{\Lambda}
\def\th{\Theta}

\def\tr{ {}^{t}}      
\def\diag{\mathop{\rm diag}}
\def\mod{\mathop{\rm mod}}
\def\rank{\mathop{\rm rank}}
\def\rad{\mathop{\rm rad}}
\def\eps{\varepsilon}
\def\La{\Lambda}
\def\th{\Theta}

\def\tr{ {}^{t}}      
\def\diag{\mathop{\rm diag}}
\def\mod{\mathop{\rm mod}}
\def\rank{\mathop{\rm rank}}
\def\rad{\mathop{\rm rad}}

\def\congr#1#2#3{{#1}\,\equiv\,{#2}\ (\mathop{\rm mod}{#3})}
\def\ncongr#1#2#3{{#1}\,\not\equiv\,{#2}\ (\mathop{\rm mod}{#3})}

\def\y#1{<\!#1\!>} 

\def\co#1{c_{p}^{-1}{#1}}
\def\xa#1{\chi_{#1}(p)}
 
\def\repr#1#2#3#4{{{r^*(\q_{#2},p^{#1}\q_{#3})}\over{e(\q_{#4})}}}
\def\Repr#1#2#3{E(\q_{#2})\backslash{R^*(\q_{#2},p^{#1}\q_{#3})}}
\def\reprn#1#2#3#4{{{r^*(\n_{#2},p^{#1}\n_{#3})}\over{e(\n_{#4})}}}
\def\Reprn#1#2#3{E(\n_{#2})\backslash{R^*(\n_{#2},p^{#1}\n_{#3})}}

\def\rep#1#2{{{r(\q_{#1},{#2})}\over{e(\q_{#1})}}}
\def\Rep#1#2{E(\q_{#1})\backslash{R(\q_{#1},{#2})}}
\def\repn#1#2{{{r(\n_{#1},{#2})}\over{e(\n_{#1})}}}

\def\rey#1#2{{\frak r}_{#1}({#2})}
\def\Rey#1#2{{\frak R}_{#1}({#2})}
\def\brey#1{\bar{\frak r}({#1})}

\def\eich#1#2{{\frak t}^*_{#1}(p^{#2})}
\def\Eich#1#2{{\frak T}^*_{#1}(p^{#2})}
\def\beich#1{\bar{\frak t}^*(p^{#1})}
\def\heich#1{\hat{\frak t}^*_{a}(#1)}

\def\Z{{\bb Z}}
\def\C{{\bb C}}
\def\Q{{\bb Q}}
\def\N{{\bb N}}
\def\F{{\bb F}}
\def\H{{\bb H}}
\def\A{{\bb A}}
\def\D{{\bb D}}
\def\M{{\bb M}}

\def\P{\hfill\break}

\def\litem#1{\item{\hbox to \parindent{\enspace#1\hfill}}}
%
%
\def\rharr#1#2{\smash{\mathop{\hbox to .5in{\rightarrowfill}}
    \limits^{\scriptstyle#1}_{\scriptstyle#2}}}
\def\lharr#1#2{\smash{\mathop{\hbox to .5in{\leftarrowfill}}
    \limits^{\scriptstyle#1}_{\scriptstyle#2}}}
\def\dvarr#1#2{\llap{$\scriptstyle #1$}\left\downarrow
    \vcenter to .5in{}\right.\rlap{$\scriptstyle #2$}}

\def\diagram#1{{\normallineskip=8pt
    \normalbaselineskip=0pt \matrix{#1}}}

\def\Lemma#1{\par\bigskip{\sc Lemma #1. }}
\def\Theorem#1{\par\bigskip{\sc Theorem #1. }}
\def\Corollary#1{\par\bigskip{\sc Corollary #1. }}
\def\Proposition#1{\par\bigskip{\sc Proposition #1. }}
\def\Remark#1{\par\smallskip{\sc Remark #1. }}
\def\Proof{\par\medskip{\sl Proof. }}
\def\qed{\vbox{\hrule\hbox{\vrule\kern3pt\vbox{\kern6pt}\kern3pt\vrule}
         \hrule}}
\def\flushqed{\unskip\nobreak\hfil\penalty50\hskip.75em\null\nobreak\hfil
              \qed {\parfillskip=0pt \finalhyphendemerits=0 \par\bigskip}}
%
%
\global\evenpagehead={\line{{\tenbf\folio}\hfil
  \sc\the\runningauthor\hfil}}
\global\oddpagehead={\line{\hfil\sc
  \the\runningtitle\hfil\tenbf\folio}}
\runningauthor={f. andrianov}
\runningtitle={clifford algebras and shimura's lift}
\global\evenpagefoot={\hfil}
\global\oddpagefoot={\hfil}

{\parindent=0pt\obeylines
The University of Chicago \hfill \today
Department of Mathematics
5734 S. University avenue
Chicago, Il 60637, USA.
\smallskip
fedandr@math.uchicago.edu}
%
%
\null
\bigskip\bigskip\bigskip
\bigskip\bigskip\bigskip
\centerline{\bf CLIFFORD ALGEBRAS AND}
\bigskip
\centerline{\bf SHIMURA'S LIFT FOR THETA-SERIES}
\bigskip\medskip                                    
\centerline{\sl FEDOR ANDRIANOV}
\bigskip\bigskip
%
%
\bigskip
{\parindent=2cm\narrower\noindent\sevenrm
ABSTRACT. 
Automorph class theory formalism is developed for the case of integral
nonsingular quadratic forms in an odd number of variables. As an application,
automorph class theory is used to construct a lifting of
similitudes of quadratic $\Z$-modules of arbitrary nondegenerate
ternary quadratic forms to morphisms between certain subrings of
associated Clifford algebras. The construction explains and generalizes
Shimura's correspondence in the case of theta-series of positive definite
ternary quadratic forms. Relations between associated zeta-functions are
considered.
\par}
\footnote{}{{\bf Key words}: Clifford algebras, Hecke operators,
Shimura's correspondence, quadratic forms, theta-series, zeta-functions.}
\bigskip\bigskip
\bigskip\bigskip
%
%
\centerline{\sl Contents}
\bigskip
Introduction \widedotfill 2\par
\S 1. Isotropic sums \widedotfill 5\par
\S 2. Action of Hecke operators on theta-series \widedotfill 23\par
\S 3. Matrix Hecke rings of orthogonal groups \widedotfill 28\par
\S 4. Shimura's lift for theta-series \widedotfill 31\par
\S 5. Clifford algebras and automorph class lift \widedotfill 34\par
\S 6. Factorization of automorph lifts \widedotfill 44\par
\S 7. Concluding remarks \widedotfill 54\par
References \widedotfill 57\par
\bigskip\bigskip
%
%
\bigskip\noindent
{\bf Introduction}

\smallskip
In his important paper [12] G. Shimura discovered a remarkable correspondence
between modular forms of half-integral weight and forms of integral weight 
-- the so-called `` Shimura lift '' (see [12, Main Theorem]). 
P. Ponomarev [10] determined the effect of Shimura's lifting on the 
theta-series $\th (z,\q )$ of certain positive definite ternary 
quadratic forms $\q$. Using purely arithmetical approach he expressed 
the lift of $\th (z,\q )$ as an explicit linear combination of 
theta-series associated to a certain system of quaternary quadratic forms
coming from a quaternion algebra (see [10], Theorem 1). Considering the
Fourier coefficients of corresponding theta-series one can interpret these
results as a link between the numbers of representations of integers by
positive definite quadratic forms in three and in four variables. But
according to the theory of automorph class rings (i.e. matrix Hecke rings
of orthogonal groups) developed by A. Andrianov in [1] and [2], the majority
of known multiplicative properties of the numbers of integral 
representations  by quadratic forms turns out to be merely a reflection of 
certain relations between the representations themselves. Following ideas of
A. Andrianov and P. Ponomarev one can try to define a direct algebraic
correspondence between the representations of integers by quadratic forms
in three and four variables which would underlie Shimura's lift in case of
theta-series. 
In the present paper we examine such a correspondence constructed by means
of Clifford algebras. Our main result is contained in Theorems 6.5 and 6.7
together with Corollary 6.8 which can be summarized as follows (see notations
in the end of this section):
\Proposition{$\!$}
{\sl
Let $\q$ be an integral nonsingular ternary quadratic form and let
$\n$ be the integral quadratic form defined (up to integral equivalence)
by the norm on the even subalgebra of the Clifford algebra $C(\q)\,$.
Take $\{\q_1,\dots ,\q_h\}$ and $\{\n_1,\dots ,\n_H\}$ to be
complete systems of representatives of different equivalence classes
of the similarity classes of $\q$ and of $\n$ respectively. Let $p$ be a
prime number coprime to $\det\q\,$. Denote by
$$
{\cal T}^*_{\q}(p^{\iota}) = \bigcup\limits_{i=1}^h E(\q_i)\backslash
R^*(\q_i,p^{\iota}\q)
\quad\hbox{and}\quad
{\cal T}^*_{\n}(p^{\iota}) = \bigcup\limits_{j=1}^H E(\n_j)\backslash
R^*(\n_j,p^{\iota}\n)
$$
the sets of classes of (primitive) automorphs of $\q$ and of $\n$
respectively (with multiplier $p^{\iota}$). Then one can define an injection
$\Psi:A\mapsto\Psi_A,\ {\cal T}^*_{\q}(p^2)\rightarrow{\cal T}^*_{\n}(p^2)\,$,
which can be naturally continued to an injection 
${\Upsilon}:A\mapsto{\Upsilon}_A,\ 
{\cal T}^*_{\q}(p^2)\rightarrow{\cal T}^*_{\n}(p)\,$,
in such a way that ${\Upsilon}_A$ divides $\Psi_A\,$ from the right.
Conversely, for any ${\cal M}\in{\cal T}^*_{\n}(p)$ there exists a unique
$A\in{\cal T}^*_{\q}(p^2)$ such that ${\cal M}$ divides $\Psi_A\,$
from the right. This construction turns the set $\cup_j E(\n_j)\backslash
R(\n_j,p\,\n)$ of classes of quaternary automorphs into a 2-fold covering
of the set $\cup_i E(\q_i)\backslash R^*(\q_i,p^2\q))\,$ of classes of 
ternary automorphs .}
\par
\bigskip
In view of Eichler's commutation relations (2.5), (3.3) and (3.4) the
construction corresponds to Shimura's lift (4.1) for generic theta-series
of ternary positive definite quadratic forms and generalizes the latter
in a purely arithmetical way to the case of arbitrary  indeterminate
(nonsingular) ternary quadratic forms:
\Corollary{}
{\sl With the notation and under the assumptions of the above Proposition, 
let in addition
$\q$ be positive definite, then Shimura's lift of the generic 
theta-series $\th_{\{\q\}}(z)$ and the generic theta-series
$\th_{\{\n\}}(z)$ of the norm $\n$ on the even Clifford subalgebra 
$C_0(E,\q)$ have the same eigenvalues $\,p+1\,$ for all Hecke operators 
$T(p)$ with $p{\not|}\det\q$.}
\bigskip
\noindent
(See also Theorem 6.3 and Corollary 6.4.)
\par
Thus defined {\it automorph class lift}
(5.28) rises numerous questions which still await solution. Some of the
possibilities, including those concerning relations between
associated zeta-functions, will be discussed in Concluding Remarks (Section 7).
\par
I am greatly indebted to Professor A. Andrianov without whose 
inspiration this paper would never be possible. I would like to
express my deep gratitude to Professor W. Baily for many stimulating
conversations and to Professor P. Ponomarev for his friendly support.
\par
%
%
\bigskip\noindent
{\bf Notations}
\smallskip\noindent
As usual, the letters $\Z, \Q\hbox{ and }\C$ denote the ring
 of rational integers,
the field of rational numbers and the field of complex numbers respectively.
We put $\H=\{z\in\C\  ;\ Im\ z >0\}$ to be the upper half plane.
$\A^{m}_{n}$ is the set of all $(m\times n)$-- matrices with entries in a
set $\A,\ \A^{m}=\A^{m}_{1},\ \A_{n}=\A^{1}_{n}$.
We let $\La^{m}=GL_{m}(\Z)$ denote the group of all integral invertible
matrices of oder $m$ and $1_{m}$ stand for the unit matrix of oder $m$.
If $M$ is a matrix, then $\tr M$ denotes its transposed, $\tilde M$
is its adjoint and $|M|$ stands for the absolute value of its determinant
(for square $M$).
We write
$$
L[M]=\tr MLM \ ,
$$
if the product of matrices on the right makes sense.
We consider a quadratic form
$$
\q (X)=\sum\limits_{1\le i\le j\le m} q_{ij}x_{i}x_{j}\ ,\quad
X=\pmatrix{x_1 \cr \vdots \cr x_m}\ ,
$$
in $m$ variables with rational integral coefficients $q_{ij}$.
The symmetric matrix
$$
Q=(q_{ij})+\tr(q_{ij})
$$
is called the {\it matrix} of the form $\q$ and its determinant 
$\det Q=\det\q$
is the {\it determinant} of $\q$. Quadratic form $\q$ is {\it nonsingular}
if $\det\q\ne0$\ . Clearly,
$$
2\q (X) = Q[X]\ .
$$
Note that $Q$ is an {\it even} matrix (i.e. an integral symmetric matrix
whose diagonal elements are even). The {\it level} of $\q$ (and of $Q$)
is defined to be the least positive integer $\ell$ such that $\ell Q^{-1}$ is 
even. It should be noted that $\ell\hbox{ and }\det\q$ have the same prime
divisors.
If $A\in\Q^{m}_{n}$ with some $m,n\in\N$, then $\q[A]$ denotes the quadratic
form 
$$
\q'(Y)={1\over 2}Q'[Y]=\q(AY)={1\over 2}Q[AY],\ 
Y=\pmatrix{y_1 \cr \vdots \cr y_n}\ ,
$$
resulting from $\q$ by the linear change of variables $X=AY$ .
The matrix $A$ is called a {\it representation} (over $\Q$ ) of the form
$\q'$ by the form $\q$. In case $A\in\Z^{m}_{n}$ , we call it an
{\it integral representation} of $\q'\hbox{ by }\q$ and denote by
$$
R(\q,\q')=\{A\in\Z^{m}_{n}\ ;\ \q[A]=\q'\}=\{A\in\Z^{m}_{n}\ ;\ Q[A]=Q'\}\ ,
$$
the set of all such integral representations. In particular, if $n=1$ and
$\q'=ay^2$ then the set $R(\q,\q')$ coincides with the set $R(\q,a)$ of
(integral) representations of the number $a\hbox{ by }\q$ :
$$
R(\q,ay^2)=R(\q,a)=\{X\in\Z^m\ ;\ \q(X)=a\}\ .
$$
In case $\q=\q'$ we get $E(\q)=R(\q,\q)\cap\La^m$ -- the 
{\it group of units} of the form $\q$. Subset of $R(\q,\q')$ consisting of
matrices whose entries are coprime is denoted by $R^{*}(\q,\q')$ and its
elements are called {\it primitive (integral) representations of $\q'$ by 
$\q$} . We set
$$\displaylines{
r(\q,\q')=|R(\q,\q')|\ ,\ r^{*}(\q,\q')=|R^{*}(\q,\q')|\ ,
\ r(\q,a)=|R(\q,a)|\ , \cr
e(\q)=|E(\q)| \cr 
}$$
to be corresponding cardinalities. Finally, the abbreviation $\gcd$ stands
for {\it greatest common divisor}.
\vfill\break
%
%
\bigskip\noindent
{\bf 1. Isotropic sums}
\smallskip\noindent
If $\q\hbox{ and }\q'$ are integral quadratic forms in the same number of 
variables with equal determinants, we say that $\q'$ is {\it similar} to 
$\q$ if $R(\q, a\q')$ is not empty for some $a$ prime to $\det\q$. The set of
all quadratic forms similar to $\q$ is called the {\it similarity class} of
$\q$. The Minkowski reduction theory of quadratic forms (see [3] for example)
shows that the similarity class of a nonsingular integral quadratic form in
$m$ variables is a finite union
$$
\bigcup\limits_{i=1}^{h}\{\q_{i}[U]\ ;\ U\in GL_m(\Z)\}\ ,
$$
of mutually disjoint classes of integrally equivalent quadratic forms. 
(One can compare the similarity class with the {\it genus} of $\q$.
Although in general the two sets are different, but both consist of quadratic
forms with arithmetical invariants equal to those of $\q$. Moreover, 
both sets are finite disjoint unions of classes modulo integral equivalence.)
\par
For similar quadratic forms $\q$ and $\q'$ the elements of the set 
$R(\q,a\q')$ are called {\it automorphs} (of $\q$ to $\q'$) with {\it
multiplier} $a$. 
\Lemma{1.1}
{\sl For any integral nonsingular quadratic forms $\q$ and $\q'$ and for
any nonzero number $a\in\Z$ the set $R(\q,a\q')$ of automorphs with 
multiplier $a$ is a finite union of left cosets modulo $E(\q)\,$.}
\Proof
Let $m$ denote the rank of $\q$ and $\q'$ and let $Q,Q'$ be their respective
matrices. If $M\in R(\q,a\q')$ then $Q[M]=aQ'\,$. Considering corresponding
determinants we have 
$|M|=\left(a^{m}\det\q\,\slash\det\q'\right)^{1\slash2}=|a|^{m\slash 2}\in\Z$
which means that $R(\q,a\q')\subset\Delta^m\left(|a|^{m\slash 2}\right)\,$,
where 
$$\Delta^m(d)=\{D\in\Z^m_m\ ;\ |\det D|=d\}
\eqno(1.1)
$$
for a positive integer $d\,$. We claim that $\Delta^m(d)$ is a finite union of
left cosets modulo $\Lambda^m\,$ (see Notations). To establish this we note
first that according to the theory of elementary divisors (see [3], Lemma 3.2.2)
each double coset $\Lambda^m D\Lambda^m\,,\,D\in\Delta^m(d)$ has a unique
representative of the form $\mathop{\rm ed}(D)=\diag(d_1,\dots,d_m)$ with
$d_i>0\,,\,d_i|d_{i+1}\,$ and so the cardinality
$|\Lambda^m\backslash\Delta^m(d)\slash\Lambda^m|$ is bounded above by 
number of factorizations $d=\prod_{i=1}^m d_i,\ d_i>0,\,d_i|d_{i+1}\,$, which
means that $\Delta^m(d)$ is a finite disjoint union of double cosets modulo
$\Lambda^m\,$. Next we note that each of these double cosets is a finite 
(disjoint) union of left cosets:
$$
\Lambda^m D\,\Lambda^m=\bigcup_{D'\in D\cdot
\left((D^{-1}\Lambda^m D\ \cap\ \Lambda^m)\backslash\Lambda^m\right)}
\Lambda^m D' \quad , 
$$
because $\Lambda^m DU=\Lambda^m D$ with $U\in\Lambda^m$
if and only if $U\in D^{-1}\Lambda^m D\cap\Lambda^m$ and the latter
subgroup has finite index in $\Lambda^m\,$. To justify the last claim
consider the principle congruence subgroup of level $d$ of $\Lambda^m$:
$$
\Lambda^m(d)=\{M\in\Lambda^m\ ;\ \congr{M}{1_m}{d}\}\quad .
$$ 
Clearly $\Lambda^m(d)$ has a finite index in $\Lambda^m$ as the kernel of
reduction modulo $d: \Lambda^m\rightarrow GL_m(\Z\slash d\Z)\,$.
On the other hand $\Lambda^m(d)\subset D^{-1}\Lambda^m D$ since 
entry-wise $D\,\Lambda^m(d)\tilde D\equiv\congr{D\tilde D}{0}{d}$ and so
$D\,\Lambda^m(d)D^{-1}=(\pm d)^{-1}D\,\Lambda^m(d)\tilde D\subset\Lambda^m\,$,
here $\tilde D$ is the matrix adjoint to $D\,$. Thus 
$D^{-1}\Lambda^m D\cap\Lambda^m$ contains subgroup $\Lambda^m(d)$ of finite
index in $\Lambda^m\,$, which means that $D^{-1}\Lambda^m D\cap\Lambda^m$
is itself of finite index in $\Lambda^m\,$.
Summing up, we see that $\Delta^m(d)$ is a finite (disjoint) union of double
cosets each of which is a finite (disjoint) union of left cosets modulo
$\Lambda^m$ and thus $\Delta^m(d)$ is indeed a finite union of left cosets
modulo $\Lambda^m\,$. 
\par
It remains to note that $R(\q,a\q')\subset\Delta^m\left(|a|^{m\slash2}\right)$
is therefore also a finite union of left cosets modulo $\Lambda^m\,$.
But each left coset $\Lambda^m M\cap R(\q,a\q')$ consists of a single left
coset $E(\q)M$ (if not empty), because if $UM=M'\in R(\q,a\q')$ with
$U\in\Lambda^m$ then $Q[U]=Q[M'M^{-1}]=aQ'[M^{-1}]=Q\,$, i.e. $U\in E(\q)\,$.
We conclude that $R(\q,a\q')$ is a finite union of left cosets modulo $E(\q)\,$.
\flushqed
The above property of automorphs will allow us to construct certain Hecke
algebras of orthogonal groups in section 3. Then the obvious inclusion
$$
R(\q,a\q')\cdot R(\q',b)\subset R(\q,ab)
$$
for any quadratic forms $\q,\q'$ and for any integers $a,b$ (the dot refers
to the usual matrix multiplication) will help us to define an action of these
Hecke algebras on representations of integers by quadratic forms. In order to
prepare for later quantitative investigation of that action, we first need
to address the question of inverse inclusion, i.e. we need to factor a given
element $K\in R(\q,ab)$ into a product of an automorph $M\in R(\q,a\q')$ and
a representation $L\in R(\q',b)\,$. We also want to find total number of
such factorizations. For our purposes it is enough to restrict attention
to the case when the automorph multiplier $a$ is a power of a fixed prime
$p$ and quadratic forms $\q,\q'$ are similar. Moreover, it appears that
consideration of just two automorph multipliers $a=p$ or $a=p^2$ is completely
sufficient for all situations (see [2]). The case of quadratic forms in an
even number of variables and automorph multiplier $p$ was treated in great
detail in [1]. Our first goal is to conduct a similar investigation for the
case of quadratic forms in and odd number of variables and automorph 
multiplier $p^2$ (it is easy to see that $p^2$ is the least power of an
automorph multiplier possible for quadratic forms in an odd number of
variables, see the proof of the Lemma 1.1). 
\par
Thus, given $K\in R(\q,p^2b)$ we seek to find total number of factorizations
$K=ML$ with $M\in R(\q,p^2\q'),\,L\in R(\q',b)\,$, where $\q,\q'$ are similar
integral quadratic forms in an odd number of variables $m=2k+1\geq 3$ and $p$
is a prime coprime to $\det\q\,$. Following general method of investigation of
questions of this type developed in [1], we will view a matrix 
$M\in R(\q,p^2\q')$ as a solution of quadratic congruence 
$$
\congr{\q[M]}{0}{p^2}\quad ,
$$
whose matrix of elementary divisors $D_p$ has a specific form.
\Lemma{1.2}
{\sl Let $\q,\q'$ be similar integral quadratic forms in an odd number of
variables $m=2k+1\geq 3$ and let $p$ be a prime coprime to $\det\q\,$. Then
$$
R(\q,p^2\q')=\cases{
R^*(\q,p^2\q') &, if $\q$ and $\q'$ are not integrally equivalent,\cr
R^*(\q,p^2\q')\cup pR(\q,\q') &, otherwise,\cr}
$$
where $R^*(\q,p^2\q')$ is the set of primitive automorphs, 
$R(\q,\q')\subset\Lambda^m\,$, and the union on the right hand side is 
disjoint. Furthermore, the matrix of elementary divisors of an automorph
from $pR(\q,\q')$ is equal to $p1_m$ and the matrix of elementary divisors
of a primitive automorph has the form
$$
D_p=D_p(d)=\pmatrix{
1_d &           &   \cr
    & p1_b &   \cr
    &     & p^2 1_d \cr}
\eqno(1.2)    
$$
for some $d,\, 1\leq d\leq k\,$, here $k=(m-1)\slash 2$ and $b=m-2d\,$.
}
\Proof
Let $M\in R(\q,p^2\q')$ be an automorph and let $\mu$ denote the greatest
common divisor of its entries. Since $\q[M]=p^2\q'$ then $|\det M|=|M|=p^m\,$,
which implies that $\mu|p$ and so either $\mu=1$ and $M\in R^*(\q,p^2\q')$
or $\mu=p$ and $p^{-1}M\in R(\q,\q')\,$. We also note that any elementary
divisor of $M$ is a power (at most $m^{\hbox{th}}$) of $p\,$.
Next, since $\q$ and $\q'$ are similar, and in particular have equal
determinants, then $R(\q,\q')\subset
\Lambda^m$ and therefore this set is empty if $\q$ and $\q'$ are not
integrally equivalent. It is also clear that any double coset from
$\Lambda^m pR(\q,\q')\Lambda^m\subset p\Lambda^m$ coincides with 
$\Lambda^m p1_m\Lambda^m\,$. Furthermore, if $M\in R^*(\q,p^2\q')\,$, then
$p^2 M^{-1}=(p^2\slash\det M)\cdot\tilde M=(\det\q)^{-1}\tilde{Q'}\,\tr M Q$
is an integral matrix since $\det M$ and $\det\q$ are coprime. Therefore
the matrix $D_p$ of elementary divisors of $M$ has the form
$\diag(1_{s_0},p1_{s_1},p^2 1_{s_2})$ with $s_0+s_1+s_2=m$ and 
$s_1+2s_2=m\,$, which leaves us with $s_o=s_2=d$ and $1\leq d\leq 
(m-1)\slash 2\,$.
\flushqed
Thus if $M\in R^*(\q,p^2\q')\,$, then $\congr{\q[M]}{0}{p^2}$ and 
$M\in\Lambda^m D_p\Lambda^m$ for a matrix $D_p$ of the form (1.2).
Next, since at the moment we do not want to distinguish a particular form
$\q'$ in its equivalence class $\{\q'[U]\ ;\ U\in\Lambda^m\}$, we will be
interested only in different cosets $M\Lambda^m\in\Lambda^m D_p \Lambda^m
\slash\Lambda^m$ (and thus only in cosets $R^*(\q,p^2\q')\slash E(\q')$ 
after a choice of particular form $\q'$). Finally, because our interest
is in factorizations of a particular column $K\in\Z^m$ such that
$\congr{\q(K)}{0}{p^2}\,$, we will consider only those $M$ which divide
from the left such a $K\,$. To this end we introduce {\it isotropic sums}
of the form 
$$
{\cal S}_{p^2}(\q, D_p, K)=
\sum_{\textstyle{{ M\in\Lambda^m D_p\Lambda^m\slash\Lambda^m}\atop
{\congr{\q[M]}{0}{p^2}\,,\, M|K}}} 1 \qquad .
\eqno(1.3)
$$
We will use geometric methods to compute the above sum. The following
notions and results will help us to reinterpret it in geometric terms.
\Lemma{1.3}{\sl Let $\delta$ be a positive integer and $p$ be a rational
prime. For a matrix $M\in\Z^m_n$ with $m>n$ the following statements are
equivalent:
{\parindent=0.7cm
\smallskip
\item{\ \ i)} There exists a primitive matrix $M'\in\Z^m_n$ such that
$\congr{M'}{M}{p^\delta}\,$.
\smallskip
\item{ ii)} Columns of $M$ are linearly independent modulo $p\,$.
\smallskip
\item{iii)} The matrix $M$ can be complemented to a matrix $(MM'')\in\Z^m_m$
such that $\congr{(MM'')}{U}{p^\delta}$ with $U\in\Lambda^m\,$.
\par}}
\Proof First of all, recall that an integer matrix is called {\it primitive}
if and only if the greatest common divisor of its principal minors is equal
to 1. Note that if $M'$ is primitive then not all of its principal minors
are zero and so its $\rank_{\Z}$ is equal to $n\,$. Then 
$\rank_{(\Z\slash p^\delta\Z)}M=\rank_{(\Z\slash p^\delta\Z)}M'=n\,$, i.e.
columns of $M$ are linearly independent modulo $p^\delta\,$. It remains to
note that for an arbitrary $\delta\,$, linear independence modulo $p^\delta$
is equivalent to linear independence modulo $p\,$. Thus {\sl i)} implies
{\sl ii)}. Conversely, suppose that columns of $M$ are linearly independent
modulo $p$ (and thus modulo $p^\delta$). Then 
$\rank_{(\Z\slash p^\delta\Z)}M=n$ and using elementary transformations over
$\Z\slash p^\delta\Z$ one can find $A\in\Lambda^m$ such that
$$
\congr{AM}{
\pmatrix{
a_1 & \dots & 0 \cr
\vdots & \ddots  & \vdots \cr
0   & \dots & a_n \cr
p^\delta & \dots  & p^\delta \cr
0 & \dots  & 0 \cr
\vdots & \ddots & \vdots \cr
0   & \dots  & 0  \cr
}=L
}{p^\delta}
$$
with $\ncongr{\prod a_i}{0}{p}\,$. Moreover, choosing appropriate 
representatives modulo $p^\delta$ one can assume that $a_i$'s are
pairwise coprime. Clearly such $L$ is primitive. Then $M'=A^{-1}L$ 
is congruent to $M$ modulo $p^\delta$ and also primitive since
multiplication by an invertible matrix preserves primitiveness.
(Indeed, using Binet-Cauchy formula, one can see that any principal minor
of $L$ has the form
$$
\left| (A\cdot M')^{1,\dots,n}_{\gamma_1,\dots,\gamma_n} \right|=
\left| (A)^{1,\dots,m}_{\gamma_1,\dots,\gamma_n}
\cdot(M')^{1,\dots,n}_{1,\dots,m} \right|=
\sum_{1\leq\alpha_1<\dots<\alpha_n\leq m} 
\left| (A)^{\alpha_1,\dots,\alpha_n}_{\gamma_1,\dots,\gamma_n} \right|
\cdot\left|(M')^{1,\dots,n}_{\alpha_1,\dots,\alpha_n} \right|
$$
i.e. it is an integral linear combination of minors of oder $n$ of $M'\,$.
Therefore $\gcd$ of principal minors of $M'$ divides $\gcd$ of those of
$L$ which is equal to 1 since $L$ is primitive.) This establishes equivalence
of {\sl i)} and {\sl ii)}.
\par
To prove equivalence of {\sl i)} and {\sl iii)} it is enough to show that
a matrix $U'\in\Z^m_n$ is primitive if and only if it can be complemented
to $U=(U',U'')\in\Lambda^m\,$. This was done in a much more general setting
(for matrices over arbitrary Dedekind rings) in Lemmas 2.3, 2.4 and 2.5 of
[1]. Here we will present another explanation of this fact for our particular
case of integral matrices. If $U'\in\Z^m_n$ can be complemented to
$U=(U',U'')\in\Lambda^m\,$, then there exists $U^{-1}=\tr(\tr A\tr B)$
with $A\in\Z^n_m, B\in\Z^{m-n}_m\,$. So we have:
$$
U^{-1}\cdot U=
\pmatrix{
AU' & AU'' \cr
BU' & BU'' \cr
} = 1_m\quad ,
$$
in particular $AU'=1_n\,$. Applying Binet-Cauchy formula one has:
$$
\sum_{1\leq\gamma_1<\dots<\gamma_n\leq m}
\left|(A)^{\gamma_1,\dots,\gamma_n}_{1,\dots,n}\right|\cdot
\left|(U')_{\gamma_1,\dots,\gamma_n}^{1,\dots,n}\right|
=\det(A\cdot U')=1 \quad ,
$$
which means that an integral linear combination of principal minors of
$U'$ is equal to 1, i.e. their $\gcd$ is 1 and $U'$ is primitive.
\par
Conversely, assume first that $n=1$ and $U'\in\Z^m$ is a primitive column.
Using induction on $m\geq 2\,$ we want to prove that $U'$ can be complemented
to an invertible matrix. In the base case $m=2$ one has $U'=\tr(u_1,u_2)$ with
$\gcd(u_1,u_2)=1\,$. Therefore $xu_1+yu_2=1$ for some $x,y\in\Z$ and then
$$
\pmatrix{ u_1 & -y \cr u_2 & x \cr }\in\Lambda^2
$$
as intended. If $m>2$ and $U'=\tr(u_1,\dots,u_m)$ is primitive, then either
$u_2=\dots=u_m=0$ and we can take 
$$
U=\pmatrix{
u_1 & \dots & 0 \cr
\vdots & 1_m & \cr
0 &  & \cr
}\in\Lambda^m
$$
since $u_1=\pm 1\,$, or some of $u_2,\dots,u_m$ are not zero and by induction
hypothesis we can complement the primitive column
$V'=\tr(u_2\slash d,\dots,u_m\slash d)\in\Z^{m-1}$ to an invertible matrix
$(V'V'')\in\Lambda^m\,$, here $d=\gcd(u_2,\dots,u_m)\,$. Note that
$\gcd(u_1,d)=1$ as $U'$ is primitive, so $xu_1+yd=1$ for some $x,y\in\Z$ and
then
$$
U=\pmatrix{
u_1 & 0 & y \cr
dV' & V'' & xV' \cr
}\in\Lambda^m
$$
since $\det U=  xu_1(-1)^{m-1} + yd(-1)^{m+1} =\pm 1\,$. Thus any primitive
column can be complemented to an invertible matrix. Next assume that $n>1\,,
\, U'\in\Z^m_n$ is primitive and $\tr(u_1,\dots,u_m)$ is the first column of
$U'\,$. For any selection $1\leq i_1<\dots<i_n\leq m\,,\,\gcd(u_{i_1},\dots,
u_{i_n})$ divides the minor $\left|(U')^{1,\dots,n}_{i_1,\dots,i_n}
\right|$ and so $\gcd(u_1,\dots,u_m)=1$ as $\gcd$ of minors of oder $n$ of
$U'$ is equal to 1. Therefore $\sum_{i=1}^m x_i u_i=1$ for some $x_i\in\Z\,$.
The column $\tr(x_1,\dots,x_m)$ is primitive, and as we already know it can 
be complemented to a matrix $X\in\Lambda^m\,$. Then $\tr XU'$ is primitive
with $(\tr XU')_{11}=1\,$. Applying elementary transformations to rows of
$\tr XU'$ we can find matrix $Y\in\Lambda^m\,$, such that
$$
YU'=\pmatrix{
1 & \dots & * \cr
\vdots & V' & \cr
0 &  & \cr
}\in\Z^m_n
$$
is again a primitive matrix whose minors of oder $n$ are either $0$
(in case the minor does not include the first row), or have the form
$\left|(V')^{1,\dots,n-1}_{i_1,\dots,i_{n-1}}\right|\,,\,1\leq i_1<\dots<
i_{n-1}\leq m\,$. The latter implies that $\gcd$ of principal minors of $V'$
is 1, i.e. $V'\in\Z^{m-1}_{n-1}$ is also primitive and we can recursively
apply the above procedure to it. Continuing in this fashion one can find a
matrix $Z\in\Lambda^m$ such that
$$
ZU'=\pmatrix{
1 & \dots  & * \cr
\vdots & \ddots & \vdots  \cr
0 & \dots & 1 \cr
\vdots  &  & \vdots  \cr
0 & \dots & 0 \cr
}\in\Z^m_n \quad ,
$$
then the matrix $W=\left( ZU'\ |\ {0\atop{1_{m-n}}}\right)$ is clearly
invertible and we can take $U=(U',U'')=Z^{-1}\cdot W$ to complement $U'$ to an
invertible matrix.
\flushqed
We will consider quadratic modules of the form
$$
V_{p^\delta}(1_m)=
\left((\Z\slash p^{\delta}\Z)^m\,,\,\q\mathop{\rm mod}p^{\delta}\right)
\eqno(1.4)
$$
obtained by reduction modulo $p^{\delta}$ of quadratic module 
$\left(\Z^m,\q\right)$ whose (scalar product) matrix with respect to the
standard basis of $\Z^m$ is $Q\,$, the matrix of $\q\,$. For a matrix
$M\in\Z^m_n$ we will denote by $V_{p^\delta}(M)$ the quadratic submodule
of 
$\left((\Z\slash p^{\delta}\Z)^m\,,\,\q\mathop{\rm mod} p^{\delta}\right)$
spanned by the columns of $M$ modulo $p^{\delta}\,$.
\Lemma{1.4}{\sl Let $D_p=D_p(d)$ be a matrix of the form (1.2) and let
$M\in\Lambda^m D_p\Lambda^m\,$. Then the map $M\mapsto V_{p^2}(M)$ defines
a bijection of the set of right cosets 
$M\Lambda^m\subset\Lambda^m D_p\Lambda^m\slash\Lambda^m$ to the set of all
subgroups of $\left(\Z\slash p^2\Z\right)^m$ whose invariants are
$(\underbrace{\mathstrut p^2,\dots,p^2}_{d}\,,
\,\underbrace{\mathstrut p,\dots,p}_{b})\,$,
here $b=m-2d\,$.}
\Proof
First of all we note that if $M'\in M\Lambda^m$ then clearly
$V_{p^2}(M')=V_{p^2}(M)$ as generators of any of these two groups are integral
linear combinations of generators of the other. Thus our map is defined
coset-wise. Next, let $M=UD_pV$ with $U,V\in\Lambda^m$ and let $U_i$'s denote
the columns of $U\,$. Then each $V_{p^2}(U_i),\,1\leq i\leq d$ is a cyclic
group of oder $p^2$ and each $V_{p^2}(pU_i),\,d\leq i\leq d+b$ is a cyclic
group of oder $p\,$. (To see this note that if $n$ is the oder of a column
$W\in\left(\Z\slash p^2\Z\right)^m$ then $\congr{nW}{0}{p^2}$ and so
$\congr{W}{0}{p^2\slash n}\,$. But the $\gcd$ of elements of any of $U_i$'s
is equal to $1$ since $U in\Lambda^m\,$. Thus the oder of $W=U_i$ is $p^2$ and
the oder of $W=pU_i$ is $p\,$.) Next we claim that
$$
V_{p^2}(M)=V_{p^2}(U_1)\oplus\dots\oplus 
V_{p^2}(U_d)\oplus V_{p^2}(pU_{d+1})\oplus\dots\oplus V_{p^2}(pU_{d+b})\quad ,
$$
i.e. $V_{p^2}(M)$ is a direct sum of $d$ cyclic subgroups of oder $p^2$ and 
$b$ cyclic subgroups of oder $p\,$. Indeed,
$$\openup2pt\displaylines{
\qquad V_{p^2}(M)=V_{p^2}(UD_p)=V_{p^2}(U_1,\dots U_d,pU_{d+1},\dots,pU_{d+b})=
\hfill\cr
\hfill
V_{p^2}(U_1)+\dots+V_{p^2}(U_d)+V_{p^2}(pU_{d+1}+\dots+V_{p^2}(U_{d+b})\quad\cr
}$$
and we just need  to show that the above sum is direct. To see this consider
the following: if
$$
\congr{\sum_{i=1}^d\alpha_i U_i+p\sum_{j=d+1}^{d+b}\alpha_j U_j}{0}{p^2}
$$
then $\congr{U^{-1}\cdot(U_1,\dots,U_{d+b})\cdot
\tr(\alpha_1,\dots,\alpha_d,p\alpha_{d+1},\dots,p\alpha_{d+b})}{0}{p^2}$
and so $p^2|\alpha_i$ for $1\leq i\leq d$ and $p|\alpha_j$ for
$d\leq j\leq d+b\,$, thus the zero of the group can be written only as
the sum of zeros of the subgroups, i.e. the sum is indeed direct. Therefore
the invariants of $V_{p^2}(M)$ are $(p^2,\dots,p^2,p,\dots,p)\,$.
\par
Conversely, assume that $G\subset\left(\Z\slash p^2\Z\right)^m$ is a subgroup
with this set of invariants, i.e.
$$
G=V_{p^2}(K_1)\oplus\dots\oplus V_{p^2}(K_d)\oplus
V_{p^2}(L_1)\oplus\dots\oplus V_{p^2}(L_b)\quad ,
$$
where the oder of each $K_i$ is $p^2$ and the oder of each $L_i$ is $p\,$.
Since $\congr{pL_i}{0}{p^2}$ then $p|L_i$ for $1\leq i\leq b\,$. The integral
columns $K_1,\dots,K_d,p^{-1}L_1,\dots,p^{-1}L_b$ are linearly independent
modulo $p$ (because the sum of the cyclic groups generated by $K_i$'s and
$L_j$'s is direct). Therefore, according to the Lemma 1.3 there exists a
primitive matrix $A\in\Z^m_{d+b}$ such that
$\congr{A}{K_1,\dots,K_d,p^{-1}L_1,\dots,p^{-1}L_b)}{p^2}$ which can be
complemented to an integral invertible matrix $A'\in\Lambda^m$ for which
$$
G=V_{p^2}(K_1,\dots,K_d,L_1,\dots,L_b)=
V_{p^2}(A_1,\dots,A_d,pA_{d+1},\dots A_{d+b})=V_{p^2}(A'D_p)\quad .
$$
Thus any $G\subset\left(\Z\slash p^2\Z\right)^m$ with invariants
$(p^2,\dots,p^2,p,\dots,p)$ has the form $V_{p^2}(A'D_p)$ for some
$A'\in\Lambda^m\,$.
\flushqed
The above Lemma allows us to replace summation over $M\in\Lambda^m
D_p\Lambda^m\slash\Lambda^m$ in (1.3) by summation over different
submodules $V_{p^2}(M)\subset\left(\Z\slash p^2\Z\right)^m$ with
$M\in\Lambda^m D_p\Lambda^m\,$. Furthermore, conditions
$\congr{\q[M]}{0}{p^2}\,,\,M|K$ in (1.3) mean exactly that $V_{p^2}(M)$
should be an isotropic submodule of quadratic module $V_{p^2}(1_m)=
\left((\Z\slash p^2\Z)^m,\,\q\mathop{\rm mod} p^2\right)$ containing a
fixed vector $K\,$. In order to find the total number of such submodules
(and thus to compute the isotropic sum (1.3)), we will first determine
the total number of bases of certain kind which generate modules of
type $V_{p^2}(M)\,,\, M\in\Lambda^m D_p\Lambda^m$ and then factor this
number by the number if different bases generating the same module.
We start with exhibiting the bases of $V_{p^2}(M)$ we are interested in.
We choose
$$
\{U\cdot D_p\ ;\ U\in\Lambda^m\slash(D_p\Lambda^m D_p^{-1}\cap\Lambda^m)\}
$$
as a complete system of representatives from classes $M\Lambda^m\in
\Lambda^m D_p\Lambda^m\slash\Lambda^m$ and let $U_i\,,\,1\leq i\leq m$
denote the columns of matrix a $U\,$. Then any isotropic module 
$$
V_{p^2}(M)=V_{p^2}(UD_p)=V_{p^2}(U_1,\dots,U_d,pU_{d+1},\dots,pU_{d+b}),
\hbox{ where }
b=m-2d\ .
$$
Using Lemma 1.3 we replace condition $U\in\Lambda^m$ with requirement that
$U_1,\dots,U_{d+b}$ should be linearly independent modulo $p\,$. Then the
condition $\q[M]\,\equiv\,\congr{\q[UD_p]}{0}{p^2}$ should be replaced with
$$
\cases{ \congr{\tr U_iQU_j}{0}{p^2}\,,& if $1\leq i\leq j\leq d\,$, \cr
\congr{\tr U_iQU_j}{0}{p}\,,& if $1\leq i\leq d$ and $d+1\leq j\leq d+b\,$.\cr
}
\eqno(1.5)
$$
Conversely, it is easy to see that any (ordered) collection $\{U_1,\dots,
U_{d+b}\}\subset (\Z\slash p^2\Z)^d\times(\Z\slash p\Z)^b$ of linearly
indepemdent modulo $p$ vectors satisfying (1.5) defines a basis
$$
U_1,\dots,U_d,pU_{d+1},\dots,pU_{d+b}
\eqno(1.6)
$$
of an isotropic module of the form $V_{p^2}(M),\,M\in\Lambda^m D_p\Lambda^m\,$.
From now on, when referring to a basis of the form (1.6), we will always
assume that vectors $\{U_1,\dots,U_{d+b}\}\subset(\Z\slash p^2\Z)^d\times
(\Z\slash p\Z)^b$ are linearly independent modulo $p$
and satisfy (1.5). Let us find the total number of such bases.
\Lemma{1.5}{\sl Under the above notations there exist
$$
p^{d(m+b-1)}\cdot\left|GL_b(\F_p)\right|\cdot
\prod_{s=0}^{d-1}(p^{2(k-s)}-1)\,,\hbox{ where }\ k=(m-1)\slash 2 \quad,
$$
different (ordered) bases of the form (1.6).}
\Proof We note first that $U_1,\dots,U_d$ should form an {\it isotropic system}
of vectors of the quadratic module $V_{p^2}(1_m)\,$, see (1.4), i.e. a system
of linear independent modulo $p$ vectors spanning an isotropic submodule. In
particular they should also form an isotropic system of $d$ vectors in 
$V_p(1_m)\,$. The number of isotropic systems in $V_p(1_m)$ consisting of $d$
vectors was computed in [3], Proposition A.2.14 and equals to
$$
p^{d(d-1)\slash 2}(p^{2k}-1)(p^{2(k-1)}-1)\cdot\dots\cdot(p^{2(k-d+1)}-1)\quad.
$$
It remains to compute in how many ways such a system modulo $p$ can be lifted 
to an isotropic system in $V_{p^2}(1_m)\,$. Note that any vector $L\in
(\Z\slash p\Z)^m$ splits into $p^m$ different vectors of the form $L+pX$ in
$(\Z\slash p\Z)^m\,$. Therefore, if $\{L_1,\dots,L_d\}$ is an isotropic system
in $V_p(1_m)\,$, then $\{L_1+pX_1,\dots,L_d+pX_d\}$ would form an isotropic
system in $V_{p^2}(1_m)$ if and only if
$$
\tr(L_i+pX_i)Q(L_j+pX_j)\,\equiv\,
\congr{\tr L_i QL_j+p(\tr L_iQX_j+\tr X_iQL_j)}{0}{p^2}
$$
for $1\leq i\leq j\leq d\,$.
(We do not worry about linear independence of the vectors because linear
independence modulo $p^{\delta}$ is equivalent for any $\delta$ to linear
independence modulo $p\,$, as we already noted in the proof of Lemma 1.3.)
The above congruences form a system of $d(d+1)\slash 2$ linear equations
(over $\F_p$) with respect to $md$ unknowns -- elements of columns $X_i\,$.
Since columns $L_i\,,\,1\leq i\leq d$ are linearly independent modulo $p\,$,
the rank of the matrix of this linear system as well as the rank of its
extended matrix is equal to $d(d+1)\slash 2$ and therefore the system has
$$
p^{dm-d(d+1)\slash 2}
\eqno(1.7)
$$
solutions in $\F_p\,$. (Note that this is exactly the number of different 
isotropic systems in $V_{p^2}(1_m)$ which lie above a fixed isotropic system
in $V_p(1_m)\,$.) We conclude that there exist
$$
p^{md-d(d+1)\slash 2}\cdot p^{d(d-1)\slash 2}\cdot
\prod_{s=0}^{d-1}(p^{2(k-s)}-1)=p^{d(m-1)}\cdot\prod_{s=0}^{d-1}(p^{2(k-s)}-1)
$$
different isotropic systems $\{U_1,\dots,U_d\}\subset V_{p^2}(1_m)\,$.
In how many ways such a system can be extended to an isotropic basis of
the form (1.6)? From (1.5) one can see that each $U_j\,,\,d\leq j\leq d+b$
should be orthogonal to the linear span of $U_1,\dots,U_d\,$, i.e. it should
satisfy the following system of linear equations
$$
\congr{
\pmatrix{\tr U_1 Q \cr \vdots \cr \tr U_d Q \cr}\cdot
\pmatrix{y_1 \cr \vdots \cr y_m \cr}
}{0}{p}
$$
over $\F_p\,$. The above system has $m-d=d+b$ linear independent (modulo $p$)
solutions. Besides, each of $U_i\,,\,1\leq i\leq d$ is its solution. Therefore
the remaining (ordered) set $U_i\,,\,d+1\leq i\leq d+b$ can be chosen in 
$$
p^{db}\cdot\left|GL_b(\F_p)\right|
\eqno(1.8)
$$
different ways. Multiplying (1.8) by the number of different isotropic systems
of $d$ vectors in $V_{p^2}(1_m)$ computed above, we deduce the result of the
Lemma.
\flushqed
Next we need to determine which of the bases in question generate the same
isotropic module. Assume that $\{U_1,\dots,pU_{d+b}\}$ and
$\{V_1,\dots,pV_{d+b}\}$ are two (ordered) bases of the form (1.6) such that
$V_{p^2}(U_1,\dots,pU_{d+b})=V_{p^2}(V_1,\dots,pV_{d+b})\,$. Then
$$
\congr{(V_1,\dots,V_d,pV_{d+1},\dots,pV_{d+b})}
{(U_1,\dots,U_d,pU_{d+1},\dots,pU_{d+b})\cdot\pmatrix{ A & B \cr C & D \cr}}
{p^2}\quad ,
$$
where $A\in\Z^d_d\,,\,B\in\Z^d_b\,,\,C\in\Z^b_d\,,\,D\in\Z^d_d\,$. Because
$\congr{(U_1,\dots,U_{d+b})\cdot\pmatrix{B \cr pD \cr}}{0}{p}\,$, then
necessarily $\congr{B}{0}{p}\,$. Since vectors
$$
\congr{(V_1,\dots,V_{d+b})}
{(U_1,\dots,U_{d+b})\cdot\pmatrix{A & p^{-1}B \cr pC & D \cr}}{p^2}
$$
are linearly independent modulo $p\,$, then the system of congruences
$$
\congr{
\pmatrix{A & p^{-1}B \cr pC & D \cr}\cdot
\pmatrix{\gamma_1 \cr\vdots\cr \gamma_{d+b}}
}{0}{p}
$$
should have only trivial solution, for which it is necessary and sufficient
that
$$
\pmatrix{ A & p^{-1}B \cr pC & D \cr }\in GL_{d+b}(\F_p)\quad,
$$
i.e. $A\in GL_d(\F_p)$ and $D\in GL_b(\F_p)\,$. Therefore two sets of vectors
$\{U_1,\dots,pU_{d+b}\}$ and $\{V_1,\dots,pV_{d+b}\}$  generate the same 
isotropic module only if they differ by a matrix from the group
$$\openup2pt\displaylines{
\qquad\Gamma=\biggl\{
\pmatrix{A & B \cr C & D \cr}\in GL_{d+b}(\Z\slash p^2\Z)\ ;\hfill\cr
\hfill
A\in GL_d(\Z\slash p^2\Z),\, D\in GL_b(\Z\slash p^2\Z),\, \congr{B}{0}{p}
\biggr\}\,.\qquad(1.9)\,\cr
}$$
Clearly the converse is also true. Thus $\Gamma$ acts transitively on the set
of bases of the form (1.6) of an isotropic module of type $V_{p^2}(UD_p)\,$.
Let $\Gamma'$ denote stabilizer of a basis 
$(U_1,\dots,U_d,pU_{d+1},\dots,pU_{d+b})$ under this action and let 
$\pmatrix{A & B \cr C & D \cr}\in\Gamma'\,$, i.e.
$$
\congr{(U_1,\dots,U_d,pU_{d+1},\dots,pU_{d+b})\cdot
\pmatrix{1_d-A & B \cr C & 1_b-D \cr}}{0}{p^2}\quad.
$$
Because $U_i$'s are linearly independent over $\F_p\,$, the last congruence
is equivalent to the set of conditions: $\congr{A}{1_d}{p^2}\,,\,
\congr{B}{0}{p^2}\,,\,\congr{C}{0}{p}\,,\,\congr{D}{1_b}{p}\,$. Combining
this with (1.9) one can see that the cardinality of $\Gamma$ modulo the
stabilizer $\Gamma'$ is given by
$$
\left|\Gamma\slash\Gamma'\right|=p^{d(d+2b)}\cdot
\left|GL_b(\F_p)\right|\cdot\left|GL_b(\F_p)\right|
$$
and this is the number of different bases of form (1.6) generating some
fixed isotropic module of the type $V_{p^2}(UD_p),\,U\in\Lambda^m\,$.
Using Lemma 1.5 we conclude that the number of different isotropic modules 
of the type we are interested in is equal to
$$
{\cal S}_{p^2}(\q, D_p, 0)=
p^{d(d-1)}\cdot\left|GL_d(\F_p)\right|^{-1}\cdot
\prod_{s=0}^{d-1}(p^{2(k-s)}-1)\quad,
\eqno(1.10)
$$
where $k=(m-1)\slash 2$ and $0\in\Z^m$ is the zero vector. We note also that
since the zero vector belongs to any submodule of $V_{p^2}(1_m)\,$, then
${\cal S}_{p^2}(\q, D_p, K)={\cal S}_{p^2}(\q, D_p, 0)$ in case
$\congr{K}{0}{p^2}\,$, see (1.3).
\par
In order to finish our computation of isotropic sums (1.3), we still need to
find the number of isotropic modules of the form $V_{p^2}(UD_p),\,U\in
\Lambda^m\,$, which contain a fixed nonzero isotropic vector 
$K\in(\Z\slash p^2\Z)^m\,$.
It will be more convenient to consider two separate cases: either
$\ncongr{K}{0}{p}\,$, or $\congr{K}{0}{p}$ but $\ncongr{K}{0}{p^2}\,$.
\Lemma{1.6}{\sl Under the above notations, assume that $\ncongr{K}{0}{p}\,$.
Then there exist
$$
{\cal S}_{p^2}(\q, D_p, K)=
(p^d-1)\cdot p^{(d-1)^2}\cdot\left|GL_d(\F_p)\right|^{-1}\cdot
\prod_{s=1}^{d-1}(p^{2(k-s)}-1)
$$
different isotropic modules of the form $V_{p^2}(UD_p),\,U\in\Lambda^m\,$,
which contain the vector $K\,$.}
\Proof If $K\in V_{p^2}(UD_p)=V_{p^2}(U_1,\dots,U_d,pU_{d+1},\dots,pU_{d+b})$
and $\ncongr{K}{0}{p}$ then this vector can be complemented to a basis
$\{V_1,\dots,K,\dots,V_d,pV_{d+1},\dots,pV_{d+b}\}$ of $V_{p^2}(UD_p)\,$.
As we already know, such a basis differs from $\{U_1,\dots,pU_{d+b}\}$ only
by a matrix from $\Gamma\,$, see (1.9). Moreover, if
$$
\congr{K}{\sum_{i=1}^d \gamma_i U_i + p\sum_{j=d+1}^{d+b} \gamma_j U_j}{p^2}
\quad,
\eqno(1.11)
$$
then among the first $d$ columns of a matrix from $\Gamma$ which transforms
$\{U_1,\dots,pU_{d+b}\}$ into $\{V_1,\dots,K,\dots,V_d,pV_{d+1},\dots,pV_{d+b}\}$
there should be a column of the form 
$$
\tr(\gamma_1,\dots,\gamma_d,\gamma'_{d+1},\dots,\gamma'_{d+b})\,\hbox{ with }
\congr{\gamma'_j}{\gamma_j}{p}\,\hbox{ for }d+1\leq j\leq d+b\quad.
$$
Since there exist exactly
$$
{{d}\over{p^d-1}}\cdot\left|GL_d(\F_p)\right|\cdot\left|GL_b(\F_p)\right|\cdot
p^{d^2+b^2+3db-d-b}
$$
different matrices from $\Gamma$ with this property, then dividing by the
oder of the stabilizer $|\Gamma'|$ we deduce that a fixed isotropic module 
$V_{p^2}(UD_p)$ contains 
$$
{{d}\over{p^d-1}}\cdot\left|GL_d(\F_p)\right|\cdot\left|GL_b(\F_p)\right|\cdot
p^{d^2+2db-d-b}
\eqno(1.12)
$$
different bases of the form 
$\{V_1,\dots,K,\dots,V_d,pV_{d+1},\dots,pV_{d+b}\}\,$.
Now we just need to find the total number of bases of the form (1.6) which
contain $K$ among their first $d$ vectors. We will proceed in a fashion
similar to the proof of Lemma 1.5: first we compute the number of isotropic
systems of $d$ vectors in $(\Z\slash p^2\Z)^m$ containing the column $K$ and
then multiply the result by the number of ways such an isotropic system can
be extended to an isotropic basis of the form we are interested in (the latter
was already computed in the course of proof of Lemma 1.5, see (1.8)).
\par
We already have noted above that any isotropic system modulo $p^2$ is also an
isotropic system modulo $p\,$. Therefore we start with computation of the 
number of isotropic systems of $d$ vectors in $(\Z\slash p\Z)^m$ containing
$K\mod p\,$. We use induction: assume that 
$X_1,\dots,X_{n-1}\in(\Z\slash p\Z)^m$ 
form an isotropic system. In how many ways it can be
complemented to an isotropic system of $n$ vectors? Since $X_1,\dots,X_{n-1}$
form a basis of the isotropic subspace $V_p(X_1,\dots,X_n)$ of nondegenerate
quadratic space $V_p(1_m)\,$ then there exist vectors $X'_1,\dots X'_{n-1}\in
V_p(1_m)$ such that each pair $X_i,X'_i\,,\,1\leq i\leq i-1,$ is hyperbolic,
i.e. vectors $X_i,X'_i$ are isotropic, linear independent and 
$\congr{\tr X_i Q X'_i}{1}{p}\,$, see [3], proposition A.2.12 or [9], 
Satz 2.24. Then the spaces $V_p(X_i,X'_i)\,,\,1\leq i\leq n-1$ are 
(nondegenerate) hyperbolic planes and
$$
V_p(1_m)=V_p(X_1,X'_1)\oplus V_p(X_{n-1},X'_{n-1})\oplus V
$$
is a direct sum of pairwise orthogonal subspaces. Therefore a vector
$Y\in V_p(1_m)$ can be isotropic and orthogonal to each of $X_i$'s
only if
$$
\congr{Y}{\sum_{i=1}^{n-1}\alpha_i X_i+v}{p}\quad,
$$
where $v$ is an isotropic vector of $V\,$. The vectors $X_1,\dots,X_{n-1},Y$
will be linearly independent only if $\ncongr{v}{0}{p}\,$. Since $V$ is
nondegenerate and $\dim V=m-2(n-1)\,$ then the number of nonzero isotropic
vectors in it equals to $p^{m-2n+1}-1\,$ (see [3], Proposition A.2.14).
Thus the vector $X_n$ complementing $\{X_1,\dots,X_{n-1}\}$ to an
isotropic system
of $n$ vectors can be chosen in $p^{n-1}(p^{m-2n+1}-1)$ ways. Using induction
on $n$ we conclude that a fixed isotropic vector $\ncongr{K}{0}{p}$ can be
complemented to
$$
p^{d(d-1)\slash2}\cdot\prod_{s=1}^{d-1}(p^{2(k-s)-1})
$$
different (ordered) isotropic systems $\{K,X_2,\dots,X_d\}\subset V_p(1_m)\,$.
Clearly, the total number of isotropic systems of $d$ vectors in $V_p(1_m)$
containing $K$ is $d$ times that. Next we need to lift our solution to
$V_{p^2}(1_m)$ keeping in mind that $K$ is fixed modulo $p^2\,$. Since any
vector $X\in(\Z\slash p\Z)^m$ splits in $(\Z\slash p^2\Z)^m$ into $p^m$ 
different vectors $X+pY\,,\,Y\in(\Z\slash p\Z)^m$, we need to look for 
isotropic systems of the form
$$
\{X_1+pY_1,\dots,K=X_i,\dots,X_d+pY_d\}\subset(\Z\slash p^2\Z)^m\quad.
$$
Such a system will be isotropic modulo $p^2$ if and only if
$$\openup2pt\displaylines{
\congr{\tr KQY_j}{-{1\over p}\tr X_j Q K}{p}\quad;\cr
\congr{\tr X_s Q Y_j+\tr X_j Q Y_s}{-{1\over p}\tr X_s Q X_j}{p}
\hbox{ for }s,j\not=i\quad.\cr
}$$
This is a system of $(d-1)(d+2)\slash2$ linear nonhomogeneous equations in
$m(d-1)$ variables (elements of columns $U_j,\,j\not=i$) over $\F_p\,$.
Since the columns $X_1,\dots,X_i=K,\dots,X_d$ are linearly independent modulo
$p\,$, so are the columns $QX_1,\dots,QX_d\,$, and therefore  the $\rank$ of
the matrix of the above system as well as the $\rank$ of its extended matrix
is equal to $(d-1)(d+2)\slash2\,$, which implies that it has
$p^{(d-1)(m-(d+2)\slash2)}$ solutions. We conclude that $V_{p^2}(1_m)$ contains
$$
dp^{(m-1)(d-1)}\cdot\prod_{s=1}^{d-1}(p^{2(k-s)})
$$
different isotropic systems of $d$ vectors one of which is $K\,$. According to
(1.8), each of these systems can be complemented to a basis of the form (1.6)
of some isotropic module $V_{p^2}(UD_p),\,U\in\Lambda^m$ in 
$p^{db}|GL_b(\F_p)|$ different ways. Multiplying the last two numbers we get
the total number of bases of the form (1.6) which contain $K$ among their
first $d$ vectors. Dividing the latter by (1.12) we finally deduce the result
of the Lemma.
\flushqed
It remains to consider the isotropic sum ${\cal S}_{p^2}(\q, D_p, K)$ in the
case when $\congr{K}{0}{p}$ but $\ncongr{K}{0}{p^2}\,$. Let 
$K\in V_{p^2}(UD_p)$ for some $U=(U_1,\dots,U_m)\in\Lambda^m$ and assume
(1.11). Then $\congr{\tr(\gamma_1,\dots,\gamma_d)}{0}{p}$ and we are faced with
two possibilities: 
$$
\matrix{
(i)\quad &\quad \ncongr{(\gamma_{d+1},\dots,\gamma_{d+b})}{0}{p}\,;\hfill \cr
(ii)\quad &\quad \congr{(\gamma_{d+1},\dots,\gamma_{d+b})}{0}{p}\,,\,
\ncongr{(\gamma_1,\dots,\gamma_d)}{0}{p^2}\,.\hfill \cr
}$$
In the first case the column $\tr(\gamma_1,\dots,\gamma_{d+b})$ can be 
complemented to a matrix from $\Gamma\,$, see (1.9), and in the second
case this is impossible. Note that the type of representation (1.11) of
the vector $K$ does not depend on the choice of basis of $V_{p^2}(UD_p)\,$.
Therefore, if $\congr{K}{0}{p}$ but $\ncongr{K}{0}{p^2}$ then there exists
two distinct types of isotropic modules which contain $K\,$: in modules
of the first type there exist a basis of the form 
$\{V_1,\dots,V_d,pV_{d+1},\dots,K,\dots,pV_{d+b}\}\,$, but in the modules
of the second type there are no such bases. In the two following Lemmas we
compute the number of isotropic modules of each these two types.
\Lemma{1.7}{\sl Under the above notations, assume that $\congr{K}{0}{p}$
but $\ncongr{K}{0}{p^2}\,$.  Then the number of isotropic modules of the
form $V_{p^2}(UD_p),\,U\in\Lambda^m$ for which (1.11) implies case $(i)$
above is equal to
$$\openup2pt\cases{
\eta_p\cdot p^{-d}(p^k-\xa{\q}\eps_p)(p^{k-d}+\xa{\q}\eps_p)
\prod_{s=1}^{d-1}(p^{2(k-s)}-1)\,,& 
if $\ \ncongr{\q(p^{-1}K)}{0}{p}\,$,\cr
\eta_p\cdot\prod_{s=1}^d(p^{2(k-s)}-1)\,,& 
if $\ \congr{\q(p^{-1}K)}{0}{p}\,$.\cr
}$$
Here $\eta_p=p^{d^2+b-1}(p^b-1)|GL_{b-1}(\F_p)||GL_d(\F_p)|^{-1}
|GL_b(\F_p)|^{-1}\,$, the character $\xa{\q}$ of the form $\q$ is defined
via (2.4) below and $\eps_p=\eps_p\left(\q(p^{-1}K)\right)$ is given by 
the Legendre symbol
$$
\eps_p(a)=\left({{2a}\over{p}}\right)\,,\ a\in\N\quad.
\eqno(1.13)
$$
}
\Proof As we already noted above, if an isotropic module satisfies conditions
of the Lemma, then it has a basis of the form 
$\{V_1,\dots,V_d,pV_{d+1},\dots,K,\dots,pV_{d+b}\}\,$. Thus we need to compute
the total number of bases of the form (1.6) which have $K$ among the last $b$
vectors, and then divide it by the number of such bases generating the same
isotropic module. The latter is equal to the number of matrices in $\Gamma$ 
(1.9) one of whose last $b$ columns is fixed modulo 
$(\Z\slash p^2\Z)^d\times(\Z\slash p\Z)^b$ divided by the oder of stabilizer
$|\Gamma'|\,$:
$$
b(p^b-1)^{-1}|GL_d(F_p)|\cdot|GL_B(\F_p)|p^{d^2+2db-d}\quad.
\eqno(1.14)
$$
Now we need to compute the total number of bases of the form (1.6) which have
$K$ among their last $b$ vectors. Assume first that
$\ncongr{\q(p^{-1}K)}{0}{p}\,$, i.e. the column $p^{-1}K$ is anisotropic. Then 
$V_p(1_m)=V_p(p^{-1}K)\oplus V\,$,
where $V$ is an orthogonal complement of $V_p(p^{-1}K)\,$ (see [9], Satz 1.15).
The module $V$ is nondegenerate and its dimension $\dim V=m-1$ is even. Any
vector orthogonal to $p^{-1}K$ belongs to $V$ and in particular is linear
independent (modulo $p$) with $p^{-1}K\,$. Therefore first $d$ columns of a
basis of the type in question should form an isotropic system in $V\,$.
According to Proposition A.2.14 of [3], the number of such systems is equal to
$$
p^{d(d-1)\slash2}(p^k-\xa{\q}\eps_p)(p^{k-d}+\xa{\q}\eps_p)
\prod_{s=1}^{d-1}(p^{2(k-s)}-1)\quad,
\eqno(1.15)
$$
where $\eps_p=\eps_p\left(\q(p^{-1}K)\right)$ is the Legendre symbol (1.13).
Each of these systems splits
into $p^{dm-d(d+1)\slash2}$ different isotropic systems in $V_{p^2}(1_m)\,$,
see (1.7), any of which can be chosen as first $d$ columns of a basis of 
the form (1.6). Now we need to find in how many ways an isotropic system
$U_1,\dots,U_d,K$ can be complemented by $b-1$ vectors from $(\F_p)^m$ to form
a basis of the type (1.6). Orthogonality conditions imply that elements of
each of the complementary columns should satisfy the following system of
linear homogeneous equations:
$$
\congr{
\pmatrix{\tr U_1 Q \cr \vdots \cr \tr U_d Q \cr}\cdot
\pmatrix{y_1 \cr \vdots \cr y_m \cr}
}{0}{p}
$$
This system of $d$ equations in $m$ variables over $\F_p$ has $d+b$ linear
independent modulo $p$ solutions which can be chosen in the form
$\vartheta=(U_1,\dots,U_d,p^{-1}K,V_2,\dots,V_b)\,$. Any other set of
linear independent solutions has the form $\vartheta\cdot M$ with
$M\in GL_{d+b}(\F_p)\,$. Among such matrices $M$ there are exactly
$b|GL_{b-1}(\F_p)|p^{(b-1)(d-1)}$ matrices which do not change the first
$d$ columns of the basis and leave $p^{-1}K$ intact (possibly changing its
position). And this is exactly the number of different (modulo $p$) sets of
vectors $pU_{d+1},\dots,K,\dots,pU_{d+b}$ which complement a given isotropic
system $U_1,\dots,U_d$ to a basis of the form (1.6). Multiplying this number
by (1.15) times (1.7) and dividing by (1.14) we deduce the statement of the
Lemma in the case when $\ncongr{\q(p^{-1}K)}{0}{p}\,$. 
\par
We turn to the case of an isotropic vector $p^{-1}K\in\Z^m\,$. Assume that
$\congr{\q(p^{-1}K)}{0}{p}$ and let $V\subset V_p(1_m)$  again denote its
orthogonal complement. We see that $\dim V=m-1$ again but now $V$ is degenerate:
it has nontrivial radical $\rad V= V_p(p^{-1}K)\,$. Nevertheless,
$V=\rad V\oplus W\,$, where $W$ is a nondegenerate subspace orthogonal to
$\rad V\,$. The first $d$ vectors of a basis of the type we are interested in
should form an isotropic system in $V\,$. Suppose that $U_1,\dots,U_d$ is a 
set of elements of $V\,$. Each of these vectors can be uniquely represented
in the form
$$
\congr{U_i}{\alpha_i p^{-1}K+w_i}{p}\,,\,\hbox{ where }\,
\alpha_i\in\F_p\,,\,w_i\in W\quad.
$$
Naturally, each of $U_i$'s is orthogonal to $p^{-1}K\,$. They will be 
isotropic an pairwise orthogonal if and only if
$\congr{\tr w_iQw_j}{0}{p}$ for $1\leq i\leq j\leq d\,$. Furthermore, the
columns $U_1,\dots,U_d,p^{-1}K$ will be linearly independent (modulo $p$) if
and only if $w_1,\dots,w_d\in W$ are linearly independent. We conclude that
$U_1,\dots,U_d$ will be an isotropic system if and only if $w_1,\dots,w_d$
form an isotropic system. The number of isotropic systems of $d$ vectors in
$W$ is equal to
$$
p^{d(d-1)\slash2}\cdot\prod_{s=1}^{d}(p^{2(k-s)}-1)
$$
(see [3], Proposition A.2.14). Multiplying this number by $p^d$ we get the
number of different isotropic systems in $V$ consisting of $d$ vectors linear
independent with $p^{-1}K\,$. The rest of the computations completely coincide
with the case of anisotropic $p^{-1}K$ and lead to the result stated in the
Lemma.
\flushqed
Let us turn to the second case.
\Lemma{1.8} {\sl Under the above notations, assume that $\congr{K}{0}{p}$ but
$\ncongr{K}{0}{p^2}\,$. Then the number of isotropic modules of the form 
$V_{p^2}(UD_p),\,U\in\Lambda^m$ for which (1.11) implies case $(ii)$ above is
equal to
$$\openup2pt\cases{
0\,,&
if $\ \ncongr{\q(p^{-1}K)}{0}{p}\,$,\cr
p^{(d-1)^2}|GL_{d-1}(\F_p)|^{-1}\cdot\prod_{s=1}^{d-1}(p^{2(k-s)}-1)\,,& 
if $\ \congr{\q(p^{-1}K)}{0}{p}\,$.\cr
}$$
}
\Proof Let $\{U_1,\dots,pU_{d+b}\}$ be a basis of the form (1.6). Assume that
(1.11) holds for some column $\tr(\gamma_1,\dots,\gamma_{d+b})\,$. Then the
conditions $\congr{(\gamma_1,\dots,\gamma_{d+b})}{0}{p}$ and 
$\ncongr{(\gamma_1,\dots,\gamma_{d})}{0}{p^2}$ will hold if and only if
$p^{-1}K\in V_p(U_1,\dots,U_d)\,$. But the subspace $V_p(U_1,\dots,U_d)
\subset V_p(1_m)$ is isotropic, therefore the column 
$\tr(\gamma_1,\dots,\gamma_{d+b})$ satisfies $(ii)$ if and only
$\congr{\q(p^{-1}K)}{0}{p}\,$, otherwise there exist no isotropic modules
which would satisfy the conditions of the Lemma. This proves its first case.
\par
Suppose now that $\congr{\q(p^{-1}K)}{0}{p}\,$. In order to find the number of
isotropic modules we are interested in, we need to divide the number  of 
different bases of the form (1.6) for which $p^{-1}K\in V_p(U_1,\dots,U_d)$
by the number of different bases of the form (1.6) which generate the same
isotropic module, i.e. by $|\Gamma\slash\Gamma'|$ computed above, see
discussion immediately following (1.9).  Next, the number of bases of the
form (1.6) with the property in question is equal to the product of the
number  of different isotropic systems $\{U_1,\dots,U_d\}\in V_{p^2}(1_m)$ for
which $p^{-1}K\in V_p(U_1,\dots,U_d)$ multiplied by the number of ways an
isotropic system of $d$ vectors can be complemented to a basis of the form
(1.6), i.e. by $p^{db}\cdot|GL_b(\F_p)|\,$, see (1.8). According to (1.7),
the number of isotropic systems $\{U_1,\dots,U_d\}\in V_{p^2}(1_m)$ for which
$p^{-1}K\in V_p(U_1,\dots,U_d)$ is equal in turn to the number of isotropic
systems $\{U_1,\dots,U_d\}\in V_p(1_m)$ for which $p^{-1}K\in 
V_p(U_1,\dots,U_d)$ times $p^{dm-d(d+1)\slash2}\,$. Thus it remains to find
the number of isotropic systems $\{U_1,\dots,U_d\}\in V_p(1_m)$ for which 
$p^{-1}K\in V_p(U_1,\dots,U_d)\,$. Since $\ncongr{K}{0}{p^2}$ then $p^{-1}K$
is a nonzero isotropic vector of $V_p(1_m)\,$, and so if $p^{-1}K\in
V_p(U_1,\dots,U_d)$ then $p^{-1}K$ can be complemented to a basis of 
$V_p(U_1,\dots,U_d)\,$.  Therefore the number of different isotropic systems
in $V_p(1_m)$ with the property in question is equal to the number of
isotropic systems of $d$ vectors in $V_p(1_m)$ containing the vector
$p^{-1}K$ (which was computed in the proof of Lemma 1.6) divided by
$dp^{d-1}\cdot|GL_{d-1}(\F_p)|\,$, i.e. by the number of such isotropic
systems generating the same module, and multiplied by $|GL_d(\F_p)|\,$, i.e.
by the number of different bases of a $d$-dimensional linear space over 
$\F_p\,$. This finishes the proof of the Lemma.
\flushqed
Combining the results of Lemmas 1.4 and 1.6--1.8 along with formulas (1.3)
and (1.10), we deduce the following Theorem:
\Theorem{1.9} {\sl Let $\q$ be an arbitrary integral nonsingular quadratic
form in an odd number of variables $m=2k+1\geq3\,$, and let $p$ be a rational
prime not dividing the determinant $\det\q$ of the form $\q\,$. Then for
each integral column $K\in\Z^m$ such that $\congr{\q[K]}{0}{p^2}\,$ the value of
the isotropic sum (1.3) is given by
$$\openup2pt\displaylines{
\qquad {\cal S}_{p^2}(\q, D_p(d), K)=\hfill\cr
\hfill\alpha_p(d)\cdot\left\{
1+\left(\eps_p \xa{\q}p^{k-1}+\kappa_p(d)\right)\cdot\delta_{(p^{-1}K)}+
p^{m-2}\cdot\delta_{(p^{-2}K)} \right\}\,,\quad(1.16)\,\cr
}$$
where $D_p(d)\,,\,1\leq d\leq k$ is a matrix of the form (1.2), $\eps_p=
\eps_p\left(\q(p^{-1}K)\right)$ is
given by the Legendre symbol (1.13), the character $\xa{\q}$ of the form 
$\q$ is defined via (2.4), $\delta_{(X)}$ is the generalized Kronecker symbol
(2.7) and
$$\openup2pt\displaylines{
\alpha_p(d)=p^{1-d}\cdot\prod_{s=1}^{d-1}(p^{2(k-s)}-1)\slash(1-p^{-s})\quad,\cr
\kappa_p(d)={{p^{m-2}-p^{d-1}}\over{p^d-1}}-1\quad.\cr
}$$
}
\Proof It is well-known that for $n\geq1$
$$
|GL_n(\F_p)|=\prod_{s=0}^{n-1}(p^n-p^s)=p^{n^2}\cdot\prod_{s=1}^{n}(1-p^{-s})
\quad.
$$
We also put $|GL_0(\F_p)|=1$ and note that $|GL_{n-1}(\F_p)|^{-1}\cdot
|GL_n(\F_p)|=p^{n-1}\cdot(p^n-1)\,$. Assume that $\ncongr{K}{0}{p}\,$, then 
according to Lemma 1.6
$$
{\cal S}_{p^2}(\q, D_p(d), K)=
p^{(d-1)^2-d^2}(p^d-1)\cdot(1-p^{-d})^{-1}\cdot
\prod_{s=1}^{d-1}(p^{2(k-s)}-1)\slash(1-p^{-s})=\alpha_p(d)\quad,
$$
which coincides with the right-hand side of (1.16) for such $K$ since
$\delta_{(p^{-1}K)}=\delta_{(p^{-2}K)}=0$ in this case.
Next, assume that $\congr{K}{0}{p}$ but $\ncongr{K}{0}{p^2}\,$. Combining
results of Lemmas 1.7, 1.8 and employing a rather straightforward calculation
one can see that ${\cal S}_{p^2}(\q, D_p(d), K)$ is equal to
$$
(p^d-1)^{-1}(p^k-\eps_p\xa{\q})(p^{k-d}+\eps_p\xa{\q})\cdot
\prod_{s=1}^{d-1}(p^{2(k-s)}-1)\slash(1-p^{-s})\quad,
$$
if $\ncongr{\q(p^{-1}K)}{0}{p}$ and
$$
(p^d-1)^{-1}(p^{a+b-1}-1)\cdot\prod_{s=1}^{d-1}(p^{2(k-s)}-1)\slash(1-p^{-s})
\quad,
$$
if $\congr{\q(p^{-1}K)}{0}{p}\,$. Note that in the latter case the symbol
$\eps_p=\eps_p\left(\q(p^{-1}K)\right)$ in (1.13) is equal to zero and
that otherwise $\eps_p=\pm1\,$. Therefore
$$
\openup2pt\displaylines{
\qquad {\cal S}_{p^2}(\q, D_p(d), K)=\hfill\cr
(p^d-1)^{-1}\cdot\left(p^{d+b-1}-1+\eps_p\xa{\q}p^k(1-p^{-d})\right)\cdot
\prod_{s=1}^{d-1}(p^{2(k-s)}-1)\slash(1-p^{-s})=\cr
\hfill\alpha_p(d)\cdot\left\{1+\eps_p\xa{\q}p^{k-1}+\kappa_p(d)\right\}
\qquad\cr
}$$
which coincides with the right-hand side of (1.16) for $\congr{K}{0}{p}\,$,
$\ncongr{K}{0}{p^2}\,$. Finally, if $\congr{K}{0}{p^2}$ then $\eps_p=0$ and,
using (1.10) we have
$$
{\cal S}_{p^2}(\q, D_p(d), K)=
\alpha_p(d)\cdot p^{d-1}\cdot(p^{2k}-1)\slash(p^d-1)\quad,
$$
which coincides with the right-hand side of (1.16) for such $K\,$.
\flushqed
Summing up the equalities (1.16) over $d$ ranging from $1$ to $k=(m-1)\slash2$
we immediately deuce the following
\Corollary{1.10} {\sl Under the notations and assumptions of Theorem 1.8,
$$
\sum_{d=1}^k{\cal S}_{p^2}(\q, D_p(d), K)=
c_p\cdot\left\{1+\left(\eps_p\xa{\q}p^{k-1}+c_p^{-1}\beta_p\right)
\cdot\delta_{(p^{-1}K)}+p^{m-2}\cdot\delta_{(p^{-2}K)}\right\}\quad,
$$
where $c_p=c_p(\q)$ and $\beta_p=\beta_p(\q)$ are given by (2.3).
}
\par\bigskip
We have finished computation of the isotropic sums of type (1.3). In order
to apply the results of Theorem 1.9 and Corollary 1.10 to investigation of
multiplicative arithmetic of integral representations by quadratic form
$\q\,$, we need to rewrite the isotropic sums in a slightly different terms.
Recall that the Minkowski reduction theory of quadratic forms (see [3] for
example) shows that the similarity class of a nonsingular integral quadratic 
form in $m$ variables is a finite union of mutually disjoint classes of 
integrally equivalent quadratic forms.  Let us fix a complete system
$$
\{\q_1,\dots,\q_h\}
$$
of representatives of different equivalence classes contained in the similarity
class of $\q\,$. Then the isotropic sum (1.3) can be rewritten as follows. The
condition $\congr{\q[M]}{0}{p^2}$ means that $\q[M]=p^2\q'\,$, where $\q'$ is
an integral quadratic form in $m$ variables. Since $\det M= \det D_p= p^m$
then $\det\q'=\det\q$ and thus quadratic form $\q'$ is similar to $\q\,$.
Since $M$ is defined only modulo right multiplication by $\Lambda^m\,$, we
can replace $\q'$ by any (integrally) equivalent form $\q'[U]\,,
\,U\in\Lambda^m\,$. In particular we can assume that $\q[M]=p^2\q_i$ for
exactly one of $\q_i$'s, $1\leq i\leq h$ and
$M\in\left(R(\q,p^2\q_i)\cap\Lambda^m D_p(d)\Lambda^m\right)\slash\Lambda^m$
for some $d\,$. Furthermore, if $M,MU\in R(\q,p^2\q_i)\cap
\Lambda^m D_p(d)\Lambda^m$ where $U\in\Lambda^m$ then $U\in R(\q_i,\q_i)=
E(\q_i)$ and so $M\in R(\q,p^2\q_i)\slash E(\q_i)\cap\Lambda^m 
D_p(d)\Lambda^m\,$. Therefore
$$
{\cal S}_{p^2}(\q, D_p(d), K)=
\sum_{i=1}^h\sum_{\textstyle{{ M\in\left(R(\q,p^2\q_i)\cap
\Lambda^m D_p(d)\Lambda^m\right)\slash E(\q_i)Lambda^m}\atop
{M|K}}} 1 \qquad.
\eqno(1.17)
$$
Next we note that according to the Lemma 1.2 the set of primitive automorphs
$R^*(\q,p^2\q_i)$ is a disjoint union
$$
R^*(\q,p^2\q_i)=\bigcup_{1\leq d\leq (m-1)\slash2}
\left(R(\q,p^2\q_i)\cap\lambda^m D_p(d)\Lambda^m\right)\quad.
$$
Therefore summing both sides of (1.17) over $d$ ranging from $1$ to 
$(m-1)\slash2$ and applying Corollary 1.10 we deduce the following Theorem:
\Theorem{1.11} {\sl Let $\q$ be an arbitrary integral nonsingular quadratic
form in an odd number of variables $m=2k+1\geq3\,$, and let $p$ be a rational
prime not dividing the determinant $\det\q$ of the form $\q\,$. Let 
$\{\q_1,\dots,\q_h\}$
be a complete system of representatives of different equivalence classes 
contained in the similarity class of $\q\,$. Then for each integral column 
$K\in\Z^m$ such that $\congr{\q[K]}{0}{p^2}\,$ we have
$$
\openup2pt\displaylines{
\qquad\sum_{i=1}^h\sum_{\textstyle{{ M\in R^*(\q,p^2\q_i)\slash E(\q_i)}\atop
{M|K}}} 1=\hfill\cr
\hfill c_p\cdot\left\{1+\left(\eps_p\xa{\q}p^{k-1}+c_p^{-1}\beta_p\right)
\cdot\delta_{(p^{-1}K)}+p^{m-2}\cdot\delta_{(p^{-2}K)}\right\}\,,\quad
(1.18)\,\cr
}$$
where the symbol $\eps_p=\eps_p\left(\q(p^{-1}K)\right)$
is defined in (1.13), $\xa{\q}$ is the character (2.4), $\delta_{(X)}$ is 
the generalized Kronecker symbol (2.7) and $c_p=c_p(\q),
\beta_p=\beta_p(\q)$
are given by (2.3).}
\vfill\break
%
%
\bigskip\noindent
{\bf 2. Action of Hecke operators on theta-series}
\smallskip\noindent
In case when $\q$ is positive definite, the numbers $r(\q,n)$ are finite 
and one can consider the theta-series
$$
\th(z,\q)=\sum\limits_{X\in\Z^m} e^{2\pi iz\q(X)}\ =\ 
\sum\limits_{n\geq0} r(\q,n) e^{2\pi inz}\ ,\ z\in\H
$$
and show that it is a (classical) modular form of weight $m/2$ and some 
character with respect to a certain congruence subgroup (see, for example
[3], Theorem 2.2.2). Although theta-series generally are not
eigenfunctions of Hecke operators, but the spaces spanned by theta-series
which come from a fixed similitude class of quadratic forms are invariant
under the action of Hecke operators. The fundamental discovery of M. Eichler
[6] is that in this situation the corresponding eigenmatrices (``Anzahlmatrizen''
in Eichler's terminology) are purely arithmetical and can be defined 
without any reference to modular forms. (The same phenomenon occurs in case
of Siegel modular forms of arbitrary degree, as was shown by E. Freitag [7]
and A. Andrianov [3].)

More specifically, let $\{\q_1,\dots,\q_h\}$ be a full system of
representatives of different equivalence classes of the similarity class of
$\q$ and let $p$ be a prime number which does not divide $\det\q$, then
$$
\th(z,\q)\bigm|_{m\over 2} T(p)=\co{}\cdot\sum_{1\leq j\leq h}
{\repr{}{}jj}\th(z,\q_{j})\ ,\eqno(2.1)
$$
if $m$ (the number of variables of $\q$) is even and
$$
\th(z,\q)\bigm|_{m\over 2} T(p^2)=\co{}\cdot\Bigl(\sum_{1\leq j\leq h}
{\repr2{}jj}\th(z,\q_{j})-\beta_p\th(z,\q)\Bigr)\ ,\eqno(2.2)
$$
if $m$ is odd. Here
$$
c_{p}=c_{p}(\q)=\cases{1 & if $m=2$, \cr
       \prod\limits_{i=0}^{k-2}(1+\xa{\q}p^i) & if $m=2k>2$, \cr
       \sum\limits_{a=1}^{k}\prod\limits_{i=1}^{a-1}
            {{p^{2(k-i)}-1}\over{1-p^{-i}}}\cdot p^{1-a} 
& if $m=2k+1\geq3$, \cr}
\eqno(2.3)
$$
$$
\beta_p=\beta_p(\q)=\sum_{a=1}^k\biggl(\prod\limits_{i=1}^{a-1}
{{p^{2(k-i)}-1}\over{1-p^{-i}}}\cdot 
{{p^{m-2}-p^{a-1}(p+1)+1}\over{p^{a-1}(p^a-1)}}\biggr)\ ,
$$
and the quadratic character $\xa{\q}$ is given by the Legendre symbol
$$
\xa{\q}=\Bigl({{(-1)^k\det\q}\over p}\Bigr)\ .
\eqno(2.4)
$$
For a detailed proof of formula (2.1) for the action of Hecke operators $T(p)$
on theta-series of integral weight see Chapter 5 of [3] and specifically
Theorem 5.2.5 there. The second formula (2.2) for the action of Hecke operators
$T(p^2)$ on theta-series of half-integral weight easily follows from combination
of Theorem 1.7 of [12] and Theorem 1.11 above:
\par\medskip{\sl Proof of Formula (2.2). }
Assume that $\q$ is a positive definite (integral) quadratic form in an odd
number of variables $m=2k+1\geq3\,$. Then using Theorem 1.7 of [12] we have
$$
\th(z,\q)\bigm|_{m\over 2} T(p^2)=\sum_{n\geq0}
\left\{r(\q,p^2n)+\xa{q}\left({{2n}\over{p}}\right)p^{k-1}r(\q,n)
+p^{m-2}r(\q, n\slash p^{2})\right\}\cdot e^{2\pi inz}\,,
$$
where we understand that $r(\q, n\slash p^2)=0$ if $p^2\not|\ n\,$.
On the other hand, summing up both sides of (1.18) over $K\in R^*(\q,p^2n)$
for some $n\in\N\cup\{0\}$ we conclude that the left-hand side
$$
\sum_{i=1}^h \sum_{K\in R(\q,p^2n)}
\sum_{\textstyle{{ M\in R^*(\q,p^2\q_i)\slash E(\q_i)}\atop{M|K}}} 1=
\sum_{i=1}^h \sum_{\textstyle{{ M\in R^*(\q,p^2\q_i)}\atop{L\in R(\q_i,n)}}}
{{1}\over{e(\q_i)}}=
\sum_{i=1}^h {{r^*(\q,p^2\q_i)}\over{e(\q_i)}}\cdot r(\q_i,n)
$$
is equal to the right-hand side
$$
\openup2pt\displaylines{
\quad\sum_{K\in R(\q,p^2n)}
c_p\cdot\left\{1+\left(\eps_p(n)\xa{\q}p^{k-1}+c_p^{-1}\beta_p\right)
\cdot\delta_{(p^{-1}K)}+p^{m-2}\cdot\delta_{(p^{-2}K)}\right\}=\hfill\cr
\hfill c_p\cdot\left\{r(\q,p^2n)+\left(\eps_p(n)\xa{\q}p^{k-1}+
c_p^{-1}\beta_p
\right)r(\q,n)+p^{m-2}r(\q,n\slash p^2)\right\}\,.
\quad\cr
}$$
Thus
$$
c^{-1}_p\cdot\left(
\sum_{i=1}^h {{r^*(\q,p^2\q_i)}\over{e(\q_i)}}r(\q_i,n)-\beta_p r(\q,n)\right)
$$
coincides with the $n^{th}$ Fourier coefficient of
$\Theta(z,\q)|_{m\over2}T(p^2)\,$, 
and on the other hand it is equal to the $n^{th}$ Fourier coefficient of
the series on the right-hand side of (2.2). This finishes the proof of (2.2)
because of the uniqueness of Fourier decomposition.
\flushqed
One can easily see that in case of odd number of variables the formulas are
not as nice as when $m$ is even. (Unfortunately this is always the case.)
Since we can substitute for $\q$ any of $\q_{i}\ ,\ 1\leq i\leq h$,
the above formulas lead to the so-called Eichler's commutation relation:
$$
\pmatrix{\vdots\cr{{\th(z,\q_i)}\over{e(\q_i)}}\cr\vdots\cr}
\Bigm|_{m\over2}T(p^\iota)=
\co{}\cdot\Bigl(\eich{}{\iota}-\beta_p\delta_{({{m+1}\over2})} 1_h \Bigr)
\cdot\pmatrix{\vdots\cr{{\th(z,\q_j)}\over{e(\q_j)}}\cr\vdots\cr}\ ,
\eqno(2.5)
$$
where
$$
\iota=\cases{1 & if $m=2k$ is even, \cr
             2 & if $m=2k+1$ is odd, \cr}
$$
the square matrix of oder $h$
$$
\eich{}{\iota}=\eich{\bf q}{\iota}=\Bigl(\repr{\iota}iji\Bigr)
\eqno(2.6)
$$
is the Eichler's ``Anzahlmatrix'' (with multiplier $p^{\iota}$)
and $\delta$ is generalized Kronecker symbol defined by
$$
\delta_{(X)}=\cases{1 &, if $X$ is an integral matrix, \cr
                    0 &, if matrix $X$ is not integral.\cr}
\eqno(2.7)
$$
\par
Finally, following C. L. Siegel one can consider the {\it generic} theta-series
$$
\th_{\{\q\}}(z)=\sum_{j=1}^{h} e^{-1}(\q_j)\th(z,\q_j)
$$
and show that for $p{\not|}\det\q$
it is an eigenfunction of corresponding Hecke operators $T(p^{\iota})$
with eigenvalues $\co{}\bigl(\sum_j\repr{\iota}{}jj-(\iota-1)\beta_p\bigr)$
, $\iota$ is either $1$ or $2$ depending on the parity of $m$.
In case $m$ is even and $\det\q=1$ the generic theta-series is proportional
to the Eisenstein series of weight $m\over2$. (This implies the famous Siegel
theorem on mean numbers of integral representations in the case of quadratic
forms of determinant $1\,$, see for example Exercise 5.1.18 in [3].)
\par
The above formulas together with remarkable multiplicative properties of
corresponding Hecke algebras (see [3]) reveal multiplicative relations between
the Fourier coefficients of the theta-series which are best expresses in terms
of associated zeta-functions:
$$
\sum_{(n,\det\q)=1} {{\brey{na}}\over{n^s}}=
\prod_{p\not{\,|}\det\q}
\Bigl({1\over{1-\co{} \beich{} p^{-s}+\xa{\q}p^{k-1-2s}}}\Bigr)
\cdot\brey{a}
\eqno(2.8)
$$
in case $m=2k$ is even and $a$ is any non-zero integer, or
$$
\sum_{(n,\det\q)=1} {{\brey{n^2 a}}\over{n^s}}=
\prod_{p\not{\,|}\det\q}
\Bigl({{1-\xa{\q}({{2a}\over p})p^{k-1-s}}\over
{1-\co{}(\beich{2}-\beta_p)p^{-s}+p^{m-2-2s}}}\Bigr)
\cdot\brey{a}
\eqno(2.9)
$$
in case $m=2k+1\geq3$ is odd and $a$ is a square-free integer.\P
For $\brey{n}$
we can substitute either the column-vector $\tr(\dots,\rep{j}n,\dots)$ or
the mean number $\sum_{j=1}^{h}\rep{j}n$. The first case corresponds to 
Fourier coefficients of the vector-valued modular form 
$\tr(\dots,{{\th(z,\q_j)}\over{e(\q_j)}},\dots)$
the second --- to Fourier coefficients of the generic theta-series
$\th_{\{\q\}}(z)$. Then
$\beich{\iota}$ stands for either Eichler's matrix (2.6) or the average
$\sum_{j=1}^{h}\repr{\iota}j{}j$ respectively. In any case the formulas
(2.8) and (2.9) give us explicit expressions of multiplicative relations
between the numbers $r(\q,n)$ of integral representations of integers by
quadratic forms. The well-known formula (2.8) easily follows from (2.5) and
Exercise 4.3.6 of [3] for example. A proof of the formula (2.9) follows below.
\par\medskip{\sl Proof of Formula (2.9). }
Let $\brey{n}=\tr(\dots,\rep{j}n,\dots)$ and let $\beich{2}$ be 
Eichler's matrix (2.6). Substituting quadratic form $\q_j\,,1\leq j\leq h$ 
(from the full system of representatives $\{\q_1,\dots,\q_h\}$) 
in place of $\q$ in (1.18), 
dividing both sides by $e(\q_j)$ and summing up over $K\in R^*(\q_j,p^2n)$ in
the same way we did for the proof of (2.2), we conclude that
$$\openup2pt\displaylines{
\qquad
c_p^{-1}\cdot\sum_{i=1}^h {{r^*(\q_j,p^2\q_i)}\over{e(\q_j)}}
\cdot{{r(\q_i,n)}\over{e(\q_i)}} =\hfill\cr
\hfill\left\{r(\q_j,p^2n)+\left(\eps_p(n)\xa{\q}p^{k-1}+c_p^{-1}\beta_p
\right)r(\q_j,n)+p^{m-2}r(\q_j,n\slash p^2)\right\}\,,
\qquad\cr
}$$
for any $j,\,1\leq j\leq h$ and any $n\in\N\cup{0}\,$. In matrix notations:
$$
c_p^{-1}\cdot\beich{2}\cdot\brey{n}=\brey{p^2n}+
\left(\eps_p(n)\xa{\q}p^{k-1}+c_p^{-1}\beta_p \right)\cdot\rey{}{n}+
p^{m-2}\cdot\rey{}{n\slash p^2}\quad.
$$
The last equality implies the following identities for the following formal
power series in $t$ with coefficients in $\Q^h\,$:
$$\openup2pt\displaylines{
\qquad\left(1_h-c_p^{-1}(\beich{2}-\beta_p)t+p^{m-2}t^2\right)\cdot
\sum_{\delta\geq0}\brey{p^{2\delta}a}t^{\delta}=\hfill\cr
\qquad\brey{a}+\left(\brey{p^2a}-c_p^{-1}(\beich{2}-\beta_p)\brey{a}\right)+
\hfill\cr
\hfill\sum_{\delta\geq2}\left\{\brey{p^{2\delta}a}-c_p^{-1}(\beich{2}-\beta_p)
\brey{p^{2(\delta-1)}a}+p^{m-2}\brey{p^{2(\delta-2)}a}\right\}t^{\delta}=
\qquad\cr
\hfill\left(1_h-\eps_p(a)\xa{\q}p^{k-1}t\right)\cdot\brey{a}\,,\qquad\cr
}$$
where $a$ is an arbitrary integer not divisible by $p^2\,$. Define formal power
series in $t$ with coefficients in $\Z^h_h$ by
$$
\sum_{\delta\geq0}\heich{p^{2\delta}}t^{\delta}=
{{1_h-\eps_p(a)\xa{\q}p^{k-1}t}\over{1_h-c_p^{-1}(\beich{2}-\beta_p)t+
p^{m-2}t^2}}\quad,
$$
in other words
$$\left\{\matrix{
\heich{1} &=& 1_h\,,\hfill\cr
\heich{p^2} &=& \left(c_p^{-1}(\beich{2}\beta_p)-\eps_p(a)\xa{\q}p^{k-1}
\right)\cdot1_h\,,\hfill\cr
\heich{p^{2\delta}} &=& c_p^{-1}(\beich{2}-\beta_p)\heich(p^{2(\delta-1)})
-p^{m-2}\heich{p^{2(\delta-2)}}\,,\hbox{ for }\delta\geq2\,.\hfill\cr
}\right.
$$
Then formally
$$
\sum_{\delta\geq0}\brey{p^{2\delta}a}t^{\delta}=
\sum_{\delta\geq0}\heich{p^{2\delta}}t^{\delta}\cdot\brey{a}
$$
and so $\brey{p^{2\delta}a}=\heich{p^{2\delta}}\cdot\brey{a}$ for 
$\delta\geq0$ and $a$ not divisible by $p^2\,$. Define
$$
\heich{n^2}=\prod_{s=1}^l\heich{p_s^{2{\delta}_s}}\,,\hbox{ if }\ 
n=p_1^{{\delta}_1}\cdot\dots\cdot p_l^{{\delta}_l}\quad,
$$
here $p_j$'s are distinct prime factors of $n$ and ${\delta}_j$'s are their
corresponding orders. Note that the definition does not depend on the oder of
the factors. To see this it is enough to note that according to Lemmas 3.3--3.6
of [4], Eichler matrices $\beich{2}\,,\,\bar{\frak t}^*(u^2)$ commute for
different primes  $p,u$ not dividing $\det\q\,$. Next we claim that
$$
\brey{n^2a}\heich{n^2}\cdot\brey{a}
$$
for any integer $n$ coprime to $\det\q\,$, and any integer $a$ not divisible
by the square of any prime factor of $n\,$. Indeed, let $n=p_1^{{\delta}_1}
\cdot\dots\cdot p_l^{{\delta}_l}$ be coprime to $\det\q\,$. We will use 
induction on $l$ to justify our claim. If $l=1$ then the statement is obvious
because of the definition of $\heich{n^2}=\heich{p_1^{2{\delta}_1}}\,$. Assume
now that if $b$ is a product of powers of $p_1,\dots,p_{l-1}$ then
$\brey{b^2a}\heich{b^2}\cdot\brey{a}$ for any $a$ not divisible by $p_i^2\,$,
$i$ ranging from $1$ to $l-1\,$. Fix an integer $a$ not divisible by any of
$p_i^2 ,\,1\leq i\leq l\,$, then
$$
\openup2pt\displaylines{
\qquad
\brey{n^2a}=\brey{\prod\nolimits_{i=1}^{l-1}p_i^{2{\delta}_i}\cdot
p_l^{2{\delta}_l}a}=
\hat{\frak t}^*_{(p_l^{2{\delta}_l}a)}
(\prod\nolimits_{i=1}^{l-1}p_i^{2{\delta}_i})
\cdot\brey{p_l^{2{\delta}_l}a}=\hfill\cr
\hfill\heich{\prod\nolimits_{i=1}^{l-1}p_i^{2{\delta}_i}}\cdot
\heich{p_l^{2{\delta}_l}}
\cdot\brey{a}=\heich{n^2}\cdot\brey{a}\qquad\cr
}$$
since $\hat{\frak t}^*_{(u^2a)}(p^2)=\heich{p^2}$ if $u$ and $p$ are coprime.
By induction on $l$ we conclude that $\brey{n^2a}\heich{n^2}\cdot\brey{a}$ for
any integer $n$ coprime to $\det\q$ and any integer $a$ not divisible by the
square of a prime factor of $n\,$ (in particular for any square-free integer
$a\,$). Combining the above equalities we deduce the following identities for
the formal zeta-functions:
$$
\openup2pt\displaylines{
\qquad\sum_{(n,\det\q)=1} {{\brey{n^2 a}}\over{n^s}}=
\sum_{(n,\det\q)=1} {{\heich{n^2}}\over{n^s}}\cdot\brey{a}=\hfill\cr
\hfill\prod_{p\not{\,|}\det\q}\Bigl(\sum_{\delta\geq0}
{{\heich{p^{2\delta}}}\over{p^{s\delta}}}\Bigr)\cdot\brey{a}=
\prod_{p\not{\,|}\det\q}
\Bigl({{1-\xa{\q}({{2a}\over p})p^{k-1-s}}\over
{1-\co{}(\beich{2}-\beta_p)p^{-s}+p^{m-2-2s}}}\Bigr)
\cdot\brey{a}\qquad\cr
}$$
for any square-free integer $a\,$ which finishes the proof of (2.9) for the 
case of column-vectors $\brey{n}=\tr(\dots,\rep{j}n,\dots)\,$. Proof of (2.9)
for the case of the mean numbers $\sum_j r(\q_j,n)\slash e(\q_j)$ proceeds in
a similar fashion. One just needs to multiply both sides of the matrix 
equalities above by the row $(1,\dots,1)$ of length $h$ and observe that
$$
(1,\dots,1)\cdot\left(\repr{\iota}iji\right)=
\left(c_p(1+p^{m-2})+\beta_p\right)\cdot(1,\dots,1)
$$
according to comparison of the $0$th Fourier coefficients in (2.2).
\flushqed
Let us note also that the zeta-functions on the right hand side of (2.8)
are closely related to Epstein's zeta-functions
$$
\zeta(s,\q_{j})=\sum_{n\geq1} {{r(\q_{j},n)}\over{n^s}}
\eqno(2.10)
$$
which are Mellin transforms of corresponding theta-series $\th(z,\q_{j})$
and therefore have good analytical properties (i.e. meromorphic continuation on
the whole $s$-plane and standard functional equation, for more details see
[8], [12] or exercises in Section 4.3.1 of [3]).
\vfill\break
%
%
\bigskip\noindent
{\bf 3. Matrix Hecke rings of orthogonal groups}
\smallskip\noindent
For a long time it was somewhat of a mystery that the numbers $r(\q,n)$ of
integral solutions of certain quadratic equations have general multiplicative
properties even though no similar relations between the solutions themselves
were known. (The sole exception to this rule was the case of so-called 
composition of quadratic forms. Then the integral representations can be
interpreted as elements of certain arithmetical rings and multiplicative
structure of these rings is reflected in relations between representations.
A classical example of this phenomenon is Gauss theory of binary quadratic
forms. Unfortunately, according to a theorem of Hurwitz, there is no 
composition unless the number of variables equals to $1,2,4$ or $8$.)

In [1] A. Andrianov developed Shimura's construction of abstract matrix Hecke
algebras (see [11]) specifically for the case of orthogonal groups and
introduced so-called {\it automorph class rings}. The latter play the role
of Hecke operators acting directly on the sets of solutions of quadratic
Diophantine equations. The basic idea behind this action is that 
$$
R(\q,a\q')\cdot R(\q',b\q'')\subset R(\q,ab\q'')
$$
for any integral quadratic forms $\q\, ,\q',\q''$ and for any integers $a,b$ 
(the dot refers to usual matrix multiplication).
This purely algebraic approach revealed existence of general
multiplicative relations between automorphs $R(\q,n\q')$ and representations
$R(\q,n)$ of integers by quadratic forms. Remarkably, these relations can be
expressed in terms of associated zeta-functions in the form almost identical
to (2.8) and (2.9). In fact, the latter formulas for numbers of representations
appear to be immediate consequences of the former relations between 
representations themselves !  This suggests existence of some correlation
between rings of Hecke operators for symplectic and for orthogonal groups.

Let us examine the construction in some detail. For an integral nonsingular
quadratic form $\q$ in $m$ variables we fix a complete system $\{\q_1,\dots,\q_h\}$
of representatives of equivalence classes contained in the similarity class
of $\q$. Consider free $\Z$-modules $\A_{ij},\ 1\leq i,j \leq h$ generated by
double cosets
$$
\bigl\{E(\q_i)AE(\q_j)\ ,\hbox{ where }
A\in E(\q_i)\backslash\bigcup\limits_{a\in\N}R(\q_i,a\q_j)/E(\q_j)\bigr\}.
$$
Since any double coset of the above type is a finite disjoint union of left
cosets modulo $E(\q_i)$, we will write elements of $\A_{ij}$ as finite linear
combinations with integral coefficients of (formal) symbols $\y{A}$ which
bijectively correspond to left cosets 
$E(\q_i)A\ ,\ A\in E(\q_i)\backslash\bigcup_{a\in\N}R(\q_i,a\q_j)$.
Note that elements of $\A_{ij}$ are invariant under the right multiplication
by matrices from $E(\q_j)$. This and the obvious inclusion
$$R(\q_i,a\q_j)\cdot R(\q_j,b\q_k)\subset R(\q_i,ab\q_k)$$
allows us to define natural multiplication
$$\openup2pt\displaylines{
\A_{ij}\times\A_{jk}\longrightarrow\A_{ik}\ , \cr
\sum_{\alpha}a_{\alpha}\y{A_{\alpha}}\cdot\sum_{\beta}b_{\beta}\y{B_{\beta}}=
\sum_{\alpha ,\beta}a_{\alpha}b_{\beta}\y{A_{\alpha}B_{\beta}}\ . \cr
}$$
Finally, we introduce the ring of $h\times h$ matrices
$$
\A=\A_{\{\q\}}=\bigl\{(A_{ij})\  ;\ A_{ij}\in\A_{ij}\, ,\, 1\leq i,j \leq h\bigr\}
$$
which is the {\it automorph class ring} (over $\Z$) of the form $\q$.
Sometimes it is more convenient to consider $\A\otimes\Q$ ---
the automorph class ring over $\Q$.

The automorph counterparts of the classical Hecke operators $T(p^{\iota})$
($\iota$ is $1$ or $2$ depending on the parity of $m$) are elements of the form
$$
\Eich{}{\iota}=\Eich{\bf q}{\iota}=\Bigl(\Eich{ij}{\iota}\Bigr)\,\hbox{, where }
\Eich{ij}{\iota}=\sum_{D\in\Repr{\iota}{i}{j}}\y{D}\ .
\eqno(3.1)
$$
Note that these matrices are close relatives of Eichler's ``Anzahlmatrizen"
$(2.6)$.

Next we should define an action of $\A$ on the integral representations by the
forms $\q_1,\dots,\q_h\,$. The idea is based on the natural inclusion:
$R(\q_i,n\q_j)\cdot R(\q_j,a)\subset R(\q_i,na)\,$. In order to suite our
formalism we introduce column-vectors of length $h$ :
$$
\Rey{}{a}=\Rey{\bf q}{a}=\pmatrix{\vdots\cr{\Rey{j}{a}}\cr\vdots\cr}\,
\hbox{, where }
\Rey{j}{a}=\sum_{L\in\Rep{j}{a}}{\mu_{\q_j}^{-1}(L)}\cdot\y{L}\ ,
\eqno(3.2)
$$
which encode all orbits modulo $E(\q_j)$ of integral representations of $a$ by
forms $\q_j$. Here $\mu_{\q_j}(L)$ is a suitably normalized measure of the
stabilizer $E_L(\q_j)=\{U\in E(\q_j)\ ;\ UL=L\}\,$.
The vectors of type (3.2) can be thought of as elements of $\Q$-module
$\D=\D_{\q}=\D_1\times\dots\times\D_h$ whose $j^{th}$ component $\D_j$ is a
free $\Q$-module generated by (formal) symbols $\y{L}$ which are in one-to-one
correspondence with orbits $E(\q_j)L \hbox{ in } E(\q_j)\backslash\Z^{m}\,$.

In this notations an automorph analog of Eichler's commutation relation $(2.5)$
take the form:
$$
\Rey{}{pa}+\xa{\q}p^{k-1}[p]\Rey{}{a\slash p}\ =\ \co{}\Eich{}{}\Rey{}{a}\, ,
\eqno(3.3)
$$
if $m=2k\,$, or
$$
\Rey{}{p^2a}+\xa{\q}p^{k-1}\Bigl({{2a}\over p}\Bigr)
[p]\Rey{}{a}+p^{m-2}[p^2]\Rey{}{a\slash{p^2}}\ 
=\ \co{}\Bigl(\Eich{}{2}-\beta_p 1_h\Bigr)\Rey{}{a}\, ,
\eqno(3.4)
$$
if $m=2k+1\geq 3\,$.
Here $p$ is a prime coprime to $\det\q$, $a$ is any positive integer and
$[p^{\iota}]=diag(\y{p^{\iota}1_m},\dots,\y{p^{\iota}1_m})\in\A$ are 
simple elements of the automorph class ring responsible for nonprimitive
automorphs. The automorph class theory formalism was initially developed
in [1], where one can find the original proof if the formula (3.3) above
as formula (1.18) of Section 1. The formula (3.4) was proved in [4], 
Theorem 3.1 (for the case of positive-definite quadratic forms) and in
(2.8) of [5] for the case of arbitrary nonsingular forms.
\par
The above formulas relating sets of representations $R(\q,p^{\iota}a)$ and
$R(\q,a)$ immediately imply analogous relations between numbers
$r(\q,p^{\iota}a)$ and $r(\q,a)$ of representations. Namely, let us consider
two `` coefficient'' homomorphisms $\pi$:
$$
\A\longrightarrow\M_{h}(\Z)\quad\hbox{ or }\quad\D\longrightarrow\Z^{h}
$$
defined entry-wise by
$$
\pi:\ \sum_{\alpha}a_{\alpha}\y{A_\alpha}\ \longmapsto\sum_{\alpha}a_{\alpha}
$$
for a corresponding finite formal linear combination of ``orbits'' $\y{\cdot}\,$.
For example,
$$
\pi\Bigl(\Eich{\bf q}{\iota}\Bigr)=
\eich{\bf q}{\iota}=
\Bigl(\bigl|\Repr{\iota}{i}{j}\bigr|\Bigr)
$$
is the Eichler's ``Anzahlmatrix'' $(2.6)\,$. In case $\q$ is positive definite
all sets of representations under consideration are finite, in particular
$$
\pi\Bigl(\Rey{\q}{n}\Bigr)=
\rey{\bf q}{n}=
\pmatrix{\vdots\cr\rep{j}{n}\cr\vdots\cr}
$$
coincides with the $n^{th}$ Fourier coefficient of the vector modular
form $\tr\bigl(\dots,{{\th(z,\q_j)}\over{e(\q_j)}},\dots\bigr)\,$. Then application
of homomorphisms $\pi$ to $(3.3)\hbox{ and }(3.4)$ leads directly to Eichler's
identities $(2.5)$ written in terms of Fourier coefficients of corresponding
theta-series. 
Furthermore, automorph class ring $\A$ itself has good multiplicative
properties (as one would expect from a Hecke algebra). Combining this 
properties with $(3.3)\hbox{ or }(3.4)$ one deduces Euler product expansions
of (formal) zeta-functions
$$
\sum\limits_{(n,\det\q)=1}{{\Rey{}{n}}\over{n^s}}\qquad
\hbox{ or }\qquad
\sum\limits_{(n,\det\q)=1}{{\Rey{}{n^2a}}\over{n^s}}
$$
very similar to Euler products $(2.8)$ or $(2.9)$ respectively (for exact
formulas and their detailed proofs see formula (1.19) of [1] and formula (1.13)
with Theorem 1.1 of [5]). This method
also provides an alternative proof of identities $(2.8)\, ,\,(2.9)$
themselves via appropriate application of the ``coefficient'' homomorphisms
$\pi\,$. One should stress here once again that all of the above considerations
are purely algebraic and remain valid for any nondegenerate quadratic form
over a Dedekind domain. (For an extensive account of multiplicative properties
of general automorph class rings and of Euler products of associated
formal zeta-functions the reader is referred to [1],[5] and especially [2].)
\vfill\break
%
%
\bigskip\noindent
{\bf 4. Shimura's lift for theta-series}
\smallskip\noindent
In [12] G. Shimura showed that if a cusp form
$f(z)=\sum_{n\geq 1}c(n)\exp(2\pi i n z)$ of weight $m/2$ (where $m=2k+1\geq3$)
is a common eigenfunction for all Hecke operators $T(p^2)$, then zeta-function
$\sum_{n\geq1} c(n^2 a)/n^s$ has Euler product decomposition of type $(2.9)$
for every square-free positive integer $a$ (one should just replace
$\bar\rey{}{\cdot}\hbox{ by }c(\cdot)\hbox{ and }\co{}(\bar\eich{}{2}-\beta_p)$
by corresponding eigenvalues of Hecke operators). Comparing denominators of
local factors in $(2.8)$ and $(2.9)$ one can note their striking similarity.
And indeed, Shimura used the Weil criterion to prove that the zeta-function
$$
\biggl(\sum_{n\geq1}{{\chi_{f}(n)\bigl({a\over n}\bigr)n^{k-1}}\over{n^s}}\biggr)
\cdot
\biggl(\sum_{n\geq1}{{c(n^2 a)}\over{n^s}}\biggr)
\eqno(4.1)
$$
defined by the Euler product of denominators of local factors of
$\sum_{n\geq1}{{c(n^2 a)}/{n^s}}$ is a Mellin transform of a 
modular form $f^{(a)}(z)$ of integral weight $m-1=2k$ (see [12], 
Main theorem). The form $c^{-1}(a)\cdot f^{(a)}(z)$ is independent
of $a$ and is called the {\it Shimura's lift} of $f(z)$. In what follows we
will use this name for $f^{(a)}(z)$ as well.

P. Ponomarev in [10] investigated the lift for theta-series $\th(z,\q)$ of
certain ternary positive definite quadratic forms. The theta-series are
neither cusp forms nor Hecke eigenforms in general,
but one can still consider lift
$\th^{(a)}(z,\q)$ defined for an integer $a$ via product of type $(4.1)$.
As the reader already expects, in order to get a complete picture one should
work not with an individual quadratic form $\q$ but rather with a complete
system
$$
\{\q_1,\dots,\q_h\}
\eqno (4.2)
$$
of representatives of different equivalence classes of the similarity class of
$\q$. Employing particular  case $m=3$ of Eichler's commutation relation $(2.5)$
P. Ponomarev showed that
$$
\th^{(a)}(z,\q_k)=\sum_{j=1}^h
\biggl(r(\q_j,a)\cdot\sum_{n\geq 0}\pi_{kj}(n^2)\,e^{2\pi inz}\biggr)\ ,
\ 1\leq k\leq h\, ,
$$
where each $\sum_{n\geq0}\pi_{kj}(n^2)\exp(2\pi inz)$ is a finite linear
combination of theta-series of quaternary quadratic forms associated with
certain lattices in a quaternion algebra (see [10], Theorem 1).
The quaternion algebra is chosen so that the reduced norm on its pure
 part is similar over $\Q$ to the form $\q$. In particular then
$\bigl(\pi_{kj}(p^2)\bigr)=\bigl(r^*(\q_k,p^2\q_j)/e(\q_j)\bigr)$ is nothing
else but transposed of Eichler's matrix $(2.6)$ . 

One of the goals of the present paper is to find an automorph analog of 
Shimura's lift for theta-series of ternary quadratic forms. Let us start
with a brief reexamination of P. Ponomarev's main result from a somewhat
different point of view. For simplicity assume that $\q$ is an integral
ternary  positive definite quadratic form of class number $h$. Fix the
system $(4.2)$ of representatives from the similitude class. Then the
generic theta-series $\th_{\{\q\}}(z)$ is an eigenform of all Hecke operators
$T(p^2)$ with $p{\not|}\det\q$ and the equality $(2.9)$ holds for associated
zeta-functions. From $(2.3)$ it follows that for ternary quadratic forms
the local constants in $(2.9)$ are given by: $\co{}=1\, ,\,\beta_p=0\,$.
Define coefficients $b_n$ by
$$
\sum_n {{b_n}\over{n^s}} =
\prod_p \bigl(1-\bar\eich{}{2}p^{-s}+p^{1-2s}\bigr)^{-1}\ ,
$$
where
$$
\bar\eich{}{2}=\bar\eich{\q}{2}=\sum_{j=1}^h \repr2j{}j\ ,
$$
and put
$$
B(z)=\sum_n b_n e^{2\pi inz}\ .
$$
Note in particular that $b_1=1\, ,\, b_p=\bar\eich{}{2}\,$.
We expect $B(z)$ to be equal to a finite linear combination
$$
B(z)\,=\,\sum_{j=1}^H x_j\th(z,\n_j)\,=\,
\sum_n\biggl(\sum_{j=1}^H x_j r(\n_j,n)\biggr)e^{2\pi inz}
$$
of theta-series associated to some quaternary (integral) quadratic forms $\n_j\,$.
We also expect $B(z)$ to satisfy
$$
B(z)\bigm|_2 T(p)\,=\,\bar\eich{\q}{2}B(z)
$$
for any prime $p\, ,\, p{\not|}\det\q\,$. A plausible candidate would be
$$
B(z)=
\Bigl(\sum_{i=1}^H\repn{i}1\Bigr)^{-1}\cdot\sum_{j=1}^H e^{-1}(\n_j)\th(z,\n_j)
\eqno (4.3)
$$
-- proportional to generic theta-series $\th_{\{\n\}}(z)$ of a quaternary
form $\n$ of class number $H\,$. In this case we would have
$$
\bar\eich{\q}2=b_p=
\Bigl(\sum_{i=1}^H\repn{i}1\Bigr)^{-1}\cdot\sum_{j=1}^H\repn{j}p\ .
\eqno (4.4)
$$
On the other hand from Eichler's commutation relations $(3.3)$ or $(2.5)$
we know that
$$
\pmatrix{\vdots\cr\repn{i}p\cr\vdots\cr}\ =\ 
\co{(\n)}\biggl(\reprn{}i{j}i\biggr)
\pmatrix{\vdots\cr\repn{j}1\cr\vdots\cr}\quad ,
$$
which implies that
$$
\sum_{i=1}^H\repn{i}p\,=\,
\co{(\n)}\sum_{j=1}^H\Bigl(\sum_{i=1}^H\reprn{}i{j}i\Bigr)\repn{j}1 \quad .
$$
But the sum $\bigl(\sum_i r(\n_i,p\n_j)/e(\n_i)\bigr)$ does not depend on $j\,$.
(In fact, the sum coincides with $c_p(\n)$ times zero Fourier coefficient of
$\th(z,\n_j)|_2 T(p)$ which is equal to $1+p$ for any $j$). Therefore,
together with hypothetical $(4.4)$, the last equality would imply that
$$
\sum_{i=1}^h\repr2i{}i=\bar\eich{\q}2=\co{(\n)}\sum_{i=1}^H\reprn{}i{}i \quad ,
$$
where $c_{p}(\n)=(1+\xa{\n})$ according to $(2.3)\,$. The above formula can
be rewritten in a slightly more general form
$$
\sum_{A\in\bigcup_{i=1}^h R^{*}(\q,p^2\q_i)/E(\q_i)} 1 \qquad =\qquad
(1+\xa{\n})^{-1}\sum_{M\in\bigcup_{i=1}^H R(\n,p\n_i)/E(\n_i)} 1
\eqno (4.5)
$$
which make sense for any {\it nonsingular} (integral) ternary form $\q\,$.
Thus according to our heuristic argument, one can expect that an automorph
$A\in R(\q,p^2\q')$ of quadratic form $\q$ can be ``lifted'' to an automorph
$M\in R(\n,p\n')$ of some quaternary form $\n$ associated to $\q\,$.
Moreover, the equality $(4.5)$ should appear as a corollary of the expected
lifting of automorphs. In the next sections we will prove that this is exactly
the case by providing an explicit construction for the ``lift'' (see (5.30),
Theorems 6.3, 6.7 and Corollary 6.8).
\vfill\break
%
%
\bigskip\noindent
{\bf 5. Clifford algebras and automorph lift}
\smallskip\noindent
Since our construction of the automorph lift is based on Clifford algebras,
we start with a brief recall of their properties. In the course we also
introduce some convenient notations.
Let
$$
\q(X)=q_1x_1^2+b_{12}x_1x_2+b_{13}x_1x_3+q_2x_2^2+b_{23}x_2x_3+q_3x_3^2
\eqno (5.1)
$$
be an integral nonsingular (ternary) quadratic form. We set
$$
Q_0=\pmatrix{
q_1 & b_{12} & b_{13} \cr
0   & q_2    & b_{23} \cr
0   & 0      & q_3    \cr
}\qquad ,\qquad
B=\pmatrix{-b_{23} \cr b_{13} \cr -b_{12} \cr}\ .
\eqno (5.2)
$$
Then the matrix $Q$ of the form $\q$ and its determinant $\det\q$
are given by
$$
Q=\tr Q_0+Q_0 \quad ,\quad 2\Delta=\det\q=2\cdot(4\,\det{Q_0}-Q_0[B])\ ,
\eqno (5.3)
$$
where $\det{Q_0}=q_1q_2q_3\,$ and $\Delta=\det\q/2\,$.
The matrices adjoint to $Q_0$ and $Q$ are
$$
\tilde{Q_0}=\pmatrix{
q_2q_3 & -q_3b_{12} & b_{12}b_{23}-q_2b_{13} \cr
0      & q_1q_3     & -q_1b_{23}             \cr
0      & 0          & q_1q_2                 \cr
}
$$
and
$$
\tilde Q=2(\tr\tilde{Q_0} +\tilde{Q_0})-B\cdot\tr B\ .
\eqno (5.4)
$$
We associate to $\q$ quadratic $\Z$-module $(E,\q)=(\Z^3,\q)$ with
corresponding symmetric bilinear form 
$$
\b(x,y)=\q(x+y)-\q(x)-\q(y)\ ,\qquad x,y\in E
\eqno (5.5)
$$
whose matrix $\bigl(\b(e_i,e_j)\bigr)$ relative to the standard basis
$e_1,e_2,e_3$ of $\Z^3$ is $Q\,$. The {\it Clifford algebra} 
$C(E)=C(E,\q)$ of $(E,\q)$ is a free $\Z$-algebra with monomorphism
$\iota:E\hookrightarrow C(E)$ such that $\iota(x)^2=\q(x)\cdot 1_{C(E)}$
for any $x\in E\,$. It is uniquely characterized by the following 
{\it universal} property: for any homomorphism $\eta:E\rightarrow D$ into
an algebra $D$ with $\eta(x)^2=\q(x)\cdot 1_D\,$, there exists unique 
homomorphism of algebras $\tau:C(E)\rightarrow D$ such that the diagram
$$
\diagram{
E & \lharr{id}{} & E \cr
\dvarr{\eta}{} && \dvarr{}{\iota} \cr
D & \lharr{\tau}{} & C(E) \cr
}
\eqno (5.6)
$$
commutes (i.e. $\eta=\tau\circ\iota\,$). The algebra $C(E)$ can be realized
as a free $\Z$-module of rank $8$ generated by products of $e_1,e_2,e_3$
(see [9] for details). We will denote the product $e_i e_j$ in $C(E)$ by
$e_{ij}$ and $e_i e_j e_k$ by $e_{ijk}\,$. The unit element of $C(E)$ is
denoted by $e_0$ and we identify $\Z$ with $\Z e_0\,$. We also identify
$E$ with $\iota(E)\subset C(E)\,$. The {\it even subalgebra}
$$
C_0(E)=\Z e_0\oplus\Z e_{23}\oplus\Z (-e_{13})\oplus\Z e_{12}
\eqno (5.7)
$$
is a free $\Z$-module of rank $4$ generated by products of even number of
vectors of $E\,$. (Note the order of the elements of the basis in (5.7) !
In spite of first strange impression, this is {\it the natural} oder. We
will refer to this basis of $C_0(E)$ as the {\it natural basis}). 
Setting further $C_1(E)=\Z (-e_{123})\oplus \Z e_1\oplus\Z e_2\oplus\Z e_3$
one has $C(E)=C_0(E)\oplus C_1(E)\,$. Sometimes it is also useful to consider
$C(E)\otimes\Q$ or $C_0(E)\otimes\Q$ -- the Clifford algebra of $(E,\q)$
over $\Q$ and its even subalgebra, in which $C(E)$ and $C_0(E)$ are orders.
The algebra $C(E)$ has an
antiautomorphism $x\mapsto\bar x$ defined on the generators by
$$
\bar{e_j}=-e_j\ ,\qquad j=1,2,3.
\eqno (5.8)
$$
The restriction of the antiautomorphism $x\mapsto\bar x$ to $C_0(E)$
can be written in coordinates as
$$
x_0 e_0 + x_1 e_{23} - x_2 e_{13} + x_3 e_{12} \mapsto
(x_0+x_3 b_{12}-x_2 b_{13}+x_1 b_{23})e_0-x_1e_{23}+x_2e_{13}-x_3e_{12}
\eqno (5.9)
$$
and gives us the {\it standard involution} on $C_0(E)\,$. That is to say,
$$
\s(x)=x+\bar x\in\Z\quad ,\quad \n(x)=x\cdot\bar x\in\Z
\eqno (5.10)
$$
for any $x\in C_0(E)\,$. In particular each element $x$ of $C_0(E)$
satisfies the equation $(z-x)(z-\bar x)=z^2-\s(x)z+\n(x)=0$ which is
unique in case $x$ and $e_0$ are linearly independent. Integers $\s(x)$
and $\n(x)$ are called respectively the {\it trace} and the {\it norm}
of $x\in C_0(E)\,$. The norm $\n$ turns the even subalgebra $C_0$ into
quadratic $\Z$-module $\bigl(C_0(E),\n\bigr)$ of rank $4$ with associated
bilinear form 
$$
\s(x\bar y)=\n(x+y)-\n(x)-\n(y)\ ,\quad x,y\in C_0(E)\ .
\eqno (5.11)
$$
The matrix $N$ of the form $\n$ with respect to the natural basis (5.7) is
$$
N=\pmatrix{
2 & -\tr B \cr
  &    \cr
-B & (\tr\tilde{Q_0}+\tilde{Q_0}) \cr
}
\eqno (5.12)
$$
One easily computes the determinant
$$
\det\n={\Delta}^2=\bigl((\det\q)/2\bigr)^2
\eqno (5.13)
$$
and the inverse matrix 
$$
N^{-1}={\Delta}^{-1}\cdot\pmatrix{
2\cdot\det{Q_0} & \tr(Q_0 B) \cr
     &   \cr
(Q_0 B) & Q \cr
}\ ,
\eqno (5.14)
$$
where $Q_0$ and $B$ are given by (5.2).
From the above
formulas one can see that the quadratic modules $(E,\q)$ and $(C_0(E),\n)$
are closely related. In fact $(C_0(E),\n)$ contains a certain quadratic
submodule of rank $3$ which is ``almost'' $(E,\q)\,$. We will present the
construction shortly, but first we need to introduce the special element
$$
t=e_{123}+\overline{e_{123}}=e_{123}-e_{321}\in C(E)
\eqno (5.15)
$$
used by M. Kneser in [9]. Straightforward (but rather tiresome) calculations
show that $t$ belongs to the center of $C(E)$, $t=\bar t$ and
$$
t=\bigl(-e_{123},e_1,e_2,e_3\bigr)\pmatrix{-2\cr B\cr}\ ,\quad
t^2=-\Delta\quad .
\eqno (5.16)
$$
Consider now $\Z$-submodule $Et$ of $C_0(E)\,$. Using (5.16) we see that
$Et$ has rank $3$ and
$$
\bigl(e_1t,e_2t,e_3t\bigr)=\bigl(e_0,e_{23},-e_{13},e_{12}\bigr)\cdot T
$$
form its $\Z$-basis, where
$$
T=\pmatrix{\tr(Q_0 B)\cr \cr Q\cr}\in{\Z}^4_3\quad .
\eqno (5.17)
$$
Let us compute the matrix of the restriction of the norm-form $\n$ onto
$Et$ with respect to the basis $e_1t,e_2t,e_3t\,$. We have
$$
N\cdot T=\pmatrix{
2\cdot\tr(Q_0 B)-\tr B\cdot Q \cr
                      \cr 
(\tr\tilde{Q_0} +\tilde{Q_0})(\tr{Q_0}+ Q_0)-B\cdot\tr B\cdot\tr{Q_0} \cr}=
\pmatrix{0\,0\,0 \cr {\scriptstyle{\Delta}}\cdot 1_3\cr}
$$
and so
$$
N[T]={\Delta}\cdot Q\quad .
\eqno (5.18)
$$
Recall that $\Delta=\det\q/2$ is an integer (see (5.3)), thus
$T\in R(\n,{\Delta}\q)$ and
$\bigl(Et,{\Delta}\q\bigr)
\hookrightarrow\bigl(C_0(E),\n\bigr)$ as we stated. Finally we also note
that
$$
-e_{123}\cdot t = \bigl(e_0,e_{23},-e_{13},e_{12}\bigr)
\pmatrix{ 2\,\det{Q_0}\cr \cr Q_0 B \cr}\quad ,
\eqno (5.19)
$$
which extends multiplication of $E$ by $t$ to injection
$C_1(E)\cdot t\hookrightarrow C_0(E)\,$, whose matrix with respect to our
choice of bases is $\Delta\cdot N^{-1}\,$.

Let $\q''$ be another integral ternary quadratic form. We associate with it
a quadratic module $(E'',\q'')$ and Clifford algebra $C(E'')$ exactly as above.
Let $A\in R(\q,\q'')$ be an integral automorph. From geometric point of view
$A$ defines an isometry $\alpha=\alpha_A$ from quadratic module
$(E'',\q'')$ into $(E,\q)\,$:
$$
\bigl(\alpha_A(e_1''),\alpha_A(e_2''),\alpha_A(e_3'')\bigr)=
\bigl(e_1,e_2,e_3\bigr)\cdot A\quad ,
$$
where $e_1'',e_2'',e_3''$ is the standard basis of $E''=\Z^3\,$. Let
$\iota$ and $\iota''$ be natural inclusions of $(E,\q)$ and $(E'',\q'')$
into corresponding Clifford algebras. 
Since $\bigl(\iota\circ\alpha(x'')\bigr)^2=\q\bigl(\alpha(x'')\bigr)\cdot
e_0=\q''(x'')\cdot e_0$ for any $x''\in E''\,$, there exists a unique
homomorphism of algebras $\varphi=\varphi_A$ which makes the diagram
$$
\diagram{
(E,\q) & \lharr{\alpha_A}{} & (E'',\q'') \cr
\dvarr{\iota}{} & & \dvarr{}{\iota''} \cr
C(E) & \lharr{\varphi_A}{} & C(E'') \cr
}
\eqno (5.20)
$$
commutative. (This follows immediately from the universal property (5.6)
of Clifford algebras.) Clearly the homomorphism $\varphi_A$ is compatible
with the involutions  (5.8) on $C(E)$ and $C(E'')\,$:
$$
\varphi_A(\bar x'')=\overline{\varphi_A(x'')}\  \hbox{for all }x''\in C(E'')\ .
\eqno (5.21)
$$
Therefore the restriction of $\varphi_A$ to the even subalgebra $C_0(E'')$ is
an isometry of quadratic modules
$\varphi_A:\bigl(C_0(E''),\n''\bigr)\to\bigl(C_0(E),\n\bigr)$ whose matrix
with respect to our choice of bases (5.7) is
$$
\Phi_A=\pmatrix{
1 & \tr A_3 Q_0 A_2 & -\tr A_3 Q_0 A_1 & \tr A_2 Q_0 A_1 \cr
0 & & &  \cr
0 & & \tr\tilde A & \cr
0 & & &  \cr
}\quad ,
\eqno (5.22)
$$
so that
$$
\bigl(\varphi_A(e_0''),\varphi_A(e_{23}''),\varphi_A(-e_{13}''),\varphi_A(e_{12}'')\bigr)=
\bigl(e_0,e_{23},-e_{13},e_{12}\bigl)\cdot\Phi_A\quad .
$$
Here $A_j$ denotes the $j^{th}$ column of the original automorph
$A\in R(\q,\q'')$ which gave rise to $\Phi_A\in R(\n,\n'')\,$.
(By $\n,\n''$ we mean the integral quaternary quadratic forms which are
specializations of the corresponding norm-forms at the natural bases (5.7)
of $C_0(E)$ and $C_0(E'')$ respectively.) To complete the picture we also
need a description of action of $\varphi_A$ on the special element
$t''\in C(E'')$ of the form (5.15):
\Lemma{5.1}{\sl Let } $A=(a_{ij})\in R(\q,\q'')$
{\sl be an integral automorph. Then }
$$
\varphi_A(t'')=(\det A)\cdot t
$$
{\sl in the notations introduced above.}
\Proof Note that elements of $E$ are multiplied in $C(E)$ according
to the following rule:
$$
\pmatrix{e_1,e_2,e_3}\pmatrix{x_1 \cr x_2 \cr x_3 \cr}
\cdot
\pmatrix{e_1,e_2,e_3}\pmatrix{y_1 \cr y_2 \cr y_3 \cr}=
\pmatrix{e_0,e_{23},-e_{13},e_{12}}
\pmatrix{\tr Y Q_0 X \cr x_2y_3-x_3y_2 \cr x_3y_1-x_1y_3 \cr x_1y_2-x_2y_1 \cr}
\quad ,
$$
Therefore
$$\openup2pt \displaylines{
\varphi_A(-e_{123}'')=-\varphi_A(e_1'')\varphi_A(e_2'')\varphi_A(e_3'')=
-\alpha_A(e_1'')\alpha_A(e_2'')\alpha_A(e_3'')= \cr
-1\cdot\pmatrix{e_1,e_2,e_3}A_1\cdot\pmatrix{e_1,e_2,e_3}A_2\cdot\pmatrix{e_1,e_2,e_3}A_3=
\pmatrix{-e_{123},e_1,e_2,e_3}\pmatrix{\delta \cr Z \cr}\ , \cr
}$$
where $A=(A_1 A_2 A_3)\ ,\delta =\det A$ and components $z_1,z_2,z_3$
of the column $Z$ are given by
$$\openup2pt\eqalign{
z_1&=-(\tr A_2Q_0A_1)a_{13} + b_{13}a_{13}\hat a_{23} + b_{23}a_{23}\hat a_{23}
   -b_{12}a_{13}\hat a_{33} - q_2a_{23}\hat a_{33} + q_3a_{33}\hat a_{23} \ ,\cr
z_2&=-(\tr A_2Q_0A_1)a_{23} - b_{13}a_{13}\hat a_{13} - b_{23}a_{23}\hat a_{13}
   + q_1a_{13}\hat a_{33} - q_3a_{33}\hat a_{13}\ , \cr
z_3&=-(\tr A_2Q_0A_1)a_{33} + b_{12}a_{13}\hat a_{13} + q_2a_{23}\hat a_{13}
   - q_1a_{13}\hat a_{23}\ , \cr
}$$
and $\hat a_{ij}$ is the algebraic cofactor of $a_{ij}$ in $A\,$.
Using (5.16) applied to $t''$ we have
$$
\varphi_A(t'')=\pmatrix{-e_{123},e_1,e_2,e_3}
\pmatrix{
\delta & 0\,0\,0 \cr
 Z & A \cr }
\pmatrix{-2 \cr B'' \cr}\quad .
$$
On the other hand, from (5.16) we also know the coordinates of $t$ in this
basis. Therefore, in order to prove the Lemma it remains to check that
$$
A\cdot\pmatrix{-b_{23}''\cr\  b_{13}''\cr -b_{12}''\cr}
-2\cdot\pmatrix{z_1\cr z_2\cr z_3\cr}
=\delta\cdot\pmatrix{-b_{23}\cr\  b_{13}\cr -b_{12}\cr}
\eqno (5.23)
$$
Before entering into computations let us make several remarks.
First, since $\tr A_2Q_0A_1+\tr A_2\tr Q_0A_1=\tr A_2QA_1=b_{12}''\,$,
then $2\cdot\tr A_2Q_0A_1-b_{12}'' = \tr A_2 (Q_0-\tr Q_0) A_1
= -b_{12}\hat a_{33} + b_{13}\hat a_{23} - b_{23}\hat a_{13}\,$.
Second, $Q[A]=Q''$ implies that ${\delta}^{-1}\cdot\tr\tilde A Q''=QA\,$.
Third, $(a_{ij})=A=(A^{-1})^{-1}=\delta^{-1}\cdot\check A
=\delta^{-1}\cdot(\check a_{ij})\,$, where $\check A$ is the double adjoint
of $A\,$. Now we proceed to calculation of entries of the column-vector on
the left hand side of (5.23).
\par
First entry:
$$
\displaylines{\qquad
-a_{11}b_{23}''+a_{12}b_{13}''-a_{13}b_{12}''-2z_1=\hfill\cr
\qquad -b_{23}(a_{13}\hat a_{13}+a_{23}\hat a_{23})-b_{13}a_{13}\hat a_{23}
+b_{12}a_{13}\hat a_{33}-b_{23}a_{23}\hat a_{23}+2q_2a_{23}\hat a_{33}
-2q_3a_{33}\hat a_{23}\hfill\cr
\hfill-(a_{11}b_{23}''-a_{12}b_{13}'')=\cr
\qquad-\delta b_{23}+\hat a_{33}(b_{12}a_{13}+2q_2a_{23}+b_{23}a_{33})
-\hat a_{23}(b_{13}a_{13}+b_{23}a_{23}+2q_3a_{33})\hfill\cr
\hfill-(a_{11}b_{23}''-a_{12}b_{13}'')=\cr
\qquad \delta^{-1}
\hat a_{33}(\hat a_{21}b_{13}''+\hat a_{22}b_{23}''+\hat a_{23}2q_3'')
-\delta^{-1}\hat a_{23}
(\hat a_{31}b_{13}''+\hat a_{32}b_{23}''+\hat a_{33}2q_3'')
-\delta b_{23}=\hfill\cr
\qquad \delta^{-1}b_{13}''(\hat a_{33}\hat a_{21}-\hat a_{23}a_{31})
+\delta^{-1}b_{23}''(\hat a_{33} \hat a_{22}-\hat a_{23}\hat a_{32})
-(a_{11}b_{23}''-a_{12}b_{13}'')-\delta b_{23}=\hfill\cr
\qquad \delta^{-1}b_{13}''(-\check a_{12})
+\delta^{-1}b_{23}''(\check a_{11})-(a_{11}b_{23}''-a_{12}b_{13}'')-\delta b_{23}=
-\det A\cdot b_{23}\,.\hfill\cr
}$$
Second entry:
$$\displaylines{\qquad
-a_{21}b_{23}''+a_{22}b_{13}''-a_{23}b_{12}''-2z_2=\hfill\cr
\qquad b_{13}(a_{13}\hat a_{13}+a_{23}\hat a_{23})
+\hat a_{13}(b_{13}a_{13}+b_{23}a_{23}+2q_3a_{33})
-\hat a_{33}(2q_1a_{13}+b_{12}a_{23})\hfill\cr
\hfill-(a_{21}b_{23}''-a_{22}b_{13}'')=\cr
\qquad \delta^{-1}\hat a_{13}(\hat a_{31}b_{13}''
+\hat a_{32}b_{23}'')
-\delta^{-1}\hat a_{33}(\hat a_{11}b_{13}''+\hat a_{12}b_{23}'')
-(a_{21}b_{23}''-a_{22}b_{13}'')+\delta b_{13}=
\det A\cdot b_{13}\,.\hfill\cr
}$$
Third entry:
$$\displaylines{\qquad
-a_{31}b_{23}''+a_{32}b_{13}''-a_{33}b_{12}''-2z_3=\hfill\cr
\qquad-b_{12}(a_{13}\hat a_{13}+a_{33}\hat a_{33})
+\hat a_{23}(2q_1a_{13}+b_{13}a_{33})
-\hat a_{11}(b_{12}a_{13}+2q_2a_{23}+b_{23}a_{33})\hfill\cr
\hfill-(a_{31}b_{23}''-a_{32}b_{13}'')=\cr
\qquad \delta^{-1}\hat a_{23}(\hat a_{11}b_{13}''
+\hat a_{12}b_{23}'')
-\delta^{-1}\hat a_{13}(\hat a_{21}b_{13}''+\hat a_{22}b_{23}'')
-(a_{31}b_{23}''-a_{32}b_{13}'')-\delta b_{12}=
-\det A\cdot b_{12}\,.\hfill\cr
}$$
\flushqed
The above lemma along with (5.16) leads us to the following criterion:
\Theorem{5.2}
{\sl Let $\varphi:C_0(E'')\to C_0(E)$ be a homomorphism of
even subalgebras of Clifford algebras associated to ternary quadratic
$\Z$-modules $(E'',\q'')$ and $(E,\q)\,$. Set $\Delta=\det\q/2\, ,\,
\Delta''=\det\q''/2$ and let $t\in C(E)\, ,\, t''\in C(E'')$ be special
elements of the form (5.15). Then $\varphi$ coincides with the lift (5.20)
of an isometry $\alpha:(E'',\q'')\to (E,\q)$ if and only if
$$
\varphi(E''\cdot t'')\subset\sqrt{\Delta''/\Delta}\cdot (E\cdot t)\ .
\eqno (5.24)
$$
In this case the isometry $\alpha$ is unique (up to sign) and is given by
$$
\alpha(x'')={{\pm\varphi(x''t'')\cdot t}\over{\sqrt{\Delta\Delta''}}}
\quad\hbox{ for }\  x''\in E''\ ,
\eqno (5.25)
$$
where one can fix either ``$+$'' or ``$-$''.}
\Proof Assume (5.24) and consider $\alpha$ as defined in (5.25). For an
$x''\in E''$ we have:
$$
\alpha(x'')\subset (\sqrt{\Delta\Delta''})^{-1}\cdot\varphi(E''t'')\cdot t
\subset (\sqrt{\Delta\Delta''})^{-1}\sqrt{\Delta''/\Delta}\cdot E\cdot t^2
\subset -|\Delta|^{-1}\Delta\cdot E\subset E
$$
Furthermore,
$$\displaylines{
\q\bigl(\alpha(x'')\bigr)=\alpha(x'')\cdot\alpha(x'')=
(\Delta\Delta'')^{-1}\cdot\varphi(x''t'')\cdot t\cdot\varphi(x''t'')\cdot t=
(\Delta\Delta'')^{-1}\cdot\varphi(x''t''x''t'')\cdot t^2=\hfill\cr
\hfill(\Delta\Delta'')^{-1}\cdot\varphi\bigl(-\Delta''(x'')^2\bigr)
\cdot(-\Delta)=\varphi\bigl(\q''(x'')\bigr)=\q''(x'')\ ,
\cr }$$
thus $\alpha:(E'',\q'')\to(E,\q)$ is indeed an isometry. According to (5.20)
it can be lifted to an algebraic homomorphism
$\varphi_{\alpha}:C(E'')\to C(E)$ the restriction of which on $C_0(E'')$
is homomorphism into $C(E)\,$. Then
$$
\varphi_{\alpha}(e_{ij}'')=\alpha(e_i'')\cdot\alpha(e_j'')=
(\Delta\Delta'')^{-1}\cdot\varphi(e_i''t''e_j''t'')\cdot t^2=
(\Delta\Delta'')^{-1}(-\Delta)\cdot\varphi(-\Delta''e_i''e_j'')=
\varphi(e_{ij}'')
$$
so $\varphi_{\alpha}$ and $\varphi$ coincide on $C_0(E'')\,$.
\par
Conversely, if $\varphi=\varphi_{\alpha}\bigm|_{C_0}$ for some lift $\varphi_{\alpha}$
of an isometry $\alpha:(E'',\q'')\to(E,\q)\,$, then
$$
\varphi(x''t'')=\varphi_{\alpha}(x''t'')=\varphi_{\alpha}(x'')\cdot
\det\alpha\cdot t=\alpha(x'')\cdot\det\alpha\cdot t\subset
\det\alpha\cdot E\cdot t\subset\sqrt{\Delta''/\Delta}\cdot(E\cdot t)\ ,
$$
so (5.24) is satisfied. The equality (5.25) is also satisfied for one of
the choices of signs.
\flushqed
To summarize what we have established so far, the universal property (5.6)
of Clifford algebras allows us to lift an integral automorph
$A\in R(\q,\q'')$ of ternary forms to an integral automorph
$\Phi_A\in R(\n,\n'')$ of quaternary forms associated with norms on the
even subalgebras of Clifford algebras. Conversely, the original automorph
$A\in R(\q,\q'')$ can be recovered from the lift $\Phi_A\in R(\n,\n'')$
with the help of the Theorem 5.2 above. (We note again that quaternary
quadratic forms $\n\,,\,\n''$ are specializations of norm-forms (5.10) at
the natural bases (5.7).)
\par
Recall that the quadratic forms we are working with are not necessarily
positive definite and so the sets $R(\q,\q'')\,,\,R(\n,\n'')\dots$ of
automorphs or the sets $R(\q,a)\,,\,R(\n,a)\dots$ of representations of
a number $a$ can be infinite in general. But quotients modulo corresponding
groups of units $\bigl(E(\q)\,,\,E(\n)\dots\bigr)$ are always finite for
nonsingular forms. Because of this it is more convenient sometimes to
consider the orbits modulo groups of units rather than individual automorphs
or representations (compare to section 3). The following lemma shows that
this transition is compatible with lift (5.22):
\Lemma{5.3}{\sl Let $\ A,D\in R(\q,\q'')\ $ be integral automorphs. Then
$\ A\in E(\q)D\ $ if and only if $\ \Phi_A\in E(\n)\Phi_D\ $.}
\Proof Let $E\in E(\q)$ be such that $A=ED\,$. Consider the following
commutative diagram:
$$\diagram{
(E,\q) & \lharr{\alpha_E}{} & (E,\q) & \lharr{\alpha_D}{} & (E'',\q'') \cr
\dvarr{}{\iota} & & \dvarr{}{\iota} & & \dvarr{}{\iota''} \cr
C(E) & \lharr{\varphi_E}{} & C(E) & \lharr{\varphi_D}{} & C(E'') \cr
}$$
where the isometries $\alpha_E,\alpha_D$ and algebraic homomorphisms
$\varphi_E,\varphi_D$ are defined analogously to (5.20). Since
$\alpha_E\circ\alpha_D=\alpha_A\,$, we know that there exists unique
algebraic homomorphism $\varphi_A=\varphi_{ED}:C(E'')\to C(E)$ such that
$\iota\circ\alpha_A=\varphi_A\circ\iota''$ as in (5.20). But
$\iota\circ\alpha_A=\iota\circ\alpha_E\circ\alpha_D=
\varphi_E\circ\iota\circ\alpha_D=\varphi_E\circ\varphi_D\circ\iota''\,$,
which implies that $\varphi_{ED}=\varphi_E\circ\varphi_D$ and thus
$\Phi_A=\Phi_E\Phi_D$ (see (5.22)). Note that since $E\in E(\q)$ then
$\Phi_E\in E(\n)\,$.
\par
Conversely, assume that $\Phi_A\in E(\n)\Phi_D\,$. Then 
$\Phi_A\Phi_D^{-1}\in E(\n)$ and in particular
$$
\pmatrix{
1 & * & * & * \cr
0 &   &   &   \cr
0 & & \tr\tilde A & \cr
0 &   &   &   \cr
}
\pmatrix{
1 & * & * & * \cr
0 &   &   &   \cr
0 & & \tr\tilde D^{-1} & \cr
0 &   &   &   \cr
}
=\pmatrix{
1 &  & *\ *\ * &  \cr
0 &   &   &   \cr
0 & &{1\over{\det D}} \tr\tilde A\cdot\tr D & \cr
0 &   &   &   \cr
}$$
is an integral matrix. But $\det D=\det A\,$, so $\tr A^{-1}\cdot\tr D$ is
also an integral matrix. Then $E=DA^{-1}$ is also integral and
$Q[E]=Q[DA^{-1}]=Q''[A^{-1}]=Q\,$, i.e. $E\in E(\q)$ and $A\in E(\q)D\,$.
\flushqed
\Remark{}{\sl An obvious modification of the above proof shows that}
$$
A\in D\cdot E(\q'') \iff \Phi_A\in\Phi_D\cdot E(\n'')\quad.
$$
\medskip
\noindent
We will also need the following technical result:
\Lemma{5.4}{\sl Assume that (ternary) quadratic forms $\q$ and $\q''$ 
belong to the same similarity class. Then $\q$ and $\q''$ are integrally
equivalent if and only if corresponding quadratic forms $\n$ and $\n''$
(see (5.12)) are integrally equivalent.}
\Proof Suppose $\q$ and $\q''$ belong to the same equivalence class, i.e.
there exists a matrix $U\in R(\q,\q'')\cap\Lambda^3\,$. Then
$\alpha_U:(E'',\q'')\to(E,\q)$ is an isometric isomorphism of quadratic
modules which gives rise to an isomorphism $\varphi_U:C(E'')\to C(E)$ of
Clifford algebras exactly as in (5.20). The restriction of $\varphi_U$ on
$C_0(E'')$ is an isomorphism of even subalgebras whose matrix with respect
to the natural bases (5.7) is $\Phi_U\in R(\n,\n'')\cap\Lambda^4\,$
(see (5.22)). Thus $\n$ and $\n''$ are integrally equivalent.
\par
Conversely, suppose that there exists
${\cal U}\in R(\n,\n'')\cap\Lambda^4\,$. The matrix $\cal U$ defines a
(linear) isometric isomorphism of quadratic modules:
$$\openup2pt \displaylines{
\upsilon:\bigl(C_0(E''),\n''\bigr)\longrightarrow\bigl(C_0(E),\n\bigr)\ ,\cr
\bigl(\upsilon(e_0''),\upsilon(e_{23}''),\upsilon(-e_{13}''),
\upsilon(e_{12}'')\bigr)=
\bigl(e_0,e_{23},-e_{13},e_{12}\bigr)\cdot{\cal U}\quad .\cr
}$$
Since $\n(\upsilon(e_0''))=\n''(e_0'')=1\,$, the element $\upsilon(e_0'')$
has an inverse $\upsilon(e_0'')^{-1}=\overline{\upsilon(e_0'')}\in C_0(E)$
whose norm is also $1$. Then the right multiplication 
$x\mapsto x\cdot\overline{\upsilon(e_0'')}$ is a (linear) isometric
automorphism of $C_0(E)\,$. The composition of these two isometries
$$
\omega:x''\mapsto\upsilon(x'')\cdot\overline{\upsilon(e_0'')}\, \qquad
x''\in C_0(E'')
$$
gives a (linear) isometric isomorphism of quadratic modules
$\bigl(C_0(E''),\n''\bigr)\to\bigl(C_0(E),\n\bigr)$ which maps $e_0''$
to $e_0\,$. Therefore its matrix with respect to natural bases (5.7) 
has the form
$$
{\cal W}=\pmatrix{
1 & & *** & \cr
0 & &     & \cr
0 & & W   & \cr
0 & &     & \cr
}\quad\hbox{ with } W\in\Lambda^3\ .
$$
Further, since $N[{\cal W}]=N''\,$, then
$N^{-1}[\tr{\cal W}^{-1}]=(N'')^{-1}$ with $N^{-1}\,,\,(N'')^{-1}$ given
by (5.14). By our assumption the forms $\q$ and $\q''$ are similar, in
particular $\Delta=\det\q/2=\Delta''\,$. Therefore (5.14) implies 
$$
\pmatrix{
2\,\det{Q_0} & \tr(Q_0 B) \cr
     &   \cr
(Q_0 B) & Q \cr
}
\Biggl[
\pmatrix{
1 & & 0\ 0\ 0 & \cr
* & &  & \cr
* & & \tr W^{-1}   & \cr
* & &  & \cr
}\Biggr] =
\pmatrix{
2\,\det{Q_0''} & \tr(Q_0'' B'') \cr
     &   \cr
(Q_0'' B'') & Q'' \cr
}\ .
$$
It follows that $Q[\tr W^{-1}]=Q''$ with $W\in\Lambda^{3}\,$, i.e. the
quadratic forms $\q$ and $\q''$ are (integrally) equivalent.
\flushqed
We now proceed directly to an automorph analog of Shimura's lift. Let
$\q'$ be a quadratic form {\it similar} to $\q$ (see section 1) and let
$p$ be a prime coprime to $\det\q=\det\q'\,$. Given an integral automorph
$A\in R(\q,p^2\q')$ we can set $\q''=p^2\q'$ and apply the construction
(5.20) to get the homomorphism $\varphi_A$ of Clifford algebras:
$$\diagram{
(E,\q) & \lharr{\alpha_A}{} & (E'',p^2\q') \cr
\dvarr{\iota}{} &  & \dvarr{}{\iota''} \cr
C(E,\q) & \lharr{\varphi_A}{} & C(E'',p^2\q') \cr
}
\eqno (5.26)
$$
The restriction of $\varphi_A$ to the even subalgebra $C_0(E'',p^2\q')$
will give us an automorph $\Phi_A\in R(\n,\n'')=R(N,N'')$ given by (5.22).
Here $N$ is the matrix (5.12) and
$$
N''=
\pmatrix{
2 & -p^2\cdot\tr B' \cr
  &  \cr
-p^2 B' & p^4(\tr\tilde Q_0'+\tilde Q_0') \cr }
=N'
\biggl[
\pmatrix{
1 & \cr
 & p^2 1_3 \cr}
\biggr]
\eqno (5.27)
$$
is the matrix (5.12) of the norm-form (5.10) of the even subalgebra
$C_0(E'',p^2\q')$ with respect to its natural basis (5.7). Note that
$C_0(E'',p^2\q')$ can be identified with the subring of $C_0(E',\q')$
spanned by
$e_0',p^2 e_{23}',-p^2e_{13}',p^2e_{12}'\,$. (Indeed, the natural isometry
of quadratic modules $E''\cong pE'$ defined by $e_i''\mapsto pe_i'$ gives
rise to an algebraic monomorphism $C(E'',p^2\q')\hookrightarrow C(E',\q')\,$,
the restriction of which to the even subalgebra is the desired ring
homomorphism.) In particular, $N'$
 on the right hand side of (5.27) is
exactly the matrix (5.12) of the norm-form $\n'$ on the subalgebra 
$C_0(E',\q')\,$. Combining the last equality with
$N[\Phi_A]=N''$ we conclude that if $A\in R(\q,p^2\q')$ then
$N[\Psi_A]=N'[(p\cdot1_4)]=p^2N'$,
where we set
$$
\Psi_A=
\Phi_A\cdot
\pmatrix{ p & \cr
 & p^{-1}1_3}=
\pmatrix{
p & p^{-1}(\ \tr A_3 Q_0 A_2 & -\tr A_3 Q_0 A_1 & \tr A_2 Q_0 A_1\ ) \cr
0 & & &  \cr
0 & & p^2\cdot\tr A^{-1} & \cr
0 & & &  \cr
}\quad .
\eqno (5.28)
$$
\Lemma{5.5}{\sl In the above notations $\Psi_A\in R(\n,p^2\n')\,$, in
particular the matrix $\Psi_A$ is integral.}
\Proof We already know that $N[\Psi_A]=p^2N'\,$ therefore to prove the
lemma it remains to show that $\Psi_A$ is integral. Since $A\in R(\q,p^2\q')$
then $Q[A]=p^2Q'$ and so $p^2 A^{-1}=(\det\q)^{-1}\tilde Q'\tr A Q$ which
implies that $\det\q\cdot p^2 A^{-1}$ is an integral matrix. But 
$\det A=p^3$ is coprime to $\det\q\,$, thus $p^2 A^{-1}$ is integral. 
Next we need to deal with the first row of $\Psi_A$ in (5.28). Denote
\smallskip
$$
\tr Z_A=
\Bigl(
\tr A_3 Q_0 A_2 , -\tr A_3 Q_0 A_1 , \tr A_2 Q_0 A_1
\Bigr)
\quad .
\eqno (5.29)
$$
Then
\smallskip
$$
N[\Phi_A]=
\pmatrix{
2 & & 2\cdot\tr Z_A - p^3\cdot\tr(A^{-1}B) \cr
* & & \cr
* & & 2^{-1}\bigl( p^4\tilde Q'+(2Z_A-p^3 A^{-1}B)\tr(2Z_A-p^3 A^{-1}B)\bigr)\cr
}=
N'\biggl[\pmatrix{
1 & \cr
 & p^2 1_3 \cr}\biggr]
$$
\smallskip
which implies that $2Z_A-p^3 A^{-1}B=\congr{-p^2 B'}{0}{p^2}\,$. But
$p^2 A^{-1}$ is integral, so $\congr{p^3 A^{-1}B}{0}{p}$ and therefore
$\congr{Z_A}{0}{p}\,$. Thus $\Psi_A$ is indeed an integral matrix. \flushqed
We conclude that the correspondence $A\mapsto\Psi_A$ gives us an injection
$R(\q,p^2\q')\hookrightarrow R(\n,p^2\n')\,$. Because of the Lemma 5.3 we
can also view this map as injection
$$
\Psi\,:
\ E(\q)\backslash R(\q,p^2\q')\lhook\joinrel\longrightarrow
E(\n)\backslash R(\n,p^2\n')\quad ,\quad
p\not|\det\q\quad .
\eqno (5.30)
$$
\Remark{5.6} The above construction does not depend on primality of $p$
and works equally well for any number $a$ coprime to $\det\q$ resulting
in injection $R(\q,a^2\q')\lhook\joinrel\longrightarrow R(\n,a^2\n')\,$.
 In particular, if $\q$ is similar to
$\q'$ then $\n$ is similar to $\n'\,$.
\Remark{5.7} The lift $\Psi$ is also compatible with the bijection
$\ R(\q,p^2\q')\longrightarrow R(\q',p^2\q)$ given by
$A\mapsto p^2A^{-1}$ in the
sense that $\Psi_{p^2A^{-1}}=p^2\Psi_A^{-1}\in R(\n',p^2\n)\,$. Indeed,
the fact that $\Psi_{p^2A^{-1}}$ and $p^2\Psi_A^{-1}$ both belong to
$R(\n',p^2\n)$ follows directly from definition of $\Psi\,$. Then using
the equality $N'[\Psi_{p^2A^{-1}}]=N'[p^2\Psi_A^{-1}]$ together with 
(5.12) and (5.28) one can see that $\Psi_{p^2A^{-1}}=p^2\Psi_A^{-1}\,$.
\vfill\break
%
%
\bigskip\noindent
{\bf 6. Factorization of automorph lifts}
\smallskip\noindent
Thus, we can lift an automorph $A\in R(\q,p^2\q')$ of ternary quadratic
forms to an automorph $\Psi_A\in R(\n,p^2\n')$ of quaternary forms (see 5.28).
But with quaternary quadratic forms we are no longer limited to multipliers
$p^2$ and can use automorphs with multipliers $p$ (compare to existence
of Hecke operators $T(p)$ rather than $T(p^2)$ for corresponding theta-series).
According to the Theorem 1.3 of [2] each integral automorph
${\cal A}\in R(\n,p^2\n')\ ,\ p\not|\det\n$ is a product ${\cal ML}$ 
of automorphs 
${\cal M}\in R(\n,p{\bf f})$ and ${\cal L}\in R({\bf f},p\n')$ for some
quadratic form ${\bf f}$ similar to $\n\,$. Our next goal is to describe
explicitly such factorizations for the constructed lift $\Psi_A\in R(\n,p^2\n')\,$.
Since the general strategy behind Clifford algebras is to replace questions
about representations by quadratic forms with questions concerning
multiplicative arithmetic of these algebras, we will try to use for our
purposes certain properties of ideals of the algebras. First of all let us
look more closely at the algebraic homomorphism $\varphi_A$ in (5.26) the
restriction of which to the even subalgebra $C_0(E'',p^2\q')$ defines the
lift $\Psi_A$ via (5.28). Recall that the restricted $\varphi_A$ can be
considered as a ring homomorphism
$$\openup2pt \displaylines{
\varphi_A:(e_0',p^2e_{23}',-p^2e_{13}',p^2e_{12}')\Z^4
\longrightarrow C_0(E,\q)\ ,\cr
\bigl(\varphi_A(e_0'),\varphi_A(p^2e_{23}'),\varphi_A(-p^2e_{13}'),
\varphi_A(p^2e_{12}')\bigr)=
\bigl(e_0,e_{23},-e_{13},e_{12}\bigr)\cdot\Phi_A\quad ,\cr
}$$
defined on the corresponding order of $C_0(E',\q')\,$. One can easily extend
this homomorphism to an isomorphism of $\Q$-algebras:
$$
\varphi_A:C_0(E',\q')\otimes\Q\longrightarrow C_0(E,\q)\otimes\Q\ ,
$$
whose matrix with respect to the natural bases is
$\Phi_A\cdot\mathop{\rm diag}(1,p^{-2},p^{-2},p^{-2})=p^{-1}\Psi_A\,$.
Using Lemma 5.5 one can see that the matrix
$\Phi_A\cdot\mathop{\rm diag}(1,p^{-1},p^{-1},p^{-1})$ is integral, which
implies that the image under $\varphi_A$ of the order of $C_0(E')$ generated
by $e_0,pe_{23}',-pe_{13}',pe_{12}'$ belongs to $C_0(E)\,$. Then from (5.28)
it follows that the submodule $(e_0,e_{23},-e_{13},e_{12})\Psi_A\cdot\Z^4
\subset C_0(E)$ is nothing else but the image $\varphi_A\bigl(p\,C_0(E')\bigr)$
of the principle ideal $p\,C_0(E')\subset C_0(E')$ and therefore is itself an
ideal of the of the order $\varphi_A\bigl(C_0(E')\bigr)\cap C_0(E)\,$.
Consider the left $C_0(E)$-ideal
$$
{\frak I}_A=C_0(E)\cdot\varphi_A\bigl(p\,C_0(E')\bigr)\ ,
\eqno (6.1)
$$
which extends $\varphi_A\bigl(p\,C_0(E')\bigr)\,$. We claim the following
\Lemma{6.1}{\sl In the above notations ${\frak I}_A$ is a proper (left)
$C_0(E)$-ideal which properly includes $\varphi_A\bigl(p\,C_0(E')\bigr)\,$:}
$$
\varphi_A\bigl(p\,C_0(E')\bigr)\subset{\frak I}_A\subset C_0(E)\ .
$$
{\sl Moreover, the norm of any element of ${\frak I}_A$ is divisible by $p\,$.}
\Proof
Validity of the above inclusions is clear and we only need to prove that both
of them are proper. To show that $\varphi_A\bigl(p\,C_0(E')\bigr)\not=
{\frak I}_A\,$, assume the opposite. Since ${\frak I}_A$ is a left 
$C_0(E)$-ideal, we would have then for any $c\in C_0(E)$ an inclusion 
$c\cdot\varphi_A\bigl(p\,C_0(E')\bigr)\subset\varphi_A\bigl(p\,C_0(E')\bigr)\,$.
Recall that as an isomorphism of $\Q$-algebras $\varphi_A$ is invertible,
the matrix of $\varphi_A^{-1}$ with respect to the natural bases is 
$p\Psi_A^{-1}=p^{-1}\Psi_{p^2 A^{-1}}\,$, in particular
$\varphi_A^{-1}\bigl(C_0(E)\bigr)\subset p^{-1}C_0(E')\,$. Thus, our
assumption would imply that $\varphi_A\bigl(\varphi_A^{-1}(c)\cdot
p\,C_0(E')\bigr)\subset\varphi_A\bigl(p\,C_0(E')\bigr)$ and so
$\varphi_A^{-1}(c)\cdot p\,C_0(E')\subset p\,C_0(E')$ for any $c\in C_0(E)\,$.
The last inclusion would be possible only if $\varphi_A^{-1}(c)\in C_0(E')$
for any $c\in C_0(E)\,$, which is definitely not the case (for example for
$c=e_{23}\,$). Therefore $\varphi_A\bigl(p\,C_0(E')\bigr)$ is a proper subset
of ${\frak I}_A\,$. Next, to prove that ${\frak I}_A\not=C_0(E)$ we will show
that the norm of any element of ${\frak I}_A$ is divisible by $p\,$. Indeed, any
element of ${\frak I}_A=C_0(E)\cdot\varphi_A\bigl(p\,C_0(E')\bigl)$ is a sum of
elements of the form $c\cdot\varphi_A(pc')$ with $c\in C_0(E)$ and
$c'\in C_0(E')\,$. Observe that $\congr{\n\bigl(\varphi_A(pc')\bigr)}{0}{p^2}$
for any $c'\in C_0(E')\,$, according to Lemma 5.5\ . It remains to note that if
$b,c\in C_0(E)$ and $b',c'\in C_0(E')$ then
$$\openup2pt \displaylines{
\n\bigl(b\,\varphi_A(pb')+c\,\varphi_A(pc')\bigr)=
\n\bigl(b\,\varphi_A(pb')\bigr)+\n\bigl(c\,\varphi_A(pc')\bigr)
+\s\bigl(b\,\varphi_A(pb')\cdot\overline{c\,\varphi_A(pc')}\bigr)\ \equiv \cr
\s\bigl(b\,\varphi_A(pb')\cdot\varphi_A(p\,\overline {c'})\,\bar c\bigr)
\,\equiv\,\s\bigl(pb\,\varphi_A(pb'\,\overline{c'})\,\bar c\bigr)\,\equiv\,
\congr{p\cdot\s\bigl(b\,\varphi_A(pb'\,\overline{c'})\,\bar c\bigr)}{0}{p} \cr
}$$
therefore $\congr{\n\bigm|_{{\frak I}_A}}{0}{p}$ as claimed and so ${\frak I}_A$ is a
proper subset of $C_0(E)\,$.
\flushqed
The above lemma immediately leads to
\Lemma{6.2}
{\sl Let ${\cal I}_A$ be a generating matrix of the extension ideal ${\frak I}_A$
(6.1) with respect to the natural basis, so that}
${\frak I}_A=(e_0,e_{23},-e_{13},e_{12})\cdot{\cal I}_A\Z^4\ .$
{\sl Then $\Psi_A\in{\cal I}_A\Z^4_4$ and ${\cal I}_A\in R(\n,p{\bf f})$ with
integral quadratic form
$$
{\bf f}_{{\frak I}_A}={{\n\bigm|_{{\frak I}_A}}\over{\sqrt{|{\frak I}_A|}}}
$$
{\it similar} to $\n\,$. Here we have set $|{\frak I}_A|=|{\cal I}_A|=
|C_0(E)\slash{\frak I}_A|\,$, the index of ${\frak I}_A$ in $C_0(E)\,$.}
\Proof
Indeed, since ${\frak I}_A$ is a proper (left) ideal of $C_0(E)\,$, we know
that $|{\cal I}_A|\not=1\,$. Moreover, for a fixed element
${\frak i}\in{\frak I}_A$ we have $C_0(E)\cdot{\frak i}\subset{\frak I}_A$ and
thus
$C_0(E)\cdot{\frak i}=
(e_0{\frak i},e_{23}{\frak i},-e_{13}{\frak i},e_{12}{\frak i})\,\Z^4=
(e_0,e_{23},-e_{13},e_{12})\cdot{\cal I}_A{\cal Y}\,\Z^4$
with some ${\cal Y}\in\Z^4_4\,$. Then, looking at the matrix of the restriction
of the norm-form $\n$ onto the submodule $C_0(E)\cdot{\frak i}\subset C_0(E)$ we
see that $\n({\frak i}) N=N[{\cal I}_A{\cal Y}]$ which implies that
$|{\frak I}_A|$ divides $\n^2({\frak i})\,$. Denote for a moment 
$|{\cal I}_A|=a^2b$ with a square-free positive integer $b\,$. Clearly then
$ab$ divides $\n({\frak i})\,$, but ${\frak i}\in{\frak I}_A$ was arbitrary and 
therefore the form ${1\over{ab}}\n|_{{\frak I}_A}$ is integral as well as its
matrix $(ab)^{-1}N[{\cal I}_A]$ with respect to the natural basis. In particular
$\det\left((ab)^{-1} N[{\cal I}_A]\right)=(|{\frak I}_A|^2\det{\n})\slash(ab)^4=
(\det{\n})\slash b^2=
(\Delta\slash b)^2$ is integral and so $b$ divides $\Delta=\det{\q}\slash 2\,$,
see (5.13). On the other hand, using Lemma 6.1 we also see that 
$\varphi_A\bigl(p\,C_0(E')\bigr)$ is a proper subset of ${\frak I}_A$ which in
terms of corresponding matrices means that $\Psi_A={\cal I}_A{\cal X}$ for some
${\cal X}\in\Z^4_4\ ,\ |{\cal X}|\not=1$ and in particular $|{\frak I}_A|$ is a
proper divisor of $|\Psi_A|=p^4\,$. Therefore $b\mid\gcd(\Delta,p^4)\,$, but
according to our assumption prime $p$ is coprime to $\det\q=2\Delta$ thus
$b=1$ and $|{\frak I}_A|=p^2\,$. We conclude that the quadratic form
$$
{\bf f}_{{\frak I}_A}(x)=
{{\n\bigm |_{{\frak I}_A}(x)}\over{\sqrt{|{\frak I}_A|}}}=
p^{-1}\cdot\n\bigm |_{{\frak I}_A}(x)
\eqno (6.2)
$$
is integral. The matrix of ${\bf f}_{{\frak I}_A}$ with respect to the
natural basis $(e_0,e_{23},-e_{13},e_{12})\cdot{\cal I}_A$ is 
$F_{{\frak I}_A}=p^{-1}N[{\cal I}_A]\in\Z^4_4$
and its determinant $\det{\bf f}_{{\frak I}_A}$ is equal to $\det\n\,$.
Summing up we see that ${\bf f}_{{\frak I}_A}$ is {\it similar} to $\n$
and that ${\cal I}_A\in R(\n,p\,{\bf f}_{{\frak I}_A})\,$.
\flushqed
Lemmas 6.1 and 6.2 in essence allow us to write a lift $\Psi_A\in R(\n,p^2\n')$
as a product of automorphs with multipliers $p$ of {\it similar} quadratic forms
$\n\,,\,\n'$ and ${\bf f}\,$. Note that the quadratic form 
${\bf f}={\bf f}_{{\frak I}_A}$ in the above factorization is associated with a
class of equivalent (left) ideals of $C_0(E,\q)$ and is defined only up to
integral equivalence, see (6.2). (Two left $C_0(E)$-ideals 
${\frak I,J}\subset C_0(E)$ are called {\it equivalent} if ${\frak I}\alpha=
{\frak J}\beta$ for some $\alpha,\beta\in C_0(E)\,$.) An analogous
factorization of $\Psi_A$ can be obtained if one considers the extension of
$\varphi_A\bigl(p\,C_0(E')\bigr)$ to the right ideal 
$\varphi_A\bigl(p\,C_0(E')\bigr)\cdot C_0(E)$ instead of the left ideal as
in (6.1). We now seek to find the total number of such factorizations (different
modulo corresponding groups of units) for a given $\Psi_A\,$. Following a
general method of investigation of questions of this type developed in [1],
we will view a matrix ${\cal M}\in R(\n,p\,{\bf f})$ for a quadratic form 
${\bf f}$ {\it similar} to $\n$ as a solution of the quadratic congruence
$$
\congr{\n[{\cal M}]}{0}{p}\ ,
$$
whose matrix of {\it elementary divisors} is equal to 
${\cal D}_p=\diag(1,1,p,p)\,$. That is to say, ${\cal M}$ belongs to the
double coset $\Lambda^4{\cal D}_p\Lambda^4$ which is a necessary condition
for ${\cal M}$ to define an automorph with multiplier $p$ of similar
quaternary quadratic forms. (Indeed, we should have $|{\cal M}|=p^2$ and
$p{\cal M}^{-1}\in\Z^4_4\,$, which leaves us with ${\cal D}_p$ as the only
possible matrix of elementary divisors.) Since at the moment we do not
want to distinguish a particular form ${\bf f}$ in its equivalent class
$\{{\bf f}[{\cal U}]\ ;\ {\cal U}\in\Lambda^4\}\,$, we will next be counting
only different left cosets 
${\cal M}\Lambda^4\subset\Lambda^4{\cal D}_p\Lambda^4\slash\Lambda^4\,$.
Finally, because of our interest in factorizations of automorphs with
multipliers $p^2\,$, we restrict our search to only those automorphs
${\cal M}\,$, which divide (from the left) a particular 
${\cal A}\in R(\n,p^2\n')\,$. Thus, following [1], in order to find the total
number of factorizations of ${\cal A}$ one introduces {\it isotropic sums}
$$
{\cal S}_p(\n,{\cal D}_p,{\cal A})=
\sum_{\textstyle{{{\cal M}\in\Lambda^4{\cal D}_p\Lambda^4\slash\Lambda^4}\atop
{\congr{\n[{\cal M}]}{0}{p}\,,\,{\cal M}|{\cal A}}}} 1 \qquad .
\eqno (6.3)
$$
General isotropic sums of the above type were computed in [1], Theorem 5.1
(see also [2] and [4]). The idea behind the computation is to associate to
a matrix ${\cal L}\in\Z^4_{*}$ the subspace $V_p({\cal L})$
of the quadratic space $\left((\Z\slash p\Z)^4,\n\mod p\right)$ spanned by
the columns of ${\cal L}$ modulo $p\,$. Then the summation conditions on
${\cal M}$ in (6.3) mean exactly that $V_p({\cal M})$ is a $2$-dimensional
isotropic subspace of $\left((\Z\slash p\Z)^4,\n\right)$ which contains 
fixed subspace $V_p({\cal A})\,$. The number of such subspaces can be found
with help of standard methods of Geometric Algebra over finite fields
(see [9]). Specializing results of [1], Theorem 5.1 to our case we find that
$$
{\cal S}_p(\n,{\cal D}_p,{\cal A})\quad =\qquad
\cases{
1 & if $\ \mathop{\rm rank}_{\F_p}{\cal A}=2,$ \cr
1+\xa{\n} & if $\ \mathop{\rm rank}_{\F_p}{\cal A}=1,\qquad$ \cr
(1+\xa{\n})(p+1) & if $\ \mathop{\rm rank}_{\F_p}{\cal A}=0,$ \cr
}
\eqno (6.4)
$$
where $p{\not|}\det\q$ and $\xa{\n}$ is the character (2.4) of the form $\n\,$.
(Note that $\xa{\n}=1$ for $p{\not|}\det\q$ because of (5.13).) The above sum
can be rewritten as follows. The condition $\congr{N[{\cal M}]}{0}{p}$ means
that $N[{\cal M}]=p\,F$ for some even matrix $F$ of oder $4\,$. Since
$\det{\cal M}=\det{\cal D}_p=p^2\,$, then
$\det F=\det N$ and thus quadratic form ${\bf f}(X)=2^{-1}F[X]$ is 
{\it similar} to $\n\,$. Since ${\cal M}$ is defined only modulo right
multiplication by $\Lambda^4\,$, we can replace $F$ by any (integrally)
equivalent form $F[{\cal U}]\,,\ {\cal U}\in\Lambda^4\,$. Let us fix a
complete system  $\{{\n}_1,\dots,{\n}_H\}$ of representatives of
different equivalence classes of the similarity class of $\n$ (analogous
to (4.2)). Note that $h\leq H$ since we can choose $\n_i$ with $1\leq i\leq h$
to be quadratic forms defined by the norms on the even Clifford subalgebras
$C_0(\Z^3,\q_i),\ 1\leq i\leq h\,$. (Indeed, according to Lemmas 5.4 and 5.5
these $n_i$'s $\,1\leq i\leq h$ belong to the different equivalence classes of
the same similarity class.)
Then for each
${\cal M}\in\Lambda^4{\cal D}_p\Lambda^4\slash\Lambda^4$ with 
$\congr{N[{\cal M}]}{0}{p}$ there exists a unique number 
$i\,,\ 1\leq i\leq H$ such that $\n[{\cal MU}]=p\,{\n}_i$
for some ${\cal U}\in\Lambda^4\,$. Next we note that according to the theory
of elementary divisors (see [3], Lemma 3.2.2) 
$R(\n,p\,{\n}_i)\subset\Lambda^4{\cal D}_p\Lambda^4\,$, therefore
$$
{\cal S}_{p}(\n,{\cal D}_p,{\cal A})\ =\quad 
\sum\limits_{i=1}^H
\sum_{\textstyle{
{{\cal M}\in R(\n,p\,{\n}_i)\slash E({\n}_i)}\atop{{\cal M}|{\cal A}}
}} 1 \qquad ,
\eqno (6.5)
$$
which already resembles the right hand side of hypothetical (4.5) we
seek to establish. Similar considerations of another isotropic sum allow
us to compute also the left hand side of (4.5). Indeed, let
$\{\q_1,\dots,\q_h\}$ be a complete system (4.2). According to the theory 
of elementary divisors the set of primitive automorphs
$R^*(\q,p^2\q_i)$ coincides with $R(\q,p^2\q_i)\cap\Lambda^3
D_p\Lambda^3\,$, where $D_p=\diag(1,p,p^2)\,$, see [4] Lemma 3.1.
Conversely, any $A\in\Lambda^3 D_p\Lambda^3\slash\Lambda^3$ such that
$\congr{Q[A]}{0}{p^2}$ defines a quadratic form $p^{-2}\,\q[A]$ which is
{\it similar} to $\q$ and therefore is integrally equivalent to exactly
one of the forms $\q_i$ so that $\q[AU]=p^2\q_i$ for a unique number 
$i\ ,\ 1\leq i\leq h\,$ and some $U\in\Lambda^3\,$. We conclude that
$$
\sum\limits_{i=1}^h
\sum_{A\in R^*(\q,p^2\q_i)\slash E(\q_i)} 1 
\quad = \quad
{\cal S}_{p^2}(\q,D_p,O_3)
\quad = \quad
\sum_{
\textstyle{
{A\in\Lambda^3 D_p\Lambda^3\slash\Lambda^3}\atop{\congr{Q[A]}{0}{p^2}}}} 1
\qquad ,
\eqno (6.6)
$$
where $O_3$ is the zero matrix of order 3. Furthermore, using Theorem 2.2 of
[4] we can provide an explicit value of the isotropic sum (6.6):
$$
{\cal S}_{p^2}(\q,D_p,O_3)\ =\ p+1 \quad .
\eqno (6.7)
$$
Combining formulas (6.3)--(6.7) we finally establish (4.5) and deduce the
following Theorem, which expresses total number of ternary automorphs of
$\q$ with multiplier $p^2$ in terms of quaternary automorphs with multiplier
$p\,$ and can be viewed as a generalization of the Shimura's correspondence
for theta-series of ternary positive definite quadratic forms to the case
of indeterminate forms:
\Theorem{6.3}
{\sl Let $\q$ be an integral nonsingular ternary quadratic form
(5.1) with determinant $\det\q=2\Delta\,$, and let $p$ be a prime not
dividing $2\Delta\,$. Fix a complete system $\{\q_1,\dots,\q_h\}$ of
representatives of different equivalence classes of the similarity class
of $\q\,$. Then
$$
\sum_{A\in\bigcup_{i=1}^h R^{*}(\q,p^2\q_i)/E(\q_i)} 1 \qquad =\qquad 
(1+\xa{\n})^{-1}\sum_{M\in\bigcup_{i=1}^H R(\n,p{\n}_i)/E({\n}_i)} 1
\qquad ,
$$
where $\n$ is integral quaternary quadratic form defined (with respect to
some basis) by the norm (5.10) on the even subalgebra $C_0(E)$ of the
Clifford Algebra $C(E,\q)$, $\xa{\n}=\left({{(-1)^2\Delta^2}\over{p}}\right)$
is the character (2.4) of $\n\,$, and $\{{\n}_1,\dots,{\n}_H\}$ is
a complete system of representatives of equivalent classes in the similarity
class of $\n\,$.}
\Proof
Indeed, taking ${\cal A}$ to be the zero matrix of order 4 in (6.4) and
using (6.7) we see that
$$
{\cal S}_{p^2}(\q,D_p,O_3)\ =
\ (1+\xa{\n})^{-1}\,{\cal S}_p(\n,{\cal D}_p, O_4)
$$
Because of (6.5) and (6.6), this relation between the isotropic sums
is equivalent to the above relation between numbers of automorphs, which
also proves (4.5).
\flushqed
Tracing back considerations of section 4, one can see that the above theorem
immediately implies
\Corollary{6.4}
{\sl With the notation and under the assumptions of Theorem 6.3, let in addition
$\q$ be positive definite. Then Shimura's lift of the generic 
theta-series $\th_{\{\q\}}(z)$ and the generic theta-series
$\th_{\{\n\}}(z)$ of the norm $\n$ on the even Clifford subalgebra 
$C_0(E,\q)$ have the same eigenvalues $\,p+1\,$ for all Hecke operators 
$T(p)$ with $p{\not|}\det\q$.}
\flushqed
Unfortunately up to the present the isotropic sums of type (6.3) or (6.6) have
not been computed in the case of prime $p$ dividing the determinant of
corresponding quadratic form. This prevents us from comparing directly
eigenvalues of $\th_{\{\n\}}(z)$ and of Shimura's lift of $\th_{\{\q\}}(z)$
on singular Hecke operators $T(p)\,$ with $p$ -- divisor of level. Thus, the
question if Shimura's lift of $\th_{\{\q\}}(z)$ for arbitrary ternary
form $\q$ is always proportional to $\th_{\{\n\}}(z)$ remains open.
\par
Nevertheless automorph class lifting $A \mapsto \Psi_A$ (5.30) constructed
in section 5 allows us to make progress in another direction and exhibit
algebraic origins of the equality (4.5) and therefore of Shimura's
correspondence for generic theta-series. Indeed, because of Eichler's
commutation relation (2.5), linear combinations of theta-series invariant
under Hecke operators (such as generic theta-series, for example)
are determined by corresponding sets of automorphs with multipliers $p$ or
$p^2\,$ (according to the weight of the theta-series). Thus, algebraic 
relations between various sets of automorphs can in principle be responsible
for correspondences between certain linear combinations of theta-series.
(Moreover, they can provide a natural generalization of such correspondences
to the case of indeterminate forms.) In particular, in the case of Shimura's
lift of generic theta-series of ternary forms we seek relations between two
sets of automorphs:
$$
\bigcup_{i=1}^h E(\q_i)\backslash R^{*}(\q_i,p^2\q) \qquad \hbox{and} \qquad
\bigcup_{j=1}^H E({\n}_j)\backslash R(\n_j,p\n) \qquad , 
$$
of ternary and quaternary quadratic forms respectively. Recall that by
factoring $\Psi_A$ in Lemmas 6.1 and 6.2, we extended our automorph lift
$A \mapsto \Psi_A$ to an inclusion $A \mapsto {\cal I}_A$ of
$\cup_i E(\q_i)\backslash R^{*}(\q_i,p^2\q)$ into
$\cup_j E(\n_j)\backslash R(\n_j,p\n)\,$.
Now we claim that this map is in fact an injection whose inverse can be
naturally extended turning the set of quaternary automorphs with multiplier
$p$ into a 2-fold covering of the set of ternary automorphs with multipliers
$p^2\,$. More precisely we have the following
\Theorem{6.5}
{\sl Let $\q$ be an integral nonsingular ternary quadratic form with 
matrix $Q$ given by (5.3) and let $N$ be the matrix (5.12). 
If $p$ is a prime number coprime to $\det\q\,$, then for any
${\cal M}\in\Lambda^4{\cal D}_p\Lambda^4\slash\Lambda^4$ such that
$\congr{N[{\cal M}]}{0}{p}$ there exists (a unique)
$A\in\Lambda^3 D_p\Lambda^3\slash\Lambda^3$ such that
$\congr{Q[A]}{0}{p^2}$ and ${\cal M}|\Psi_A$ (see (5.28)).}
\Proof
First of all, if ${\cal M}|\Psi_A$, then ${\cal M}^{-1}\Psi_A$ is an
integral automorph (of quaternary forms) with multiplier $p\,$, and thus
$p\Psi_A^{-1}{\cal M}$ is an integral matrix. Therefore, using the Remark 5.7
we see that
$$
{\cal M}|\Psi_A\Longleftrightarrow
\congr{\Psi_{p^2A^{-1}}\cdot{\cal M}}{0}{p}\quad .
\eqno (6.8)
$$
Next, let ${\cal M}={\cal UD}_p{\cal W}$ with ${\cal U,W}\in\Lambda^4\,$.
Since $\congr{N[{\cal M}]}{0}{p}$ if and only if
$\congr{N[{\cal MW}^{-1}]}{0}{p}\,$,
$\congr{\Psi_{p^2A^{-1}}{\cal M}}{0}{p}$ if and only if
$\congr{\Psi_{p^2A^{-1}}{\cal MW}^{-1}}{0}{p}\,$ and ${\cal M}$ is
defined only up to right multiplication by $\Lambda^4\,$, it is enough to
consider ${\cal W}=1_4$ and ${\cal M}={\cal UD}_p\,$.
\par
Similarly, let $A=VD_pW$ with $V,W\in\Lambda^3\,$. Clearly
$\congr{Q[A]}{0}{p^2}$ if and only if $\congr{Q[AW^{-1}]}{0}{p^2}\,$.
We also claim that $\congr{\Psi_{p^2A^{-1}}{\cal M}}{0}{p}$ if and only
if $\congr{\Psi_{p^2WA^{-1}}{\cal M}}{0}{p}\,$. Indeed, consider the
following commutative diagram:
$$\diagram{
\bigl(\Z^3,\q[p^{-1}AW^{-1}]\bigr) & \lharr{\alpha_W}{} &
\bigl(\Z^3,\q[p^{-1}A]\bigr) & \lharr{\alpha_{p^2A^{-1}}}{} &
\bigl(E,p^2\q\bigr) \cr
\dvarr{}{\iota} & & \dvarr{}{\iota} & & \dvarr{}{\iota} \cr
C\bigl(\Z^3,\q[p^{-1}AW^{-1}]\bigr) & \lharr{\varphi_W}{} &
C\bigl(\Z^3,\q[p^{-1}A]\bigr) & \lharr{\varphi_{p^2A^{-1}}}{} &
C\bigl(E,p^2\q\bigr) \cr
}$$
Using (5.20),(5.22) and (5.28) one can see that
$\Psi_{p^2WA^{-1}}=\Phi_W\cdot\Psi_{p^2A^{-1}}\in\Lambda^4\Psi_{p^2A^{-1}}$
since $\varphi_W$ is an isomorphism of Clifford algebras. Thus we can take
$W=1_3$ and $A=VD_p\,$. Next we set $V=(V_1,V_2,V_3)$ and
$$
{\cal U}=\bigl({\cal U}_1,{\cal U}_2,{\cal U}_3,{\cal U}_4\bigr)=
\pmatrix{
u_1 & u_2 & u_3 & u_4 \cr
\check U_1 & \check U_2 & \check U_3 & \check U_4 \cr} \quad .
$$
Thus given ${\cal M}={\cal UD}_p=
\bigl({\cal U}_1,{\cal U}_2,p\,{\cal U}_3,p\,{\cal U}_4\bigr)$ with
$\congr{\tr{\cal U}_i N{\cal U}_j}{0}{p}$ for $i=1,2$ we want to find
$A=VD_p=\bigl(V_1,pV_2,p^2V_3\bigr)$ such that $V\in\Lambda^3$ with
$\congr{Q[V_1]}{0}{p^2}\,,\,\congr{\tr V_2QV_1}{0}{p}$ and 
$$
\Psi_{p^2A^{-1}}\cdot{\cal M}=
\Psi_{(p^2D_p^{-1}V^{-1})}\cdot{\cal UD}_p=
\congr{
\pmatrix{
p & \tr Z \cr
0 & \tr V_1 \cr
0 & p\,\tr V_2 \cr
0 & p^2\,\tr V_3 \cr}\cdot
\pmatrix{
u_1 & u_2 & p\,u_3 & p\,u_4 \cr
\check U_1 & \check U_2 & p\,\check U_3 & p\,\check U_4 \cr}
}{0}{p}\,,
$$
where $Z=p^{-1}Z_{(p^2A^{-1})}$ depends on $V\,$, see (5.28). The last
condition is equivalent to the system of congruences
$\congr{\tr Z\check U_i}{0}{p}\,,\,\congr{\tr V_1\check U_i}{0}{p}$ for
$i=1,2\,$. Note that all of the above conditions on columns of matrix
$V\in\Lambda^3$ are modulo $p$ (or $p^2$) and therefore the choice of particular
representatives $\mathop{\rm mod}p^2$ is irrelevant as long as $V\in\Lambda^3\,$.
Recall that columns $V_1,V_2\in\Z^3$ can be complemented to a matrix
$(V_1,V_2,V_3)\in\Lambda^3$ if and only if the greatest common divisor of
principal minors of $(V_1,V_2)\in\Z_2^3$ is $1$. Using the above remarks
and Lemma 1.3 above we can replace the condition $V\in\Lambda^3$ by
linear independence of $V_1$ and $V_2$ modulo $p\,$. Thus, combining all
the conditions, we need to find $V_1,V_2\in\Z^3$ such that $V_1,V_2$ are
linearly independent modulo $p$ and
$$
\cases{
\congr{Q[V_1]}{0}{p^2} & \cr
\congr{\tr V_2 Q V_1}{0}{p} & \cr
\congr{\tr V_1\cdot\check U_i}{0}{p} & for $i=1,2$ \cr
\congr{\tr Z\cdot\check U_i}{0}{p} & for $i=1,2$ \cr
}\quad .$$
Denote
$$
{\cal W}=\tr{\cal U}^{-1}=
\bigl({\cal W}_1,{\cal W}_2,{\cal W}_3,{\cal W}_4\bigr)=
\pmatrix{
w_1 & w_2 & w_3 & w_4 \cr
\check W_1 & \check W_2 & \check W_3 & \check W_4 \cr}\in\Lambda^4\ ,
$$
then $\tr{\cal W}\cdot{\cal U}=
\bigl(\tr{\cal W}_i\cdot{\cal U}_j\bigr)=1_4$ which implies that
$\congr{w_i\cdot u_j+\tr\check W_i\cdot\check U_j}{0}{p}$ if $i\not=j\,$.
Next,since $\congr{N[{\cal UD}_p]}{0}{p}\,$, then using (5.14) we have
$$
\tilde N[\tr{\cal U}^{-1}]=
\Delta\cdot
\pmatrix{
2\,\det{Q_0} & \tr(Q_0 B) \cr
     &   \cr
(Q_0 B) & Q \cr}
\Biggl[
\pmatrix{
 & w_i & \cr
...& & ...\cr
 & \check W_i & \cr}
\Biggr] =
\pmatrix{
* & * & * & * \cr
* & * & * & * \cr
* & * & p & p \cr
* & * & p & p \cr}
$$
and therefore $\congr{\tr{\cal W}_i\tilde N{\cal W}_j}{0}{p}$
for $i,j=3,4\,$. Because ${\cal W}\in\Lambda^4\,$, the columns
${\cal W}_3,{\cal W}_4$ are linearly independent modulo $p$ and thus their
linear span (over $\F_p$) contains the column $\ncongr{\cal X}{0}{p}$
whose first component is zero modulo $p$:
$$
a{\cal W}_3+b{\cal W}_4=
\congr{\cal X}{\pmatrix{ 0 \cr \check X \cr}}{p}
\quad ,\quad
\ncongr{\check X}{0}{p}\quad .
$$
Using properties of ${\cal W}_3,{\cal W}_4$ exhibited above we see that
$\tr\check X\check U_i\,\equiv\,\congr{\tr{\cal XU}_i}{0}{p}$ for
$i=1,2$ and $Q[\check X]\,\equiv\,
\congr{\Delta^{-1}\tilde N[{\cal X}]}{0}{p}\,$. According to Lemma 2.3 of
[4] there exists $p^2$ distinct (modulo $p^2$) columns $X\in\Z^3$ such
that $\congr{X}{\check X}{p}$ and $\congr{Q[X]}{0}{p^2}\,$. For any of these
columns we can take $\congr{V_1}{X}{p^2}\,$.
\par
Next, consider the congruence $\congr{\tr V_1 QY}{0}{p}$ as a linear
equation with $3$ variables over $\F_p\,$. Since $Q\in GL_3(\F_p)$ and
$\ncongr{V_1}{0}{p}$ then $\ncongr{\tr V_1Q}{0}{p}$ and therefore the
space of solutions of the above equation is two dimensional. One 
nontrivial solution is $V_1\,$. Thus there exists $V_2\in\Z^3$ such that
$\congr{\tr V_1 Q V_2}{0}{p}$ and $V_1,V_2$ are linearly independent modulo
$p\,$. Using Lemma 1.3 we can change $V_1$ and $V_2$ modulo $p^2$ in
such a way that the matrix $(V_1,V_2)$ is primitive and therefore
complementable to an invertible matrix $V=(V_1,V_2,V_3)\in\Lambda^3\,$.
Thus given ${\cal U}\in\Lambda^4$ such that $\congr{N[{\cal UD}_p]}{0}{p}$
we can find $V\in\Lambda^3$ such that $\congr{Q[VD_p]}{0}{p^2}$ and
$\congr{\tr V_1\check U_i}{0}{p}$ for $i=1,2\,$. To prove the Theorem
we still need to show that $\congr{\tr Z\check U_i}{0}{p}$ for $i=1,2$
($Z$ depends on $V\,$, see above). Luckily, this condition is met 
automatically with our choice of $V\,$. Indeed, set $A=VD_p$ and
$\q'=p^{-2}\cdot\q[A]\,$, then $A\in R(\q,p^2\q')$ and so
$p^2\Psi_A^{-1}=\Psi_{p^2A^{-1}}\in R(\n',p^2\n)\,$, in other words
$N'[\Psi_{p^2A^{-1}}]=p^2N\,$ (see (5.12),(5.27) and Remark 5.7). Then
$N'[\Psi_{p^2A^{-1}}\cdot{\cal U}_i]=
\congr{p^2N[{\cal U}_i]}{0}{p^3}$ for $i=1,2\,$. But
$$
\Psi_{p^2A^{-1}}\cdot{\cal U}_i=
\congr{
\pmatrix{
p & \tr Z \cr
0 & \tr V_1 \cr
0 & p\tr V_2 \cr
0 & p^2\tr V_3 \cr}\cdot
\pmatrix{ u_i \cr \cr \check U_i \cr}
}{\pmatrix { \tr Z\check U_i \cr 0\cr 0\cr 0\cr}}{p}
\quad\hbox{ for }\ i=1,2
$$
by the above choice of $V\,$, therefore using (5.12) we have
$N'[\Psi_{p^2A^{-1}}\cdot{\cal U}_i]\equiv
\congr{2(\tr Z\check U_i)^2}{0}{p}$ and so $\congr{\tr Z\check U_i}{0}{p}$
for $i=1,2$ as stated (recall that $p\not=2$ since $2|\det\q\,$).
\par
Summing up, there exists $A$ such that
$\congr{\Psi_{p^2A^{-1}}\cdot{\cal UD}_p}{0}{p}\,$, i.e. ${\cal M}|\Psi_A\,$.
Finally, it is easy to see that for a given
${\cal M}\in\cup_i R(\n,p\n_i)\slash E(\n_i)$ such an
$A\in\cup_j R(\q,p^2\q_j)\slash E(\q_j)$ is unique. Indeed, first note that
since $A\in\Lambda^3 D_p\Lambda^3$ then the double coset 
$\Lambda^4\Psi_A\Lambda^4$ contains matrix of the form
$$
\pmatrix{
p & * & * & * \cr
0 &   &   &   \cr
0 &  & D_p &  \cr
0 &   &   &   \cr
}$$
and so $\mathop{\rm rank}_{\F_p}\Psi_A$ is either $1$ or $2\,$. Then (6.4)
implies that each $\Psi_A\in R(\n,p^2\n')$ is divisible from the left by
at most $2=1+\xa{\n}$ different cosets 
${\cal M}\in\cup_{i=1}^H R(\n,p\n_i)\slash E(\n_i)\,$.
(The same is true for any $A\in\cup_i R^{*}(\q,p^2\q_i)\,$). If the same
${\cal M}$ would divide $\Psi_A$ for different classes of $A$'s, then
the union over $A$
$$
\bigcup\limits_{A\in\cup_{i=1}^h R^*(\q,p^2\q_i)\slash E(\q_i)}
\{\ {\cal M}\in\cup_{j=1}^H R(\n,p\n_j)\slash E(\n_j)\ ;\ {\cal M}|\Psi_A\ \}
\eqno (6.9)
$$
would not be disjoint and therefore its cardinality would be less than
$$
(1+\xa{\n})\cdot\sum_{i=1}^h r^{*}(\q,p^2\q_i)\slash e(\q_i) \qquad =
\qquad (1+\xa{\n})\cdot(p+1)\ .
$$
On the other hand, we have just proved that any 
${\cal M}\in\cup_j R(\n,p\n_j)\slash E(\n_j)$ divides from the left some
$\Psi_A$ and thus belongs to the union (6.9), whose cardinality in view of
this should be at least $(1+\xa{\n})\cdot(p+1)\,$, see (6.4) and (6.5).
We conclude that for any ${\cal M}$ there exists exactly one 
$A\in\Lambda^3 D_p\Lambda^3\slash\Lambda^3$ such that
$\congr{Q[A]}{0}{p^2}$ and ${\cal M}|\Psi_A\,$.
Incidentally we have also proved that
$\mathop{\rm rank}_{\F_p}{\Psi_A}=1$ for any such $A\,$.
\flushqed
\Remark{6.6}
The above result means that for any class
${\cal M}E(\n_i)\in R(\n,p\n_i)\slash E(\n_i)$ there exists a unique class
$AE(\q_j)\in R^*(\q,p^2\q_j)\slash E(\q_j)$ such that any
${\cal M'}\in {\cal M}E(\n_i)$ divides from the left $\Psi_{A'}$ for any
$A'\in AE(\q_j)\,$.
Thus the established correspondence is indeed between classes modulo groups
of units rather than between individual automorphs. (To see this let
$A'=AE$ and ${\cal M'}={\cal ME}$ with $E\in E(\q_j),\ {\cal E}\in E(\n_i)\,$.
Then using (5.28), (6.8) and Lemma 5.3 we have: 
$
{\cal M'}|\Psi_{A'}\Leftrightarrow
\congr{\Psi_{p^2{A'}^{-1}}\cdot{\cal ME}}{0}{p}\Leftrightarrow
\congr{\Phi_{E^{-1}}\Psi_{p^2{A}^{-1}}\cdot{\cal M}}{0}{p}\Leftrightarrow
{\cal M}|\Psi_{A}
$
since $\Phi_{E^{-1}}\in E(\n_j)$ and so it is invertible.
\par
\bigskip
In order to obtain a generalization of Shimura's lift in terms of automorph
class theory, we need to consider correspondences between left (rather than
right) classes of automorphs. To this end we extend the automorph class
lift $\Psi$ (5.28) and define the map 
 $\Upsilon: \cup_i E(\q_i)\backslash R^*(\q_i,p^2\q)\rightarrow
\cup_j E(\n_j)\backslash R(\n_j,p\n)$ via the following chain
elements of which are already familiar (see Lemmas 5.5 and 6.2):
$$
\Upsilon\,:\ A\mapsto\Psi_A\mapsto p^2(\Psi_A)^{-1}=\Psi_{p^2A^{-1}}
\mapsto {\cal I}_{p^2A^{-1}}\mapsto 
p({\cal I}_{p^2A^{-1}})^{-1}=\Upsilon_A
\eqno (6.10)
$$
Summing up our knowledge about the above maps we arrive to the following
\Theorem{6.7}
{\sl 
Let $\q$ be an integral nonsingular ternary quadratic form and let
$\n$ be the integral quadratic form defined (up to integral equivalence)
by the norm on the even subalgebra of the Clifford algebra $C(\q)\,$.
Take $\{\q_1,\dots ,\q_h\}$ and $\{\n_1,\dots ,\n_H\}$ to be
complete systems of representatives of different equivalence classes
of the similarity classes of $\q$ and of $\n$ respectively. Let $p$ be a
prime number coprime to $\det\q\,$. Then map of classes of automorphs
$A\mapsto\Upsilon_A$ defined via (6.10), (5.28) and Theorem 6.5:
$$
\Upsilon\,:\ \bigcup\limits_{i=1}^h E(\q_i)\backslash R^*(\q_i,p^2\q)
\lhook\joinrel\longrightarrow 
\bigcup\limits_{j=1}^H E(\n_j)\backslash R(\n_j,p\n)
$$
is an injection such that $\Upsilon_A$ divides $\Psi_A$ from the right for
any $A\in\cup_i E(\q_i)\backslash R^*(\q_i,p^2\q)\,$.
}
\Proof
Indeed, the map 
$\Psi: E(\q_i)\backslash R^*(\q_i,p^2\q)\rightarrow 
E(\n_i)\backslash R^*(\n_i,p^2\n)$ is an injection according to (5.30) and
Lemmas 5.3 , 5.5. The second element of our chain is the bijection
${\cal M}\mapsto p^2{\cal M}^{-1}$ of classes of primitive automorphs
$E(\n_i)\backslash R^*(\n_i,p^2\n)\rightarrow R^*(\n,p^2\n_i)\slash E(\n_i)$
restricted to the image of $\Psi\,$. The equality 
$p^2(\Psi_A)^{-1}=\Psi_{p^2A^{-1}}$ was established in Remark 5.7.
The map $\Psi_{p^2A ^{-1}}\mapsto{\cal I}_{p^2A^{-1}}$ , 
$R^*(\n,p^2\n_i)\slash E(\n_i)\rightarrow R(\n,p\n_j)\slash E(\n_j)$ 
for some $j$ is defined via extension of ideals of corresponding Clifford
algebras in Lemmas 6.1 and 6.2. It provides us with ${\cal I}_{p^2A^{-1}}\in
\Lambda^4{\cal D}_p\Lambda^4\slash\Lambda^4\,$ such that ${\cal I}_{p^2A^{-1}}$
divides $\Psi_{p^2A^{-1}}$ from the left and therefore is an injection 
according to the Theorem 6.5 and Remark 6.6. 
The last step in the chain (6.10) is the 
bijection ${\cal M}\mapsto p{\cal M}^{-1}$ , $R(\n,p\n_j)\slash E(\n_j)
\rightarrow E(\n_j)\backslash R(\n_j,p\n)\,$. Thus, being a composition of
injections, $\Upsilon$ is certainly an injection itself. Moreover, since
${\cal I}_{p^2A^{-1}}$ divides $p^2(\Psi_A)^{-1}$ from the left, then
$\Upsilon_A=p({\cal I}_{p^2A^{-1}})^{-1}$ divides $\Psi_A$ from the right.
\flushqed
From the above theorem it follows easily that one can consider the set of 
classes of quaternary automorphs of $\n$ as a 2-fold
covering of the set of classes of ternary automorphs of $\q\,$:
\Corollary{6.8}
{\sl With the notation and under the assumptions of Theorem 6.7,}
$$
\displaylines{\quad
\bigcup_{j=1}^H E(\n_j)\backslash R(\n_j,p\n) \ =\hfill\cr 
\hfill\bigcup\limits_{A\in\cup_{i=1}^h E(\q_i)\backslash R^*(\q_i,p^2\q) }
\{\ {\cal M}\in\cup_j E(\n_j)\backslash R(\n_j,p\n) \ ;\ 
\Psi_A{\cal M}^{-1}\in\Z^4_4\}\ ,
\qquad (6.11)\cr
}$$
{\sl where for a fixed $A$ the set
$\{ {\cal M}\in\cup_j E(\n_j)\backslash 
R(\n_j,p\n) \ ;\ \Psi_A{\cal M}^{-1}\in\Z^4_4\}$
consists of exactly $2$ classes and the union over $A$ is disjoint.}
\Proof
According to the Remark 6.6, for each 
${\cal M}\in\cup_j E(\n_j)\backslash R(\n_j,p\n)$ there exists a unique
$A\in\cup_i E(\q_i)\backslash R^*(\q_i,p^2\q)$ such that
$(p{\cal M})^{-1}$ divides $\Psi_{p^2A^{-1}}$ from the left, which 
immediately implies that ${\cal M}$ divides 
$\Psi_A=p^2(\Psi_{p^2A^{-1}})^{-1}$ from the right. This establishes the
proposed equality (6.11) and shows that cardinality of the union over $A$
in the right hand side of (6.11) is equal to $(1+\xa{\n})\cdot(p+1)\,$, see
formulas (6.3),(6.4) and (6.5). On the other hand, because of (6.6) and (6.7)
we see that this cardinality can be achieved only if the union over $A$ in
(6.11) is disjoint (recall that in the proof of Theorem 6.5 we already
have established that $\mathop{\rm rank}_{\F_p}\Psi_A=1$ and so
the number of elements in the set
$\{ {\cal M}\in\cup_j E(\n_j)\backslash 
R(\n_j,p\n) \ ;\ \Psi_A{\cal M}^{-1}\in\Z^4_4\}$
is equal to $1+\xa{\n}=2$ according
to (6.4)).
\flushqed
In other words, Corollary 6.8 states that each class in
$\cup_j E(\n_j)\backslash R(\n_j,p\n)$ is a right divisor of a unique class in
$\Psi\left(\cup_i E(\q_i)\backslash R^*(\q_i,p^2\q)\right)\,$. And conversely,
each class in the image
$\Psi\left(\cup_i E(\q_i)\backslash R^*(\q_i,p^2\q)\right)$ is
divisible from the right by exactly $(1+\xa{\n})=2$ classes of automorphs in
$\cup_j E(\n_j)\backslash R(\n_j,p\n)\,$.
\vfill\break
%
%
\bigskip\noindent
{\bf 7. Concluding remarks}
\smallskip\noindent
Automorph class lift and its factorizations studied in preceding
sections sheds a new light
on relations between representations of integers by ternary and
quaternary quadratic forms. In particular case of positive definite
ternary forms one can deduce from its construction relations between
corresponding theta-series similar to Shimura's lift as defined by
P.Ponomarev in [10]. Still the theory is far from being complete
and we would like to discuss some open questions.
\par
First of all it is necessary to address the issue of singular primes,
i.e. those $p$ which divide the determinant $\det\q\,$. As we already
noted in Corollary 6.4, the isotropic sums of types (6.3) and (6.6)
have not been computed in the case of such primes. This prevents us from 
stating a stronger versions of Theorem 6.3 and Corollary 6.4 which would
satisfactorily resolve the question of whether or not Shimura's lift of
$\Theta_{\{\q\}}(z)$ is always equal to a normalized $\Theta_{\{\n\}}(z)$.
Computation of singular isotropic sums would also have other important
applications, such as a complete Euler product decomposition of general Eichler
and Epstein zeta-functions as well as investigation of their analytical
properties (see [2]).
\par
Another unclear issue concerning singular primes is proper construction of
the automorph class lift (5.28)
$ \Psi:\,\cup_i E(\q_i)\backslash R^*(\q_i,p^2\q)\lhook\joinrel\rightarrow
\cup_j E(\n_j)\backslash R^*(\n_j,p^2\n) $
for such $p$ and its further factorization into a product of quaternary
automorphs with prime multipliers (similar to results of Lemmas 6.1 and 6.2).
All our proofs use the primality of $p$ to $\det\q$ and it
would be very interesting to see to what extent the same construction works
in the singular case and what modifications (if any) should be made.
\par
The next interesting possibility is connected to investigation of effects of
the automorph class lift on individual quadratic forms. (In case of
positive definite quadratic forms this could lead to correspondences
between certain spaces of modular forms spanned by theta-series of different
weights which are invariant under Hecke operators.) For this one needs
to describe explicitly factorizations of an individual image
$\Psi\left(E(\q_i)\backslash R^*(\q_i,p^2\q)\right)$
into products of quaternary automorphs with multiplier $p\,$.
(Recall that we provided such a description only for the whole
$\Psi\left(\cup_i E(\q_i)\backslash R^*(\q_i,p^2\q)\right)$ which 
corresponds to Shimura's lift of the generic theta-series
$\Theta_{\{\q\}}(z)\,$). Inspired by results of P.Ponomarev [10], here
is a sketch of what one might expect. 
\par
The subspace of modular forms of half-integral weight $3\slash 2$ spanned by
$h$ theta-series of ternary quadratic forms (4.2) is invariant
under Hecke operators $|_{3\over 2}\, T(p^2)$ and its eigenmatrix with respect
to $\Theta(z,\q_i)\slash e(\q_i)\,,\ 1\leq i\leq h\,$ is $\eich{\q}{2}$
given by (2.6). Let $F_i(z)$ denote Shimura's lift $\Theta^{(a)}(z,\q_i)$
for some fixed $a$
and let $F_{\q}(z)=\tr\left(\dots ,F_i(z),\dots\right)$ be the associated
vector-valued modular form of integral weight $2$. Then we expect the space
spanned by $F_i(z),\ 1\leq i\leq h$ to be invariant under Hecke operators
$|_{2}\,T(p)$ whose eigenmatrix with respect to $F_i$'s is again given by
$\eich{\q}{2}$. Furthermore, let
$\Theta_{\n}(z)=\tr\left(\dots ,\Theta(z,\n_j)\slash e(\n_j),\dots\right)$
be the the vector-valued theta-series of weight 2 associated to a complete
system $\{\n_1,\dots,\n_H\}$ of quaternary quadratic forms representing 
different equivalent classes of the similarity class of quadratic form
defined by the norm on the even subalgebra of Clifford algebra
$C(\Z^3,\q_i)$ for some $i\,$. (We note again that we can take $\n_i$
to be the norm on $C_0(\Z^3,\q_i),\ 1\leq i\leq h$ and so $h\leq H$, see
Lemmas 5.4 and 5.5). We expect the space spanned by Shimura's lifts of
$\Theta(z,\q_i),\ 1\leq i\leq h$ to be a subspace of the space spanned by
the normalized theta-series $\Theta(z,\n_j)\slash e(\n_j),\ 1\leq j\leq H\,$,
which is invariant under Hecke operators $|_{2}\,T(p)$ . In other words,
we hope that there exists a matrix $X\in\Q^h_H$ (of arithmetical nature) such
that $F_{\q}(z)=X\cdot\Theta_{\n}(z)\,$. This would imply that
$$
\eich{\q}{2}\cdot X\Theta_{\n}(z)=F_{\q}(z)\bigm|_{2}T(p)=X\Theta_{\n}(z)
\bigm|_{2}T(p)=(1+\xa{\n})^{-1}X\,\eich{\n}{}\cdot\Theta_{\n}(z)
$$
Building up expectations, one hopes that
$\eich{\q}{2}\cdot X=(1+\xa{\n})^{-1}X\cdot\eich{\n}{}$ for some arithmetical
matrix $X\in\Q^h_H\,$, in other words
$$
\Bigl(|\Repr{2}ij|\Bigr)_{h\times h}\cdot X=
(1+\xa{\n})^{-1}X\cdot\Bigl(|\Reprn{}kl|\Bigr)_{H\times H}
\eqno (7.1)
$$
One can see that such an $X$ is not totally impossible using Theorem 6.3
and taking for example 
$h=1$ and $X=(\sum_j r(\n_j,1)\slash e(\n_j))^{-1}\cdot(1,\dots,1)\,$, the
particular scalar multiple is chosen in order to make the first Fourier
coefficient of $X\cdot\Theta_{\n}(z)$ equal to 1. Another, rather trivial
example of matrix $X$ with property (7.1) is
$$
X=\pmatrix{
e^{-1}(\q_1) & \dots & e^{-1}(\q_1) \cr
\vdots &  & \vdots \cr
e^{-1}(\q_h) & \dots & e^{-1}(\q_h) \cr
}_{h\times H}\qquad ,
$$
for which we have 
$\eich{\q}{2}\cdot X=(p+1)\cdot X=(1+\xa{\n})^{-1}X\cdot\eich{\n}{}$
according to the Theorem 6.3. Both examples produce the generic 
theta-series $\Theta_{\{\n\}}(z)$ as a Hecke eigenform of weight 2 
whose eigenvalue is equal to $\eich{\q}{2}\,$. It would be extremely
interesting to find nontrivial (of $\mathop{\rm rank}>1$) examples of
matrices with property (7.1), which could lead to construction of other
subspaces invariant under Hecke operators of the space spanned by the
theta-series $\Theta(z,\n_j),\ 1\leq j\leq H$ with eigenmatrix 
$\eich{\q}{2}\,$. Furthermore, and perhaps more interesting, equality
of type (7.1) would provide a direct link between various zeta-functions
associated to ternary and quaternary quadratic forms, which could
lead to insights concerning their analytical properties
(such as meromorphic continuation and functional equation). Let us
illustrate this with a simple example.
\par
If we start with $\q(x)=x_1^2+x_2^2+x_3^2$ - the sum of $3$ squares,
then corresponding quaternary form (5.12) on $C_0(E,\q)$ is the the
sum of $4$ squares $\n(y)=y_1^2+y_2^2+y_3^2+y_4^2\,$. Both forms are
positive definite and both have the class number equal to $1$ (i.e.
$h=1$ and $H=1\,$), so their theta-series
$\Theta(z,\q)=e(\q)\cdot\Theta_{\{\q\}}(z)$ and 
$\Theta(z,\n)=e(\n)\cdot\Theta_{\{\n\}}(z)$ are eigenforms of Hecke operators
$|_{3\over 2}\,T(p^2)$ or $|_2\,T(p)$ with respective eigenvalues
$\eich{3}{2}=r^*(\q,p^2\q)\slash e(\q)$ or
$\eich{4}{}=r(\n,p\n)\slash e(\n)\,$, see (2.6).
The only singular prime dividing the determinants in our example is
$p=2\,$. By the Corollary 6.4 odd-numbered Fourier coefficients of
Shimura's lift of $\Theta(z,\q)$ coincide with those of 
${1\over 8}\Theta(z,\n)\,$,we normalize with r(\n,1)=8.
Note that according to the Theorem 6.3 the equality (7.1) actually holds
in this situation with $X=1\,$: $\eich{3}{2}={1\over 2}\eich{4}{}$ 
for $p\not=2\,$.
Next we set
$$
\rey{k}{m}=
\biggl|R\left({\textstyle{\sum\limits_{i=1}^k x_i^2\,,\ m}}\right)\biggr|
\Biggm\slash
\biggl|E\left({\textstyle{\sum\limits_{i=1}^k x_i^2}}\right)\biggr|\quad ,
$$
then, using (2.8) and (2.9) we have the following chain of equalities
for any square-free positive integer $a$ and any positive integer $b\,$:
$$\openup2pt \displaylines{
\qquad
\sum_{m\hbox{\sevenrm -odd}}{{\rey{4}{mb}}\over{m^s}}\ =\hfill\cr
\prod_{p\hbox{\sevenrm -odd}}
\left(1-{\scriptstyle {1\over 2}}\eich{4}{}p^{-s}+p^{1-2s}\right)^{-1}
\cdot\rey{4}{b}\ =\ 
\prod_{p\hbox{\sevenrm -odd}}
\left(1-\eich{3}{2}p^{-s}+p^{1-2s}\right)^{-1}\cdot\rey{4}{b}\ =\cr
\sum_{m\hbox{\sevenrm -odd}}{{\chi_{\q}(m)({{2a}\over m})}\over{m^s}}\cdot
\prod_{p\hbox{\sevenrm -odd}}
\left(1-\xa{\q}\Bigl({{2a}\over p}\Bigr)p^{-s}\right)\cdot
\prod_{p\hbox{\sevenrm -odd}}
\left(1-\eich{3}{2}p^{-s}+p^{1-2s}\right)^{-1}
\cdot\rey{4}{b}\ =\cr
\hfill
L\left(s\,,({{-a}\over{\cdot}})\right)\cdot
\sum_{m\hbox{\sevenrm -odd}}{{\rey{3}{m^2a}}\over{m^s}}
\ \cdot{{\rey{4}{b}}\over{\rey{3}{a}}}
\qquad .\cr
}$$
Taking $b=1$ and setting $L\left(s,({{-a}\over{\cdot}})\right)=L_a(s)\,$,
where $({{-a}\over{\cdot}})$ is the Legendre symbol, we conclude that
$$
\sum_{m\hbox{\sevenrm -odd}}{{\rey{3}{m^2a}}\over{m^s}}\ =\ 
L^{-1}_{a}(s)\cdot\sum_{m\hbox{\sevenrm -odd}}{{\rey{4}{m}}\over{m^s}}
\ \cdot\,\rey{3}{a}
\qquad .
\eqno (7.2)
$$
The particular zeta-function on the left-hand side of (7.2) has Euler product
(2.9) and is associated to theta-series $\Theta(z,\,x_1^2+x_2^2+x_3^2)$ of 
half-integral weight. Now one can also deduce analytical properties
(meromorphic continuation and functional equation) of this zeta-function
from known properties of the Epstein zeta-function and the $L$-series on
the right-hand side of (7.2).
\par
Finally one can try to give an interpretation (and generalization) of
equality (7.1) in terms of the automorph class theory which would
provide a direct link between $\Eich{\q}{2}$ and $\Eich{\n}{}\,$, see (3.1),
i.e. between individual classes of ternary automorphs
$E(\q_i)\backslash R^*(\q_i,p^2\q_j)$ and quaternary automorphs
$E(\n_k)\backslash R^*(\n_k,p\n_l)$ with fixed $i,j,k,l\,$,
for an arbitrary integral nonsingular ternary quadratic form $\q\,$.
This would also automatically furnish a natural proof for the above
conjectures on relations between numbers of representations and, as
a consequence, between corresponding Eichler or Epstein zeta-functions.
Such an interpretation seems to require a better understanding of
arithmetic of even Clifford subalgebras $C_0(E,\q)$ in order to
characterize explicitly factorizations of individual lifts 
$\Psi_A\subset E(\n_i)\backslash R^*(\n_i,p^2\n_j)$ in terms of
automorph classes in $E(\n_k)\backslash R^*(\n_k,p\n_l)$ for
specific $k$ and $l\,$. To this day attempts to find such an 
interpretation were unsuccessful.
\par
At last we mention the most intriguing mystery concerning Shimura's
lift for theta-series -- the question of theta-series of quadratic 
forms in more than 3 variables. The apparatus of Clifford algebras 
extensively used in the present paper seems to be of no use in that 
situation since on one hand dimensions of corresponding algebras 
do not match the expected weights of the lifting and on the other
hand the standard norm on an even Clifford subalgebra of dimension
greater than $4$ does not define a quadratic form in general.
To author's knowledge there is no result describing effects of
Shimura's lift on such theta-series. The mountain is still in clouds!
\vfill\break
%
%
\bigskip\noindent
{\bf REFERENCES}
\medskip

\item{[1]} {\sc A. Andrianov}, {\it Automorphic factorizations of solutions
of quadratic congruences and their applications}, St.Petersburg Math. J.
{\bf 5} (1994), no. 5, 1--38.

\item{[2]} {\sc --- }, {\it Quadratic congruences and rationality of local
zeta-series of ternary and quaternary quadratic forms}, St.Petersburg Math.
J. {\bf 6} (1995), no 2, 199--240.

\item{[3]} {\sc --- }, {\it  Quadratic Forms and Hecke Operators},
Grundlehren Math. Wiss. {\bf 286}, Springer-Verlag, Berlin--New York, 1987.

\item{[4]} {\sc F. Andrianov}, {\it Multiplicative decompositions of integral
representations by quadratic forms}, Zap. Nauchn. Sem. St.Petersburg. Otdel.
Mat. Inst. Steklov (POMI) {\bf 212} (1994), 12--55.

\item{[5]} {\sc --- }, {\it Euler expansions of formal zeta-series of
quadratic forms in an odd number of variables}, St.Petersburg Math. J. 
{\bf 7} (1996), no. 1, 1--16.

\item{[6]} {\sc M. Eichler}, {\it Quadratishe Formen und orthogonale Gruppen},
Grundlehren Math. Wiss. {\bf 63}, Springer-Verlag, Berlin--New York, 1974.

\item{[7]} {\sc E. Freitag}, {\it Siegelsche Modulfunktionen},
Grundlehren Math. Wiss. {\bf 254}, Springer-Verlag, Berlin--New York, 1983.

\item{[8]} {\sc H. Hida}, {\it Elementary theory of L-functions and 
Eisenstein series}, London Math. Soc. Student Texts {\bf 26}, Cambridge
University Press, Cambridge, 1993.

\item{[9]} {\sc M. Kneser}, {\it Quadratishe Formen, Vorlesung},
Mathematisches Institut der Universit\"at G\"otingen, G\"otingen, 1992.

\item{[10]} {\sc P. Ponomarev}, {\it Ternary Quadratic Forms and Shimura's
correspondence}, Nagoya Math. J. {\bf 81} (1981), 123--151.

\item{[11]} {\sc G. Shimura}, {\it Introduction to the Arithmetical Theory
of Automorphic Functions}, Iwanami Publishers and Princeton University Press,
Tokyo--Princeton, N.J., 1971.

\item{[12]} {\sc G. Shimura}, {\it On modular forms of half integral weight},
Ann. of Math. {\bf 97} (1973), 440--481.
%
%
\bye